\documentclass{article}
\usepackage{arxiv}
\usepackage[utf8]{inputenc}
\usepackage{amsfonts}
\usepackage{amssymb}
\usepackage{amsthm}
\usepackage{amsmath}
\usepackage{graphicx,subfigure}
\usepackage{empheq}
\usepackage{indentfirst}
\usepackage{tikz}
\usepackage{mathrsfs}
\usepackage{lineno}
\usepackage{float}
\usepackage{caption}
\usepackage{tabularx}
\usepackage{booktabs}
\usepackage{multirow}
\usepackage{booktabs}
\usepackage[colorlinks=true]{hyperref}
\usepackage{color}
\usepackage[numbers]{natbib}
\newtheorem{theorem}{\textbf{Theorem}}[section]
\newtheorem{lemma}{\textbf{Lemma}}[section]
\newtheorem{proposition}{\textbf{Proposition}}[section]
\newtheorem{corollary}{\textbf{Corollary}}[section]
\newtheorem{remark}{\textbf{Remark}}[section]
\newtheorem{definition}{\textbf{Definition}}[section]

\allowdisplaybreaks[4]

\def\be{\begin{equation}}
\def\ee{\end{equation}}
\def\bea{\begin{eqnarray}}
\def\eea{\end{eqnarray}}
\def\bt{\begin{theorem}}
\def\et{\end{theorem}}
\def\bl{\begin{lemma}}
\def\el{\end{lemma}}
\def\br{\begin{remark}}
\def\er{\end{remark}}
\def\bp{\begin{proposition}}
\def\ep{\end{proposition}}
\def\bc{\begin{corollary}}
\def\ec{\end{corollary}}
\def\bd{\begin{definition}}
\def\ed{\end{definition}}

\usepackage{bm}
\def\Au{\bm{u}}
\def\bfu{\bm{u}}
\def\bfv{\bm{v}}
\def\baru{\overline{\bm{u}}}
\def\tu{\widetilde{\bm{u}}}
\def\tphi{\widetilde{\phi}}
\def\12{\frac{1}{2}}

\def\sw{{We}^\ast}

\def\tp{\tilde{\phi}^{n+\frac{1}{2}}}
\def\ap{{\phi}^{n+\frac{1}{2}}}

\def\bdx{\mbox{d}\bm{x}}

\title{A high-order accurate unconditionally stable bound-preserving numerical scheme for the Cahn-Hilliard-Navier-Stokes equations}

\author{Yali Gao\thanks{School of Mathematics and Statistics, Northwestern Polytechnical University, Xi'an, Shanxi, China. Email:{\texttt{gaoylimath@nwpu.edu.cn}}} \and \textbf{Daozhi Han}\thanks{Department of Mathematics, State University of New York at Buffalo, Amherst, NY, USA. Email: \texttt{daozhiha@buffalo.edu}} \and \textbf{Sayantan Sarkar}\thanks{Department of Mathematics, State University of New York at Buffalo, Amherst, NY, USA. Email: \texttt{sayantan@buffalo.edu}}} 
\date{}
\begin{document}
\maketitle
\begin{abstract}
A high-order numerical method is developed for solving the Cahn-Hilliard-Navier-Stokes equations with the Flory-Huggins potential. The scheme is based on the $Q_k$ finite element with mass lumping on rectangular grids, the second-order convex splitting method, and the pressure correction method. The unique solvability, unconditional stability, and bound-preserving properties are rigorously established. The key to bound-preservation is the discrete $L^1$ estimate of the singular potential. Ample numerical experiments are performed to validate the desired properties of the proposed numerical scheme.
\end{abstract}

\keywords {Cahn-Hilliard-Navier-Stokes,  Flory-Huggin potential,  high order accuracy, bound-preserving, unique solvability, quadrilateral element}

\section{Introduction}
{In this article, we address the solution of the Cahn-Hilliard-Navier-Stokes (CHNS) equations for binary, incompressible, and macroscopically immiscible fluids in a rectangular domain $\Omega \subseteq \mathbb{R}^d, d=2, 3$. The non-dimensional CHNS equations read as follows}
\begin{align}
&\frac{\partial \phi}{\partial t}+ \nabla \cdot (\phi \Au)=\frac{1}{Pe}\nabla \cdot(M(\phi) \nabla \mu), \quad  \text{ in } \Omega_T \label{CH} \\
&\mu=f(\phi)-\epsilon^2 \Delta \phi,  \quad \text{ in } \Omega_T \label{CP} \\
&\frac{\partial \Au}{\partial t}-\frac{1}{Re} \nabla \cdot\big(\eta(\phi) \nabla \Au\big)+\Au\cdot \nabla \Au+\nabla p=-\frac{\epsilon^{-1}}{\sw} \phi \nabla \mu,  \quad  \text{in } \Omega_T  \label{NS}\\
&\nabla \cdot \Au=0, \quad  \text{in } \Omega_T \label{div} 
\end{align}
where {$\Au$, $p$, $\phi$ and $\mu$ denote the velocity, pressure, phase field variable and  the chemical potential, respectively;  $f(\phi)=F^\prime(\phi)$ with the  Flory-Huggins free energy density  $F(\phi)$  given by  }
\begin{align}\label{poten}
	 F(\phi)=\frac{1}{2} [(1+\phi) \ln{(1+\phi)}+(1-\phi) \ln{(1-\phi)}]-\frac{\theta_c}{2}\phi^2,
\end{align}
 and $\theta_c>1$ is the relative critical temperature.
Let $\Omega_T:=\Omega \times (0,T)$ with $T>0$ be  a fixed constant. Here $\eta(\phi)=\frac{1-\phi}{2} +\frac{1+\phi}{2}\frac{\eta_2}{\eta_1}$ is the dimensionless viscosity with $\eta_1, \eta_2$ the kinematic viscosity of fluid 1 and 2, respectively;
$Re$ is the Reynolds number;  $\sw$ is the modified Weber number that measures the relative strengths of the kinetic and surface energies; {$\epsilon$ is a positive parameter proportional to the diffuse interface width;} $Pe$ is the Peclet number; mobility $M(\phi)$ is a  function of phase variable  such that  $0<m_1 \leq M(\phi) \leq m_2$. 

We assume that  the system is equipped with the following initial and boundary conditions
\begin{align}\label{BCs}
 \Au&=0, \quad \text{on } \partial \Omega \times (0,T) \\
 \nabla \phi \cdot \bm{n}=\nabla \mu \cdot \bm{n} &=0, \quad \text{on } \partial \Omega \times (0,T) \\
 (\Au, \phi)|_{t=0}&=(\Au_0, \phi_0), \quad \text{in } \Omega.
\end{align}
Here $\bm{n}$ denotes the unit outer normal vector of the boundary $\partial \Omega$. 
{It is clear that the CHNS system  satisfies  the mass conservation law } 
\begin{align}\label{Mass-c}
\frac{d}{dt}\int_\Omega \phi\, \bdx=0,
\end{align}
{and  the energy dissipation law}
\begin{align}\label{ConEnL}
\frac{d}{dt}E(\Au, \phi)=-\frac{1}{Re}\int_\Omega  |\nabla \Au|^2\,\bdx -\frac{\epsilon^{-1}}{\sw} \int_\Omega M(\phi)|\nabla \mu|^2\, \bdx,
\end{align}
where the  total energy $E$ is defined as
\begin{align}\label{Etot}
E(\Au, \phi)=\int_{\Omega}\frac{1}{2}|\Au|^2\, \bdx+ \frac{1}{\sw}\int_{\Omega} \big( \frac{1}{\epsilon} F(\phi)+\frac{\epsilon}{2}|\nabla \phi|^2\big)\, \bdx.
\end{align}


The CHNS model, also known as Model H \cite{HoHa1977, GPV1996, LiSh2003, LoTr1998}, is a diffuse interface model that describes the evolution of diffusive interface layers of finite thickness separating the two fluids. The moving and deforming layered structure of the solution—a narrow region where the solution undergoes steep changes (large gradient)—demands sufficiently fine resolution within the diffused interface to minimize spurious oscillation, which otherwise may pollute the solution and cause a blow-up of the computer code. For efficiency purposes, adaptive mesh refinement is generally adopted in practice, and unconditionally stable time-marching algorithms are preferred to avoid prohibitively small time-step constraints. Besides, it is also desirable that the numerical schemes obey a discrete version of the energy law to accurately capture the dynamics in long-term simulations. Many strategies have been proposed in the past decade for the design of energy-law-preserving (unconditionally stable) algorithms for solving the CHNS system with a polynomial potential and related phase-field models, cf. \cite{KKL2004, ShYa2010b, HGW2011, FeWi2012, GLL2014, HaWa2015, ShYa2015, Han2016, ZYGW2017, GLLW2017, GZYW2018, SXY2019, Fu2020, HaJi2020, ChZh2020, CHWZ2020, ZhHa2021}. Error analysis of these methods can be found in \cite{KSW2008, DFW2015, CLWW2016, DWWW2017, LCWW2017, CaSh2018, LiSh2020, CHWWWZ2022, CSWY2023, LMR2023}. Note that for the CHNS system with a polynomial potential, the phase field variable $\phi$ may be out of bounds between $-1$ and $1$ unless a degenerate mobility is utilized. Nonetheless, flux and slope limiting are utilized in \cite{LRTLR2022} for post-processing either linear or quadratic approximation so that there is no overshoot or undershoot in the numerical solution. A high-order limiter is recently constructed in \cite{LRSZ2023} for solving the CHNS system.

When the Flory-Huggins potential \eqref{poten} is adopted, a different bound-preserving mechanism is present for the CHNS model. Indeed, the singular function $f(\phi)$ is shown to be bounded in the $L^2$ norm, cf. \cite{Abels2009, Miranville2019}. The singularity ($f(\pm1) = \pm \infty$) implies that $\phi \in (-1, 1)$ if the solution exists. The discrete version of this nonlinear $L^2$ estimate is exploited in \cite{CoEl1992} for a Backward Euler time-marching scheme with the $C^0$ piecewise linear finite element for the Cahn-Hilliard equation. Recently, a convex-concave splitting finite difference scheme was proposed in \cite{CWWW2019}, for which a contradiction argument is devised for a direct $L^\infty$ estimate of $\phi$. The second-order scheme is designed and analyzed in \cite{CJWWW2022}. See \cite{CJWW2022, CJQWW2023} for similar schemes for the Flory-Huggins-Cahn-Hilliard-Navier-Stokes system. Other recent developments in the design of positivity-preserving schemes for gradient flows include the Lagrange multiplier approach \cite{ChSh2022}, the nonlinear change of variable \cite{HSW2022}, the JKO-type methods for optimal transport \cite{FOL2023, CWW2024}, the entropy estimate \cite{GuTi2024}, and the exponential differencing method \cite{DJLQ2021}. A recent survey on finite element methods respecting the maximum principle for linear convection-diffusion equations is provided in \cite{BJK2024}.

Current numerical methods for the Flory-Huggins-Cahn-Hilliard-Navier-Stokes equations are limited to second-order accuracy in space. In this article, we develop a spatially high-order finite element method for the CHNS system \eqref{CH}--\eqref{poten} on rectangular domains. The spatial discretization utilizes $Q_k$ finite elements with mass-lumping. The time-marching is based on the Crank-Nicolson method with convex-concave splitting for the potential. This discretization, however, does not preserve the $L^2$ estimate of the singular function. Instead, the singular function is added back as a perturbation, cf. \cite{CJWW2022, CJQWW2023}. This allows for a discrete $L^1$ estimate of the singular function, thereby preserving the desired bounds. Following the monotone argument in \cite{HaWa2015}, we establish the unconditional unique solvability of the nonlinear scheme. Compared to the second-order finite difference schemes in \cite{CJWW2022, CJQWW2023}, the discrete $L^1$ estimate enables us to construct a spatially high-order finite element method on rectangular grids. The convergence analysis can be done following \cite{DWWW2017}, which is left for future work. The scheme and its analysis can be extended to polygonal or polyhedral domains by using finite element spaces on quadrilateral meshes based on bilinear transformation \cite{ABF2002}.

The rest of the paper is organized as follows. In Section \ref{sec2}, we present the high-order stable numerical scheme, analyze its unique solvability, bounds-preserving properties, as well as energy stability. In Section \ref{sec3}, several numerical tests are performed to validate the accuracy and efficiency of the proposed numerical method.

\section{{The numerical scheme}}\label{sec2}
{Assume a uniform partition of the time interval $[0,T]$: $0=t_0<t_1<\ldots<t_{N_T}$ with time step size $\delta t=t_{n+1}-t_n$, for $0\leq n \leq N_T-1$,  $N_T=[T/\delta t]$.  The $L^2$ inner is denoted by $(\cdot,\cdot)$. For convenience we adopt the following  notations }
\begin{subequations}\label{notas}
	\begin{align}
		&\phi^{n+\frac{1}{2}}=\frac{1}{2}(\phi^{n+1}+\phi^n), \quad \tphi^{n+\frac{1}{2}}=\frac{3\phi^n-\phi^{n-1}}{2},  \\
		&\baru^{n+\frac{1}{2}}=\frac{\baru^{n+1}+\Au^n}{2}, \quad \tu^{n+\frac{1}{2}}=\frac{3\Au^n-\Au^{n-1}}{2},\\
		&\label{trilinear}
		b(\Au,\bm{v},\bm{w})=\frac{1}{2}\{(\Au\cdot \nabla \bm{v}, \bm{w})-(\Au\cdot \nabla \bm{w}, \bm{v})\}.
	\end{align}
\end{subequations}

Let $\mathcal{T}_h$ be a mesh of size $h$ consisting of rectangular cells of the rectangular domain $\Omega$.  Generalization to three dimension is straightforward. We introduce the following finite element spaces based on tensor products of $1$D polynomials of degree less than equal to $k$ on rectangular cells:
\begin{subequations}\label{notas2}
\begin{align}
&Y_h^k=\{v \in C^0{(\mathcal{T}_h)}:  v|_{K} \in Q_k(K), \forall K\in \mathcal{T}_h\}, \\
&\mathbf{X}_h^k=\{\bfv \in {Y_h^k\times Y_h^k}: \bfv|_{\partial \Omega}=0\}, \\
&M_h^{k-1}=Y_h^{k-1} \cap L^2_0(\Omega):= \{q_h \in {Y_h^{k-1}}; \int_{\Omega}q_h \bdx=0\}. 
\end{align}
\end{subequations}
Following \cite{GuQu1998a} one introduces discrete (negative) divergence operator
$B_h: \mathbf{X}_h^k   \rightarrow M_h^{k-1}$ ($\subset H^1$,  endowed with $L^2$ norm)  such that for $\Au_h \in \mathbf{X}_h^k$ and $ q_h \in M_h^{k-1}$
\begin{align}\label{dis-op}
(B_h \Au_h, q_h):=-(\nabla \cdot \Au_h, q_h)=(\Au_h, \nabla q_h):= (\Au_h, B_h^T q_h),
\end{align}
where $B_h^T$ is the transpose of $B_h$ (discrete gradient operator). It is known \cite{FBrezzi_RSFalk_SINN_1991} that the generalized Taylor-Hood pair $(\mathbf{X}_h^k, M_h^{k-1})$
is stable in the sense that it satisfies the inf-sup condition ($B_h$ is surjective).

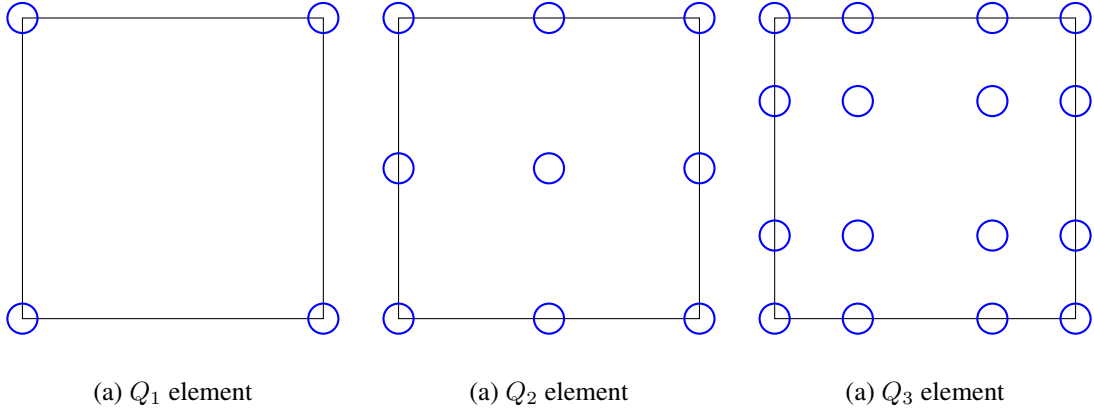
\begin{figure}
	\centering
	\begin{tikzpicture}
	\draw (0,0) rectangle (4,4); 
	\draw (5,0) rectangle (9,4);
	\draw (10,0) rectangle (14,4);
	\begin{scope}[thick, scale=1]
	\draw [blue] (0, 0) circle (0.2);
	\draw [blue] (0,4)  circle (0.2);
	\draw [blue] (4, 0)  circle (0.2);
	\draw [blue] (4, 4)  circle (0.2);
	\draw [blue] (5, 0) circle (0.2);
	\draw [blue] (7, 0) circle (0.2);
	\draw [blue] (5,2)  circle (0.2);
	\draw [blue] (5,4)  circle (0.2);
	\draw [blue] (7,4)  circle (0.2);
	\draw [blue] (9, 0)  circle (0.2);
	\draw [blue] (9, 2)  circle (0.2);
	\draw [blue] (9, 4)  circle (0.2);
	\draw [blue] (7, 2)  circle (0.2);
	\draw [blue] (10, 0) circle (0.2);
	\draw [blue] (14,0)  circle (0.2);
	\draw [blue] (14, 4)  circle (0.2);
	\draw [blue] (10, 4)  circle (0.2);
	\draw [blue] (11.1056, 0) circle (0.2);
	\draw [blue] (12.8944, 0) circle (0.2);
	\draw [blue] (14,1.1056) circle (0.2);
	\draw [blue] (14,2.8944) circle (0.2);
	\draw [blue] (12.8944, 4) circle (0.2);
	\draw [blue] (11.1056, 4) circle (0.2);
	\draw [blue] (10,2.8944) circle (0.2);
	\draw [blue] (10,1.1056) circle (0.2);
	\draw [blue] (12.8944, 1.1056) circle (0.2);
	\draw [blue] (11.1056, 1.1056) circle (0.2);
	\draw [blue] (12.8944, 2.8944) circle (0.2);
	\draw [blue] (11.1056, 2.8944) circle (0.2);
	\end{scope}
	\node (A) at (2,-1) {(a) $Q_1$ element };
	\node (B) at (7,-1) {(a) $Q_2$ element };
	\node (C) at (12.0,-1) {(a) $Q_3$ element };
	\end{tikzpicture}
	
	\caption{Spectral finite elements for quadrilateral element with degree $k=1,2,3$ in two dimensions, $\circ$ represents Gauss-Lobatto node.} \label{fig2:QK}
\end{figure}

{Three typical $Q_k$ finite elements are illustrated in Figure \ref{fig2:QK}.} 
Denote by $Z_0=\{\bm{x}_i, i\in J\}$ the set of tensor products of $(k+1)$ Gauss-Lobato quadrature nodes for all rectangular cells. 
One introduces the interpolation operator $I_h: C(\overline{\Omega}) \rightarrow Y_h^k$ such that $I_h(v)(\bm{x}_i)=v(\bm{x}_i), \forall \bm{x}_i\in Z_0$.
 One then defines the semi-inner product on $C(\overline{\Omega})$
\begin{align*}
	(u, v)_h:=\sum_{K\in \mathcal{T}_h}\sum_{j\in J_K} w_j u (\bm{x}_j) v(\bm{x}_j)=\sum_{i\in J} \tilde{w}_i u (\bm{x}_i) v(\bm{x}_i)=\int_{\Omega} I_h(uv)\, \bdx
\end{align*}
where $w_j$ is the (tensor product) Gauss-Lobatto quadrature weight ($w_j>0$), $\tilde{w}_i=\int_\Omega \chi_i\, \bdx$ and $\chi_i$ is the nodal Lagrange basis function such that $\chi_i(\bm{x}_j)=\delta_{ij}$.  The induced semi-norm, denoted by $|| \cdot||_h$,  is equivalent to the $L^2$ norm on $Y_h^k$.  One also introduces the discrete $L^1$ norm on $Y_h^k$ defined as
\begin{align*}
	||u_h||_{1,h}:=(|u_h|, 1)_h=\sum_{i\in J} \tilde{w}_i |u_h(\bm{x}_i)|.
\end{align*}It is well-known \cite{Ciarlet2002} that the Galerkin method with numerical integration using tensor product of (k+1) Gauss-Lobato quadrature nodes yields optimal convergence rates for elliptic problems.

One defines a truncation function $[x]$  such that 
\begin{align*}
	[x]=\begin{cases}
		1, \quad x\geq 1,\\
		x, \quad x \in(-1, 1), \\
		-1, \quad x\leq -1.
	\end{cases}
\end{align*}Performing the convex-concave splitting, one writes
\begin{align*}
	F(\phi)=F_v(\phi)+F_n(\phi)=\frac{1}{2}[(1+\phi)\ln(1+\phi)+(1-\phi)\ln(1-\phi)]-\frac{\theta_c}{2}\phi^2,
\end{align*}
with $f_v=F_v^\prime, f_n=F_n^\prime$, i.e. 
\begin{align} \label{singular}
	&f_v(\phi)=\frac{1}{2}\big[\ln {(1+\phi)}-\ln{(1-\phi)}\big], \quad f_n(\phi)=-\theta_c \phi.
\end{align}
For $f_v$ there holds:
\begin{lemma}
	Suppose $|\phi_h^n(\bm{x}_i)|\leq 1-\delta, \forall \bm{x}_i \in Z_0$ for sufficiently small $\delta>0$. Then there exist constants $C_1, C_2>0$ independent of  $\delta$ such that
	\begin{align}\label{Elem}
		f_v(s) (s-\phi_h^n (\bm{x}_i)) \geq C_1 \delta |f_v(s)|-C_2, \quad \forall s\in (-1, 1), \forall \bm{x}_i \in Z_0.
	\end{align}
\end{lemma}
The proof is similar to  \cite[Proposition 4.3]{Miranville2019}, and the details are omitted here. For a large $N>0$, one introduces a truncation of $f_v$ as follows \cite[pp 63]{Miranville2019}:
\begin{align*}
	f_v^N(s)=\begin{cases}
		f_v(-1+N^{-1})+f_v^\prime(-1+N^{-1}) (s+1-N^{-1}), \quad s<-1+N^{-1}, \\
		f_v(s), \quad |s|\leq 1-N^{-1}, \\
		f_v(1-N^{-1})+f^\prime_v(1-N^{-1})(s-1+N^{-1}), \quad s>1-N^{-1},
	\end{cases}
\end{align*}
and the associated convex potential is defined as $F_v^N:=\int_0^s f_v^N(x)\, dx$. One notes that $f_v^N(s)$ satisfies a similar  inequality as \eqref{Elem} for $-\infty<s<+\infty$. 

The high-order numerical scheme is as following:  find $(\phi_h^{n+1}, \mu^{n+\frac{1}{2}}_h, \baru_h^{n+\frac{1}{2}}, p_h^{n+1}, \Au_h^{n+1}) \in Y_h^k \times Y_h^k \times \mathbf{X}_h^k \times M_h^{k-1} \times \mathbf{X}_h^k$ such that for all $(v_h, \varphi_h, \bm{v}_h, q_h) \in Y_h^k \times Y_h^k \times \mathbf{X}_h^k \times Y_h^{k-1}$ there holds
\begin{align}
&\big(\phi^{n+1}_h-\phi^n_h,v_h\big)_h+\delta t \big(M([\tphi^{n+\frac{1}{2}}_h])\nabla \mu^{n+\frac{1}{2}}_h, \nabla v_h\big)-\delta t \big(\tphi^{n+\frac{1}{2}}_h\baru^{n+\frac{1}{2}}_h, \nabla v_h\big)=0, \label{d2ndCHNSwa}\\
&\big(\mu^{n+\frac{1}{2}}_h, \varphi_h \big)_h=\left(\frac{F_v(\phi_h^{n+1})-F_v(\phi_h^{n})}{\phi_h^{n+1}-\phi_h^{n}},  \varphi_h\right)_h +\big(f_n(\tp_h), \varphi\big)+\epsilon^2(\nabla \ap_h, \nabla \varphi_h)\nonumber \\
&+\delta t \big(f_v(\phi_h^{n+1})-f_v(\phi_h^n), \varphi_h\big)_h \label{d2ndCHNSwb}\\
&\big(2\baru^{n+\frac{1}{2}}_h, \bm{v}_h \big)+ \frac{\delta t}{Re}\big(\eta([\tp_h]) \nabla\baru^{n+\frac{1}{2}}_h, \nabla \bm{v}_h\big)+\delta t b\big(\tu^{n+\frac{1}{2}}_h, \baru^{n+\frac{1}{2}}_h, \bm{v}_h\big) =-\delta t\big(\nabla p^n_h, \bm{v}_h\big) \nonumber \\
&+\big(2\Au^n_h, \bm{v}_h \big)-\delta t \frac{\epsilon^{-1}}{\sw} \big(\tphi^{n+\frac{1}{2}}_h \nabla \mu^{n+\frac{1}{2}}_h, \bm{v}_h\big), \label{d2ndCHNSwc}  \\
&\big(\bfu^{n+1}_h-{\baru^{n+1}_h}, \bm{v}_h\big)+ \frac{\delta t}{2}\big(\nabla (p_h^{n+1}-p_h^n), \bm{v}_h\big)+\big(\nabla \cdot \bfu^{n+1}_h, q_h\big)=0. \label{d2ndCHNSwe}
\end{align}

One initializes the scheme with the following first order scheme
\begin{align}
	&\big(\phi^{1}_h-\phi^0_h,v_h\big)_h+\delta t \big(M(\phi^0_h)\nabla \mu^{1}_h, \nabla v_h\big)-\delta t \big(\phi^0_h\bm{u}_h^0, \nabla v_h\big)=0, \label{1stCHNSwa}\\
	&\big(\mu^{1}_h, \varphi_h \big)_h= \big(f_v(\phi_h^{1}), \varphi_h\big)_h+\big(f_n(\phi_h^0), \varphi\big)+\epsilon^2(\nabla \phi_h^1, \nabla \varphi_h),  \label{1stCHNSwb}\\
	&\big(\bm{u}^0_h, \bm{v}_h \big)+ \frac{\delta t}{Re}\big(\eta(\phi_h^0) \nabla\bm{u}_h^1, \nabla \bm{v}_h\big)+\delta t b\big(\bm{u}_h^0-\Au^0_h, \bm{u}_h^1, \bm{v}_h\big) =-\delta t\big(\nabla p^1_h, \bm{v}_h\big) \nonumber \\
	&-\delta t \frac{\epsilon^{-1}}{\sw} \big(\phi^{0}_h \nabla \mu^1_h, \bm{v}_h\big), \label{1stCHNSwc}  \\
	&\big(\nabla \cdot \bfu^{n+1}_h, q_h\big)=0.  \label{1stCHNSwd} 
\end{align}
The first order scheme is a decoupled method. The unique solvability follows from the monotone argument in the proof of Theorem \ref{USS}. Moreover, assuming $\phi_h^0 \in Y_h^k$ with $|\phi_h^0| \leq 1 $ a. e. and $\overline{\phi_h^0}:=\frac{1}{|\Omega|} \int_\Omega \phi_h^0\, dx \in (-1, 1)$, one derives from the $L^1$ estimate of Eq. \eqref{1stCHNSwb} (taking $\varphi_h=\phi_h^1-\overline{\phi_h^0}$) that there exists a small constant $\delta_1>0$ such that
\begin{align*}
		|\phi_h^1 (\bm{x}_i)|\leq  1-\delta_1, \quad \forall \bm{x}_i \in Z_0.
\end{align*}

{Now we show the proposed scheme \eqref{d2ndCHNSwa}--\eqref{d2ndCHNSwe}  is uniquely solvable, bound-preserving  and obeys discrete energy dissipation law without any time step constraint. }
\begin{theorem}
	\label{USS}
	Assume that the scheme  \eqref{d2ndCHNSwa}--\eqref{d2ndCHNSwe} is initialized by the first-order method \eqref{1stCHNSwa}--\eqref{1stCHNSwd} with $\phi_h^0 \in {Y_h^k}$, $|\phi_h^0|\leq1$ a.e. and  $\overline{\phi_h^0}:=\frac{1}{|\Omega|} \int_\Omega \phi_h^0\, dx \in (-1, 1)$. Then for any $h, \delta t>0$ the scheme \eqref{d2ndCHNSwa}--\eqref{d2ndCHNSwe} admits a unique solution such that
	\begin{align*}
		|\phi_h^n (\bm{x}_i)|< 1, \quad \forall \bm{x}_i \in Z_0, \quad n=1, 2 \ldots N_T.
	\end{align*}
	 Furthermore  the following discrete energy law holds: 
	\begin{align}\label{Ener_dis}
		E^{n+1} +	\frac{\delta t}{\epsilon We^\ast}||\sqrt{M}\nabla \mu^{n+\frac{1}{2}}_h||^2+\frac{\delta t}{Re}||\sqrt{\eta}\nabla \baru^{n+\frac{1}{2}}_h||^2 \leq E^n,
	\end{align}
where the modified energy is defined as
\begin{align*}
	E^n:=\frac{1}{\epsilon We^\ast}\Big[\big(F(\phi^{n}_h),1\big)_h+\frac{\epsilon^2}{2}||\nabla \phi^{n}_h||^2+\frac{\theta_c}{4}||\phi^{n}_h-\phi^{n-1}_h||^2\Big]+
	\frac{1}{2}||\bfu^{n}_h||^2+	\frac{\delta t^2}{8}||B_h^T p^{n}_h||^2.
\end{align*}
\end{theorem}

\begin{proof}
	One divides the proof into four steps. Without ambiguity one temporarily drops the dependence of the solution variables on $h$ and $n$.
	
	\noindent{\bf{Step 1.}}  One first recalls from \cite[Lemma 3.2]{HaWa2015} that for a given $\mu \in Y_h^k$ there exists a unique solution $\overline{\bfu} \in \mathbf{X}_h^{k}$ to Eq. \eqref{d2ndCHNSwc}, the solution is bounded and depends continuously on $\mu$.  
	
	Next, one would like to establish a similar result to Eq. \eqref{d2ndCHNSwb}.  Due to the singularity of $f_v$ at $-1$ and $1$, one considers an approximate system \cite{Miranville2019}:
	\begin{align}
		\label{che-appro}
		&(\mu, \varphi)_h=\left(G^N(\phi^N),  \varphi\right)_h +\big(f_n(\tp_h), \varphi\big)+\frac{\epsilon^2}{2}(\nabla \phi^N, \nabla \varphi)+\frac{\epsilon^2}{2}(\nabla \phi^{n}_h, \nabla \varphi)\nonumber \\
		&+\delta t \big(f_v^N(\phi^N)-f_v(\phi_h^n), \varphi\big)_h, \quad  \forall \varphi \in Y_h^k,
	\end{align}
	where $G^N(\phi^N):=\frac{F_v^N(\phi^N)-F_v^N(\phi_h^{n})}{\phi^N-\phi_h^{n}}$ is an increasing function of  $\phi^N$ thanks to the convexity of $F_v^N$. For a given $\mu \in Y_h^k$, Eq. \eqref{che-appro} is uniquely solvable since it is the Euler-Lagrange equation of the following strictly convex functional
	\begin{align*}
		&L(\phi):=\left(\int_0^\phi G^N(s)\, ds,  1\right)_h +\delta t \left( F_v^N (\phi), 1\right)_h+\int_\Omega \frac{\epsilon^2}{4}|\nabla \phi|^2\,\bdx\\
		&-\left(\mu+\delta t f_v(\phi_h^n), \phi\right)_h+\left(f_n(\tp_h)-\frac{\epsilon^2}{2} \Delta_h \langle\phi_h^{n}\rangle, \phi\right), 
	\end{align*}
	on the admissible set 
	\begin{align*}
		A_h:=\Big\{\phi \in Y_h^k,  \int_{\Omega} \phi\,\bdx=\int_{\Omega} \phi_h^n\,\bdx \Big\}.
	\end{align*}
	Here the discrete Laplacian operator $-\Delta_h: Y_h^k \cap L_0^2 \rightarrow Y_h^k\cap L_0^2$ is such that 
	\begin{align*}
		(-\Delta_h \phi, \varphi)=(\nabla \phi, \nabla \varphi), \quad \forall \varphi\in {Y_h^k,}
	\end{align*}
	and $\langle \phi\rangle:=\phi-\frac{1}{|\Omega|} \int_\Omega \phi\, \bdx$.
	
	One proceeds to deriving some a priori estimates of $\phi^N$ independent of $N$.   By the inequality \eqref{Elem} one has
	\begin{align}
		\label{L1-es0}
		\left(f_v^N(\phi^N), \phi^N-\phi_h^n\right)_h&=\sum_{i\in J}\tilde{w}_i f_v^N\big(\phi^N(\bm{x}_i)\big)\big(\phi^N(\bm{x}_i) -\phi_h^n(\bm{x}_i)\big) \nonumber \\
		&\geq \sum_{i\in J}\tilde{w}_i [C_1 \delta_n |f_v^N\big(\phi^N(\bm{x}_i)\big)|-C_2] \nonumber \\
		&=C_1 \delta_n|| I_h\big(f_v^N(\phi^N)\big)||_{1,h}-C.
	\end{align}
Since $\int_{\Omega} \phi^N\,\bdx=\int_{\Omega} \phi_h^n\,\bdx$, one deduces by Poincare's inequality
\begin{align*}
	|(\mu,  \phi^N-\phi_h^n)_h|=|(\mu-\langle \mu\rangle,  \phi^N-\phi_h^n)_h| \leq C ||\nabla \mu|| (|| \nabla \phi^N||+||\nabla \phi_h^n||).
\end{align*}
Hence by taking $\varphi=\phi^N-\phi_h^n$ in Eq. \eqref{che-appro} one obtains
	\begin{align}
		&\left(F_v^N(\phi^N), 1\right)_h+ \frac{\epsilon^2}{4}||\nabla \phi^N||^2+C \delta t \delta_n|| I_h\big(f_v^N(\phi^N)\big)||_{1,h}\leq C ||\nabla\mu||^2+\left(F_v^N(\phi_h^n), 1\right)_h\nonumber \\
		&+C(||\tp_h||^2+||\nabla \phi_h^n||^2+|| I_h\big(f_v(\phi^n_h)\big)||_h^2+1). \label{ch-es}
	\end{align}
	This  implies uniform boundedness of $ ||\nabla \phi^N||$ and hence   $||\phi^N||_h$ by Poincar\'{e}'s inequality, independent of $N$. It follows up to a subsequence
	\begin{align*}
		&\phi^N \rightarrow \phi \text{ weakly in } H^1 \text{ and strongly in } L^2, \\
		& \phi^N (\bm{x}_i) \rightarrow \phi (\bm{x}_i), \quad i\in J.
	\end{align*}
	Since  for fixed $h$ all norms are equivalent in finite dimension,  Eq. \eqref{ch-es}  implies
	\begin{align}
		\label{pot-es-s}
		||I_h\big(f_v^N(\phi^N)\big)||_{L^\infty} \leq C(h, \delta t, \mu, \delta_n).
	\end{align}
	Hence by choosing $\delta_\mu$ small enough such that $f_v(1-\delta_\mu)>C(h, \delta t,  \mu, \delta_n) $, one obtains for $N>\delta_\mu^{-1}$ and $\forall i \in J$
	\begin{align}
		&f_v^N\big(\phi^N(\bm{x}_i)\big) \leq C(h, \delta t,  \mu, \delta_n) <f_v(1-\delta_\mu)=f_v^N(1-\delta_\mu), 	\label{pot-e}\\
		&f_v^N\big(\phi^N(\bm{x}_i)\big) \geq -C(h, \delta t,  \mu, \delta_n) >-f_v(1-\delta_\mu)=f_v(-1+\delta_\mu)=f_v^N(-1+\delta_\mu). 	\label{pot-es2}
	\end{align}
	The monotonicity of $f_v^N$ yields $	|\phi^N(\bm{x}_i)| \leq 1-\delta_\mu, \forall \bm{x}_i \in Z_0$. Therefore
	\begin{align}
		\label{pot-es2}
		|\phi(\bm{x}_i)| \leq 1-\delta_\mu, i \in J.
	\end{align}
	
	By the mean value theorem and the uniform boundedness of $\phi^N$ and $\phi_h^n$, one has 
	\begin{align*}
		|| G^N(\phi^N)||_h^2=\sum_{i=1}^{N_h} m_i \Big[f_v^N\big(\theta \phi^N(\bm{x}_i)+(1-\theta)\phi_h^n(\bm{x}_i)\big)\Big]^2 \leq C_\mu.
	\end{align*}
	Since $G^N \rightarrow G, F_v^N\rightarrow F_v, f_v^N \rightarrow f_v$ as $N\rightarrow \infty$, one can pass to the limit in Eq. \eqref{che-appro} and conclude that $\phi$ is the solution. 
	
	One immediately obtains the continuous dependence of $\phi$ on $\mu$ in $H^1$  in light of the monotonicity of $G$ and $f_v$. The claim is thus proved. 
	
	\noindent{\bf{Step 2.}}  One establishes the unique existence of solution in this step.
	One defines an operator $T:Y_h^k \rightarrow {Y_h^k}$ such that for a given $\mu \in Y_h^k$
	\begin{align}\label{opT}
		\langle T(\mu), \varphi \rangle= (\phi-\phi_h^n, \varphi)_h-\delta t(\tp_h \baru, \nabla \varphi)+\delta t\big(M([\tp_h])\nabla \mu, \nabla \varphi\big), \forall \varphi \in {Y_h^k},
	\end{align}
	where $\phi$ is the solution to Eq. \eqref{d2ndCHNSwb}, $\baru$ is the  solution to Eq. \eqref{d2ndCHNSwc}.
	It is clear from Step 1  that the operator $T$ is continuous and bounded. 
	
	For any $\mu, \nu \in Y_h^k$ one calculates
	\begin{align}\label{Mono}
		&\langle T(\mu)-T(\nu), \mu-\nu \rangle = \big(\phi_\mu-\phi_\nu,\mu-\nu\big)_h+\delta t ||\sqrt{M}\nabla(\mu-\nu)|| \nonumber\\
		&-\delta t \big(\tphi^{n+\frac{1}{2}}[\baru_\mu-\baru_\nu], \nabla (\mu-\nu)\big).
	\end{align}
	In light of monotonicity of $G$ and $f_v$, one derives from Eq. \eqref{d2ndCHNSwb} 
	\begin{align}
		\label{mono1}
		&(\phi_\mu-\phi_\nu, \mu-\nu)_h=\big(G(\phi_\mu)-G(\phi_\nu), \phi_\mu-\phi_\nu\big)_h+\frac{\epsilon^2}{2}||\nabla(\phi_\mu-\phi_\nu)||^2\nonumber \\
		&+\delta t\big(f_v(\phi_\mu)-f_v(\phi_\nu), \phi_\mu-\phi_\nu\big)_h \geq 0,
	\end{align}
	and that the equality holds if only if $\mu=\nu$. It follows from  Eq. \eqref{d2ndCHNSwc} that
	\begin{align*}
		-\delta t \big(\tphi^{n+\frac{1}{2}}[\baru_\mu-\baru_\nu], \nabla (\mu-\nu)\big)=\epsilon\sw\{2||\baru_\mu-\baru_\nu||^2+\frac{\delta t}{Re}||\sqrt{\eta}\nabla(\baru_\mu-\baru_\nu)||^2\}.
	\end{align*}
	Hence
	\begin{align}
		\langle T(\mu)-T(\nu), \mu-\nu \rangle \geq 0,
	\end{align}
	with equality if only if $\mu=\nu$. This establishes the strict monotonicity of the operator $T$.
	
	Next, one has
	\begin{align}\label{Coer}
		\langle T(\mu), \mu \rangle = \big(\phi-\phi^n_h,\mu\big)_h+\delta t \big(M\nabla \mu, \nabla \mu\big)-\delta t \big(\tphi^{k+\frac{1}{2}}\baru, \nabla \mu\big), \quad \forall \mu \in Y_h^k.
	\end{align} 
Taking $\varphi=\phi-\phi_h^n$ in Eq. \eqref{d2ndCHNSwb} and noting the monotonicity of $f_v$, one deduces by Poincar\`{e}'s inequality
\begin{align}
	&(\phi-\phi^n_h, \mu)_h \geq \frac{\epsilon^2}{4}||\nabla \phi||^2+\int_\Omega I_h F_v(\phi)\, \bdx \nonumber \\
	&-\Big( C||\nabla \phi^n_h||^2+\int_\Omega I_h F_v(\phi^n_h)\, \bdx +C ||\tp_h||^2\Big). \label{coer1}
\end{align}
Likewise Eq. \eqref{d2ndCHNSwc} implies
\begin{align}\label{Coer2nd}
	-\delta t \big(\tphi^{n+\frac{1}{2}}\baru, \nabla \mu\big)&=\epsilon \sw \left[2||\baru||^2+\delta t ||\sqrt{\eta}\nabla \baru||^2- \delta t \big(2\Au^k-\nabla p^n, \baru\big)\right]\nonumber \\
	&\geq C\left[||\baru||^2+||\sqrt{\eta}\nabla \baru||^2-(||\Au^n||^2+||\nabla p^n||^2)\right].
\end{align}
 Eq. \eqref{Coer} then becomes
\begin{align}\label{Coer2}
	\langle T(\mu), \mu \rangle \geq C||\nabla \mu||^2 +\frac{\epsilon^2}{4} ||\nabla \phi||^2+C(||\baru||^2+||\sqrt{\eta}\nabla \baru||^2)- C.
\end{align}
Finally we estimate $\int_\Omega \mu\, \bdx$.  By the mean value theorem,
\begin{align*}
	G\big(\phi(\bm{x}_i)\big)=\frac{F_v\big(\phi(\bm{x}_i)\big)-F_v\big(\phi_h^n(\bm{x}_i)\big)}{\phi(\bm{x}_i)-\phi_h^n(\bm{x}_i)}=f_v\big(\theta_i \phi(\bm{x}_i)+(1-\theta_i)\phi_h^n(\bm{x}_i)\big),
\end{align*}
with $\theta_i \in (0, 1), i\in J$.  Denote $\theta_m=\max_{i\in J}{\theta_i}$.  One has
\begin{align}
	\label{G-L1}
	&\theta_m\left(G(\phi), \phi-\phi_h^n\right)_h\nonumber\\
	&=\sum_{i\in J}\tilde{w}_i f_v\big(\theta_i \phi(\bm{x}_i)+(1-\theta_i)\phi_h^n(\bm{x}_i)\big) \frac{\theta_m}{\theta_i}\big( \theta_i \phi(\bm{x}_i)+(1-\theta_i)\phi_h^n(\bm{x}_i)-\phi_h^n(\bm{x}_i)\big) \nonumber \\
	&\geq C\delta_n \sum_{i=\in J}\tilde{w}_i \big| f_v\big(\theta_i \phi(\bm{x}_i)+(1-\theta_i)\phi_h^n(\bm{x}_i)\big)\big| -C \nonumber \\
	&=C\delta_n\sum_{i\in J}\tilde{w}_i \big| G\big(\phi(\bm{x}_i)\big)\big| -C.
\end{align}
 Therefore by taking $\varphi=\theta_m(\phi-\phi_h^n)$ in Eq. \eqref{d2ndCHNSwb}  one obtains
\begin{align}
	\label{L1-es}
	\delta_n || I_h \big(G(\phi)\big) ||_{1,h}+ \theta_m \delta t || I_h \big(f_v(\phi)\big) ||_{1, h} + \frac{\epsilon^2}{4}|| \nabla \phi||^2 \leq C||\nabla \mu||^2+C.
\end{align}
{Eq. \eqref{d2ndCHNSwb} then implies  }
\begin{align*}
	\left|\int_\Omega \mu \, \bdx\right| \leq C (||\nabla \mu||+1),
\end{align*}
hence by Poincare's inequality
\begin{align}
	\label{muH1}
	||\mu||_{H^1} \leq C (||\nabla \mu||+1).
\end{align}
The inequality \eqref{Coer2} becomes
\begin{align*}
	\langle T(\mu), \mu\rangle&\geq C ||\mu||_{H^1}^2-C.
\end{align*}
This establishes the coercivity of $T$.

It follows from the Browder-Minty lemma \cite{HaWa2015} that there exists a unique solution $\mu \in Y_h^k$ such that 
\begin{align*}
	0=	\langle T(\mu), \varphi \rangle=(\phi-\phi^n_h, \varphi)_h-\delta t(\tp_h \baru, \nabla \varphi)+\delta t(M\nabla \mu, \nabla \varphi), \quad \forall \varphi \in Y_h^k.
\end{align*}

\noindent{\bf{Step 3:}} In this step one proves the energy law.   It is necessary to bring back the dependence on $h$ and $n$ in the notation.

Since
\begin{eqnarray}
	& &\big(\tphi^{n+\frac{1}{2}}_h,\phi^{n+1}_h-\phi^n_h\big)=\frac{1}{2}\big(3\phi^n_h-\phi^{n-1}_h,\phi^{n+1}_h-\phi^n_h\big) \nonumber \\
	&=&\frac{1}{2}\big(||\phi^{n+1}_h||^2-||\phi^n_h||^2\big)-\frac{1}{4}\Big[||\phi^{n+1}_h-\phi^n_h||^2-||\phi^n_h-\phi^{n-1}_h||^2\nonumber \\
	& & +||\phi^{n+1}_h-2\phi^n_h+\phi^{n-1}_h||^2\Big], \nonumber 
\end{eqnarray}
one takes $v_h= \frac{1}{\epsilon We^\ast}\mu_h^{n+\frac{1}{2}}$ in Eq. \eqref{d2ndCHNSwa}, $\varphi_h=\frac{1}{\epsilon We^\ast}\big(\phi_h^{n+1}-\phi_h^n\big)$ in Eq. \eqref{d2ndCHNSwb}, and $\bm{v}_h=\baru^{n+\frac{1}{2}}$ in Eq. \eqref{d2ndCHNSwc}, combines the results to obtain
\begin{align}\label{2ndCHf}
	&\frac{1}{\epsilon We^\ast}\Big[\big(F(\phi^{n+1}_h)-F(\phi^n_h),1\big)_h+\frac{\epsilon^2}{2}(||\nabla \phi^{n+1}_h||^2-||\nabla \phi^n_h||^2)+\frac{\theta_c}{4}(||\phi^{n+1}_h-\phi^n_h||^2-||\phi^{n}_h-\phi^{n-1}_h||^2)\Big]\nonumber \\
	&+\frac{1}{2}(||\baru^{n+1}_h||^2-||\bfu^n_h||^2)\leq -\frac{\delta t}{\epsilon We^\ast}||\sqrt{M}\nabla \mu^{n+\frac{1}{2}}||^2- \frac{\delta t}{Re}||\sqrt{\eta}\nabla \baru^{n+\frac{1}{2}}||^2-\delta t \big(\nabla p^n, \baru^{n+\frac{1}{2}}\big).
\end{align}

By utilizing the discrete divergence operator $B_h$ and the discrete gradient operator $B_h^T$ defined in \eqref{dis-op}, one writes the projection step Eqs. \eqref{d2ndCHNSwe}  as abstract equations
\begin{align}
	&\Au^{n+1}_h-\bm{\bar{u}}^{n+1}_h+ \frac{\delta t}{2}B_h^T (p_h^{n+1}-p_h^n)=0,  \text{ in } \mathbf{X}_h^k, \label{dis-pro}\\ 
	&B_h\Au^{n+1}_h=0, \text{ in } M_h^{k-1}. \label{dis-div}
\end{align}
One derives immediately
\begin{align}
	&\frac{1}{2}(||\Au^{n+1}_h||^2-||\baru^{n+1}_h||^2+||\Au^{n+1}_h-\baru^{n+1}_h||^2)=0, \label{2ndNS2}\\
	&\frac{{\delta t}^2}{8}||B_h^T(p^{n+1}_h-p^n_h)||^2=\frac{1}{2}||\Au^{n+1}_h-\baru^{n+1}_h||^2. \label{2ndNS4}
\end{align}
Since
\begin{align*}
	\frac{\Au^{n+1}_h+\Au^n_h-2\baru^{n+\frac{1}{2}}_h}{\delta t} + \frac{1}{2}B_h^T(p^{n+1}_h-p^n_h)=0, \text{ in } \mathbf{X}_h^k,
\end{align*}
testing the above equation with $\frac{\delta t^2}{2}B_h^T p^n$, using  Eq. \eqref{dis-div}, one has
\begin{align}\label{2ndNS3}
	\frac{\delta t^2}{8}\big[||B_h^T p^{n+1}_h||^2-||B_h^T p^{n}_h||^2-||B_h^T(p^{n+1}_h-p^n_h)||^2\big]=\delta t\big(\nabla p^n_h, \baru^{n+\frac{1}{2}}_h\big).
\end{align}

Summing up Eqs. \eqref{2ndCHf}, \eqref{2ndNS2}, \eqref{2ndNS4} and \eqref{2ndNS3},  one obtains
\begin{align}\label{2ndCHf-1}
	&\frac{1}{\epsilon We^\ast}\Big[\big(F(\phi^{n+1}_h)-F(\phi^n_h),1\big)_h+\frac{\epsilon^2}{2}(||\nabla \phi^{n+1}_h||^2-||\nabla \phi^n_h||^2)+\frac{\theta_c}{4}(||\phi^{n+1}_h-\phi^n_h||^2-||\phi^{n}_h-\phi^{n-1}_h||^2)\Big]\nonumber \\
	&+\frac{1}{2}(||\bfu^{n+1}_h||^2-||\bfu^n_h||^2)+	\frac{\delta t^2}{8}\big[||B_h^T p^{n+1}_h||^2-||B_h^T p^{n}_h||^2\big] \nonumber\\
	&\leq 	-\frac{\delta t}{\epsilon We^\ast}||\sqrt{M}\nabla \mu^{n+\frac{1}{2}}_h||^2- \frac{\delta t}{Re}||\sqrt{\eta}\nabla \baru^{n+\frac{1}{2}}_h||^2.
\end{align}
Take summation of \eqref{2ndCHf-1} from $n=1$ to $n=m-1$ to conclude 
\begin{align}\label{2ndCHf-2}
	&\frac{1}{\epsilon We^\ast}\Big[\big(F(\phi^{m}_h),1\big)_h+\frac{\epsilon^2}{2}||\nabla \phi^{m}_h||^2+\frac{\theta_c}{4}||\phi^{m}_h-\phi^{m-1}_h||^2\Big]+
\frac{1}{2}||\bfu^{m}_h||^2+	\frac{\delta t^2}{8}||B_h^T p^{m}_h||^2\nonumber\\
	&	\frac{\delta t}{\epsilon We^\ast}\sum_{n=1}^{m-1}||\sqrt{M}\nabla \mu^{n+\frac{1}{2}}_h||^2+\frac{\delta t}{Re}\sum_{n=1}^{m-1}||\sqrt{\eta}\nabla \baru^{n+\frac{1}{2}}_h||^2 \nonumber \\
	&\leq \frac{1}{\epsilon We^\ast}\Big[\big(F(\phi^{1}_h),1\big)_h+\frac{\epsilon^2}{2}||\nabla \phi^{1}_h||^2+\frac{\theta_c}{4}||\phi^{1}_h-\phi^{0}_h||^2\Big]
	+\frac{1}{2}||\bfu^{1}_h||^2+	\frac{\delta t^2}{8}||B_h^T p^{1}_h||^2.
\end{align}

\noindent{\bf{Step 4}.} In the final step one removes the $\mu$ dependency of $\delta$ in the estimates \eqref{ch-es}--\eqref{pot-es2}. One repeats the arguments in \eqref{L1-es0}--\eqref{pot-es2} for Eq. \eqref{d2ndCHNSwb} to obtain
		\begin{align*}
		&\left(F_v(\phi_h^{n+1}), 1\right)_h+ \frac{\epsilon^2}{4}||\nabla \phi^{n+1}_h||^2+C \delta t \delta_n || I_h\big(f_v(\phi^{n+1}_h)\big)||_{1,h}+C \delta t \delta_\mu || I_h\big(f_v(\phi^{n}_h)\big)||_{1,h}\\
		&\leq C ||\nabla\mu_h^{n+\frac{1}{2}}||^2+\left(F_v(\phi_h^n), 1\right)_h
		+C(||\tp_h||^2+||\nabla \phi_h^n||^2+1). 
	\end{align*}
In light of the energy law \eqref{2ndCHf-2}, one has
\begin{align*}
	|| I_h\big(f_v(\phi^{n+1}_h)\big)||_{1,h} \leq \frac{C \epsilon We^\ast}{\delta_n m_1 \delta t^2}, 
\end{align*}
hence 
\begin{align*}
	|| I_h\big(f_v(\phi^{n+1}_h)\big)||_{L^\infty} \leq \frac{C \epsilon We^\ast}{\delta_n m_1 h \delta t^2}. 
\end{align*}

Similar to the arguments in \eqref{pot-e} and \eqref{pot-es2}, one chooses $0<\delta_{n+1}<1$ such that $1-\delta_{n+1} >f_v^{-1}\Big(\frac{C \epsilon We^\ast}{\delta_n m_1 h \delta t^2}\Big)$ to conclude
\begin{align*}
	|\phi(\bm{x}_i)| \leq 1-\delta_{n+1}, \quad \forall \bm{x}_i \in Z_0.
\end{align*}
This completes the proof.

\end{proof}

\section{Numerical experiments}\label{sec3}
 
In this section, we perform several numerical experiments to verify the efficiency and capability of proposed numerical method. 
We first provide the numerical tests to validate the optimal convergence, energy stability as well as the boundedness of the proposed high order numerical scheme. The flexibility in simulating the evolution of binary fluid under rotational flow and lid-driven cavity flow is considered. 
The interesting Rayleigh-Taylor instability is implemented to investigate the robustness of the developed numerical scheme. 
The Newton's iteration is utilized to cope  the nonlinear term arising form the Flory-Huggins potential. Throughout, we take ${\bm{Q}}_k$--$Q_{k-1}$ Taylor-Hood  element for velocity and pressure,  $Q_k$ element for phase variable,  $\theta_c=2$.

\subsection{Convergence and accuracy}
The computational domain is chosen as a unit square, the exact solution are taken to be 
\begin{align*}
&\Au(x,y,t)=[-\cos(2\pi x)\sin(2 \pi y),\sin(2\pi x)\cos(2 \pi y)]^T\cos(\pi t),\\
&p(x,y,t)=\sin(\pi t)\sin(2\pi x),\\
&\phi(x,y,t)=\frac{1}{\pi}\cos(\pi t)\sin(2\pi x)\cos(2\pi y).
\end{align*}
The boundary conditions and source terms are modified corresponding to  the exact solution.  We set all parameters  to one.

\begin{figure}[!hbt]
	\centering
	\begin{tikzpicture}
	\draw (0,0) rectangle (5,5); 
	\draw (8,0) rectangle (13,5);
	\begin{scope}[thick, scale=1]
   \draw [blue] (0, 0) circle (0.2);
     \draw [blue] (0,2.5)  circle (0.2);
   \draw [blue] (0,5)  circle (0.2);
    \draw [blue] (2.5, 0)  circle (0.2);
    \draw [blue] (2.5, 2.5)  circle (0.2);
   \draw [blue] (5, 0)  circle (0.2);
   \draw [blue] (5, 2.5)  circle (0.2);
   \draw [blue] (5, 5)  circle (0.2);
    \draw [blue] (2.5, 5)  circle (0.2);
   \draw [blue] (8, 0) circle (0.2);
   \draw [blue] (8, 1.382) circle (0.2);
    \draw [blue] (8,3.618) circle (0.2);
   \draw [blue] (8,5)  circle (0.2);
   \draw [blue] (13, 0)  circle (0.2);
   \draw [blue] (13, 1.382)  circle (0.2);
   \draw [blue] (13, 3.618)  circle (0.2);
   \draw [blue] (13, 5)  circle (0.2);
    \draw [blue] (9.382, 0)  circle (0.2);
    \draw [blue] (11.618, 0)  circle (0.2);
     \draw [blue] (9.382, 5)  circle (0.2);
    \draw [blue] (11.618, 5)  circle (0.2);
    \draw [blue] (11.618, 1.382)  circle (0.2);
    \draw [blue] (11.618, 3.618)  circle (0.2);
     \draw [blue] (9.382, 1.382)  circle (0.2);
     \draw [blue] (9.382, 3.618)  circle (0.2);
	\end{scope}
	\node at (0, 0) {\huge${\color{red}\times}$};
	\node at (0, 5) {\Huge${\color{red}\times}$};
	\node at (5, 5) {\Huge${\color{red}\times}$};
	\node at (5,0) {\Huge${\color{red}\times}$};
	\node at (8, 0) {\Huge${\color{red}\times}$};
	\node at (10.5, 0) {\Huge${\color{red}\times}$};
	\node at (8, 2.5) {\Huge${\color{red}\times}$};
	\node at (8, 5) {\Huge${\color{red}\times}$};
	\node at (13, 0) {\Huge${\color{red}\times}$};
	\node at (13, 2.5) {\Huge${\color{red}\times}$};
	\node at (13,5) {\Huge${\color{red}\times}$};
		\node at (10.5, 2.5) {\Huge${\color{red}\times}$};
	\node at (10.5, 5) {\Huge${\color{red}\times}$};
	\node (A) at (2.5,-1) {(a) ${\bm{Q}}_2$--$Q_1$ element };
	\node (B) at (10.5,-1) {(a) ${\bm{Q}}_3$--$Q_2$ element };
    \node (D) at (15,1.5) {  ${\color{blue}\circ}$  velocity};
    \node (C) at (15,0.5) {${\color{red}\times}$ pressure};
	\end{tikzpicture}
	
	\caption{Taylor-Hood element ${\bm{Q}}_k$--$Q_{k-1}$ element,  $\circ$ denotes the velocity, $\times$ denotes the pressure.} \label{fig:QK}
\end{figure}
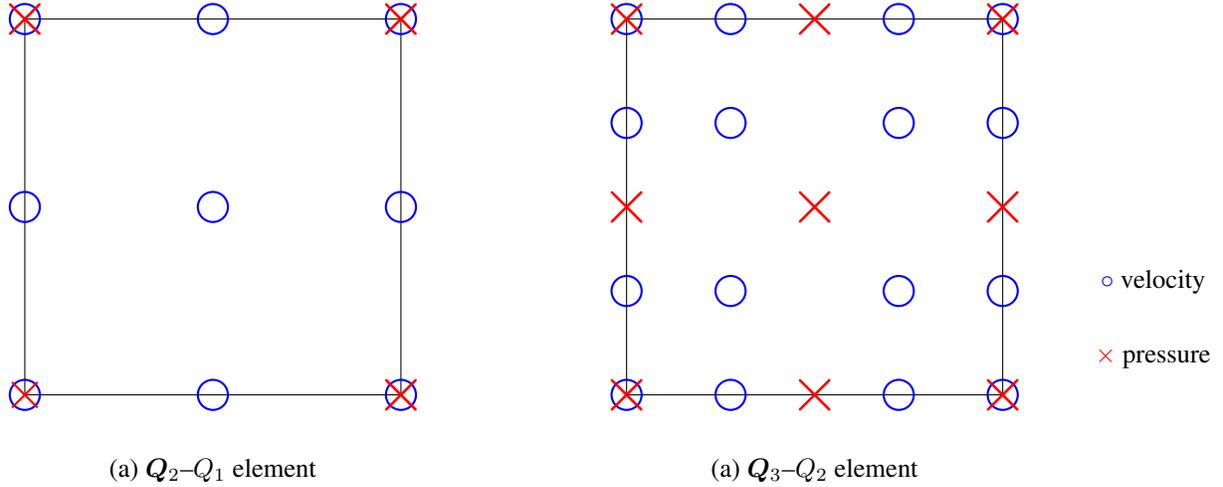

 The Taylor-Hood element are displayed in Figure \ref{fig:QK} for different type degree of freedom. Denote the numerical errors $|e_v|=|v_h^n-v(t_n)|$ $(v=\Au, p,\phi)$ between numerical solution $v_h^n$ and exact solution $v(t_n)$ with respect to $L^2$ norm $\|\cdot\|$ and $H^1$ norm  $\|\cdot\|_{H^1}$,  Tables \ref{table_CNLF_space1}-\ref{table_CNLF_space2} show the various error norms and convergence rate  at terminate time $T=1$ for ${\bm{Q}}_k$--$Q_{k-1}$--$Q_k$  defined in \eqref{notas2}   with polynomial degree $k=2,3$ corresponding to time step size $\delta t=h^\frac{3}{2}$ and $\delta t=4h^2$, respectively. 
As reported  in Tables \ref{table_CNLF_space1}-\ref{table_CNLF_space2},   the convergence order are of  the expected optimal accuracy  in space for both variables under different element types. 

For the temporal convergence, we utilized Cauchy convergence to eliminate the error from spatial discretization. Denote  Cacuhy difference  in $L^2$ norm by $\|e_v^h\|=\|v_h^{\delta t}-v_h^{\delta t/2}\|$  of numerical approximations between adjacent time steps  with fixed mesh size $h$. Choosing $h=\frac{1}{64}$, Tables   \ref{table_CNLF_time1}-\ref{table_CNLF_time2} display the numerical error and convergence rate, 
 which clearly show the proposed numerical method can achieve the excepted second order convergence rate  in time. These verify the arbitrary high order accuracy in space and second order accuracy in time of the proposed numerical scheme.

\begin{table}[!hbt]
	\begin{center}
			\caption{Numerical errors and  convergence  rates   by ${\bm{Q}}_2-Q_1-Q_2$ element pair with $\delta t=h^{3/2}$.} \label{table_CNLF_space1}
		\begin{tabular}{cp{1.9cm}<{\centering}p{0.8cm}<{\centering}p{1.9cm}<{\centering}p{0.8cm}<{\centering}p{1.9cm}<{\centering}p{0.8cm}<{\centering}p{1.9cm}<{\centering}p{0.8cm}<{\centering}p{1.9cm}<{\centering}p{0.8cm}<{\centering}p{1.9cm}}
		\toprule
	 1/h  &  $\|e_{\Au}\|$    &order       &  $\|e_p\|$    &order   &$\|e_\phi\|$   &order   \\[0.5ex]\hline
		8    &2.3102E-3   &--       &5.475E-3   &--   &5.1829E-4    &--     \\ 
		16   &2.9088E-4   &2.99        &1.0961E-3   &2.32  &6.5426E-5   &2.99      \\ 
		32   &3.6400E-5   &3.00     &1.3135E-4   &3.06  &8.1816E-6   &3.00     \\ 
		64   &4.5523E-6   &3.00       &1.9662E-5   &2.74   &1.0247E-6   &3.00    \\ 
		\bottomrule
		\end{tabular}
	\end{center}

\end{table}

\begin{table}[!hbt]
	\begin{center}
		\caption{Numerical errors and  convergence  rates   by ${\bm{Q}}_3-Q_2-Q_3$ element pair with $\delta t=4h^2$.}\label{table_CNLF_space2}
	\begin{tabular}{cp{1.9cm}<{\centering}p{0.8cm}<{\centering}p{1.9cm}<{\centering}p{0.8cm}<{\centering}p{1.9cm}<{\centering}p{0.8cm}<{\centering}p{1.9cm}<{\centering}p{0.8cm}<{\centering}p{1.9cm}<{\centering}p{0.8cm}<{\centering}p{1.9cm}}
		\toprule
		1/h  &  $\|e_{\Au}\|$    &order        &  $\|e_p\|$    &order   &$\|e_\phi\|$   &order       \\[0.5ex]\hline
		8   &1.7549E-4   &--      &1.2649E-2   &--   &4.9325E-5   &--    \\
		16   &1.1028E-5   &3.99   &7.8171E-5   &4.00   &7.7015E-4   &4.04     \\
		32   &6.9332E-7   &3.99      &4.8296E-5   &3.99   &2.0398E-7   &3.96   \\
		 64   &4.4236E-8   &3.97     &3.9805E-6   &3.60   &1.4183E-8   &3.85   \\
		\bottomrule
	\end{tabular}
\end{center}

\end{table}

\begin{table}[!hbt]
	\begin{center}
			\caption{Numerical errors and  convergence  rates   by ${\bm{Q}}_2-Q_1-Q_2$ element pair with fixed mesh size $h=\frac{1}{64}$.} \label{table_CNLF_time1}
		\begin{tabular}{cp{1.9cm}<{\centering}p{0.8cm}<{\centering}p{1.9cm}<{\centering}p{0.8cm}<{\centering}p{1.9cm}<{\centering}p{0.8cm}<{\centering}p{1.9cm}}
			\toprule
			$1/\delta t$  &  $\|e^h_{\Au}\|$    &order      &  $\|e^h_p\|$    &order   &$\|e^h_\phi\|$   &order   \\[0.5ex]\hline
		8   &4.5789E-5        &--     &3.7188E-2       &--   &1.8059E-5       &--\\
		16   &9.2878E-6   &2.30    &9.6497E-3   &1.95   &5.4505E-6   &1.73\\
		32   &1.0003E-6   &3.21   &2.3612E-3    &2.03   &1.4023E-6   &1.96\\
		64   &2.1183E-7   &2.24   &5.7955E-4   &2.03   &3.5371E-7   &1.99\\
			\bottomrule
		\end{tabular}
	\end{center}
\end{table}
\begin{table}[!hbt]
	\begin{center}
			\caption{Numerical errors and  convergence  rates   by ${\bm{Q}}_3-Q_2-Q_3$ element pair with fixed mesh size $h=\frac{1}{64}$.} \label{table_CNLF_time2}
		\begin{tabular}{cp{1.9cm}<{\centering}p{0.8cm}<{\centering}p{1.9cm}<{\centering}p{0.8cm}<{\centering}p{1.9cm}<{\centering}p{0.8cm}<{\centering}p{1.9cm}}
			\toprule
			$1/\delta t$  &  $\|e^h_{\Au}\|$    &order      &  $\|e^h_p\|$    &order   &$\|e^h_\phi\|$   &order   \\[0.5ex]\hline
			 8   &2.1434E-5   &--    &7.4171E-3   &--   &6.0206E-6   &--\\
			16   &5.0397E-6   &2.09    &1.9329E-3   &1.94   &1.7596E-6    &1.77\\
			32  &9.7396E-7   &2.37   &4.6981E-4   &2.04    &4.512E-7    &1.96\\
			64  &2.1158E-7   &2.20  &1.1084E-4   &2.08   &1.1442E-7    &1.98\\
			\bottomrule
		\end{tabular}
	\end{center}

\end{table}

\subsection{Energy dissipation}

In this test, we simulate the spinodal decomposition phenomenon of phase separation 
to validate the positivity preserving and  the energy stability.  It can be expected that this system will evolve from a random non-equilibrium state to a two-phase
state due to  the spontaneous   growth of concentration instability. 

The parameters are chosen as $\sw=50$, $Pe=100$, $Re=1$, $\eta=1$ and $\epsilon=0.02$ on computational domain $[0,1]^2$. The initial conditions for phase variable  is given by the concentration fields as $\phi_0=0.2-0.01r(\bm{x})$ with  randomly perturbed number $r(\bm{x})\in [-1,1]$. We impose no slip boundary conditions for velocity, and zero Neumann boundary conditions for $\phi$ and $\mu$ at all four boundaries. Choosing $\delta t=0.001$ and uniform quadrilateral mesh partition of mesh size $h=\frac{1}{128}$, the mobility is taken as constant  $M=0.1$.

\begin{figure}[ht]
	\begin{center}
		\subfigure[$t=0.5$] {\includegraphics[width=1.6in,trim=70 20 60 20,clip]{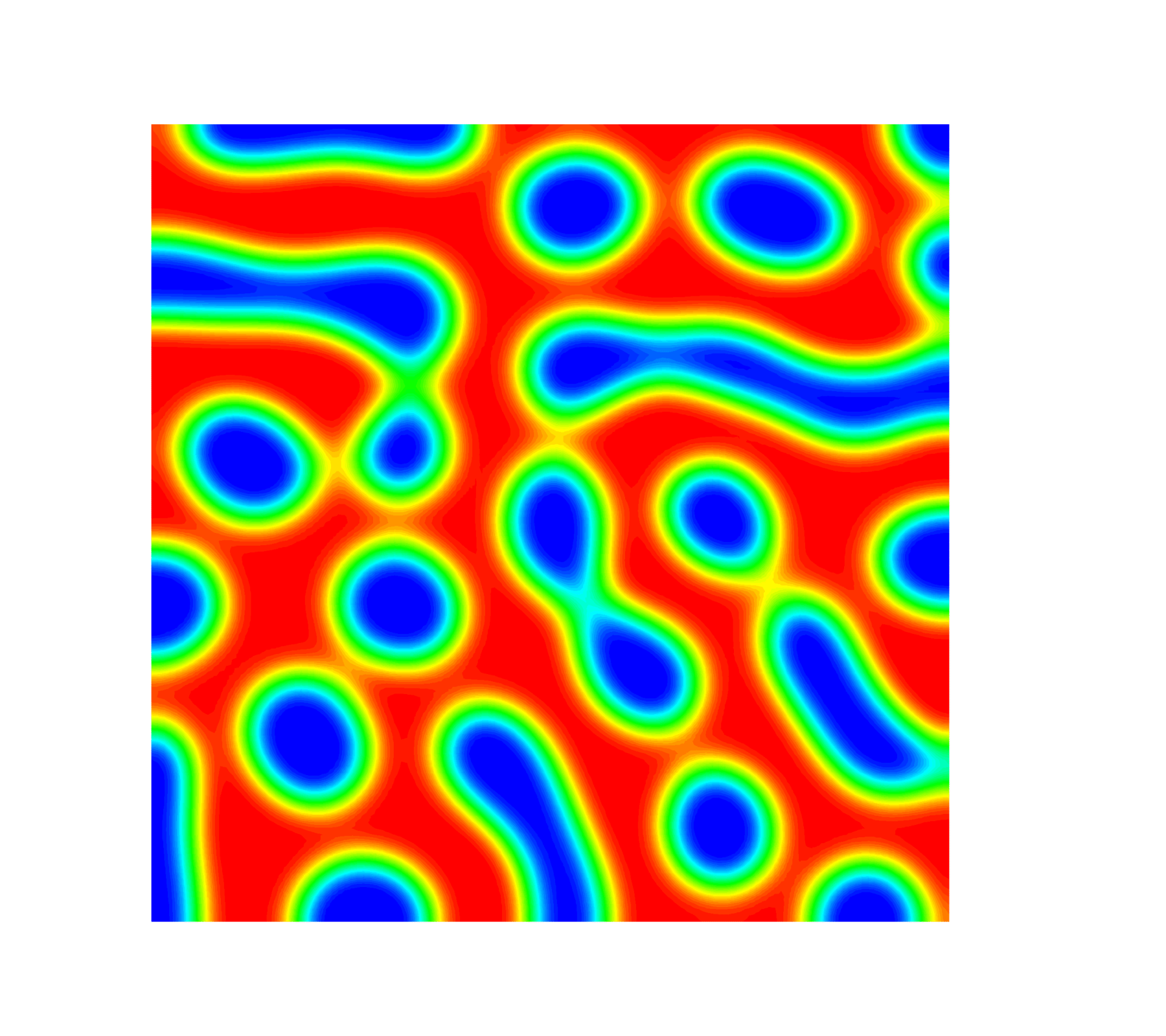}}
		\subfigure[$t=1$]   {\includegraphics[width=1.55in,trim=70 20 60 20,clip]{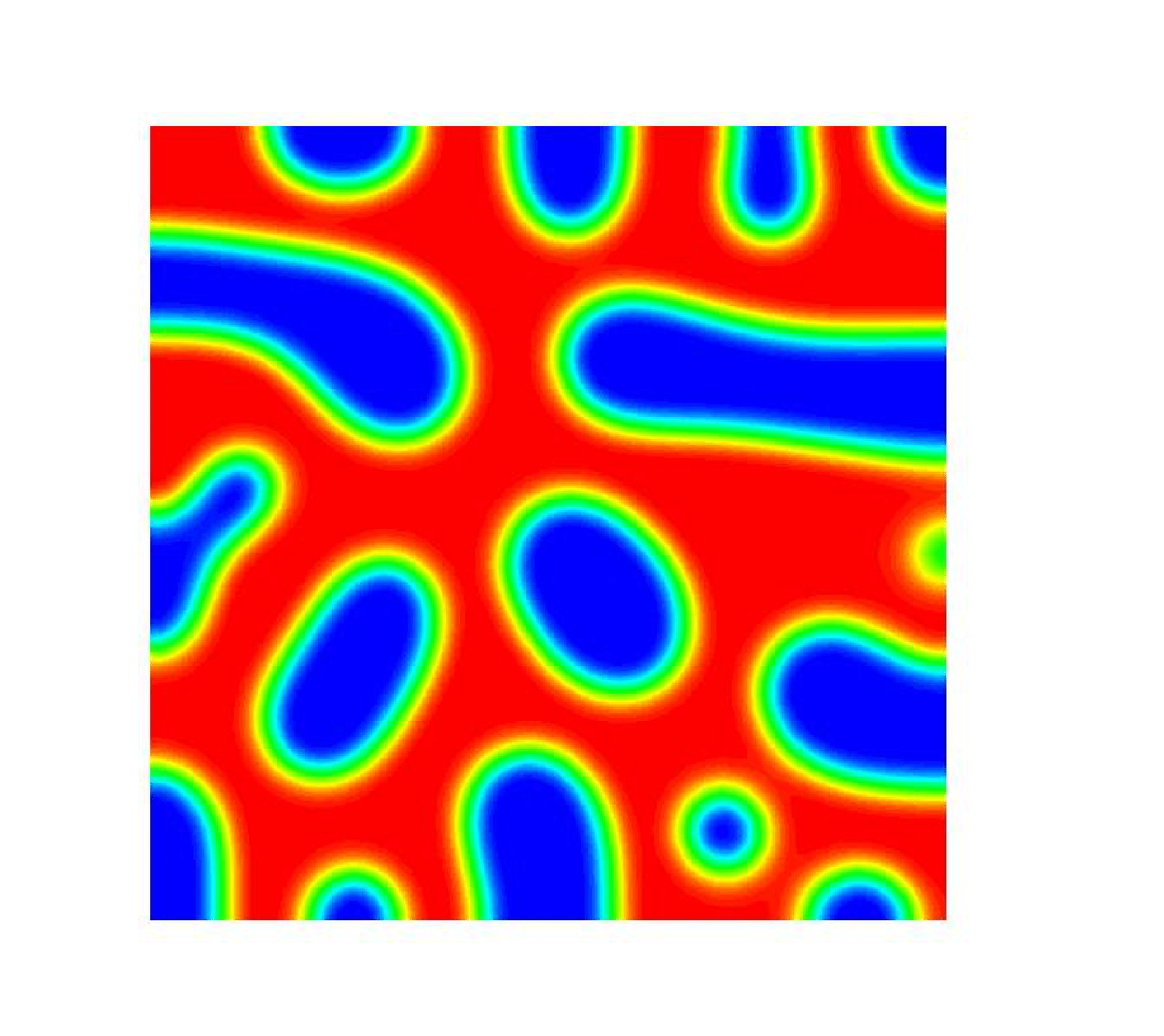}}
		\subfigure[$t=2$] {	\includegraphics[width=1.6in,trim=70 20 60 20,clip]{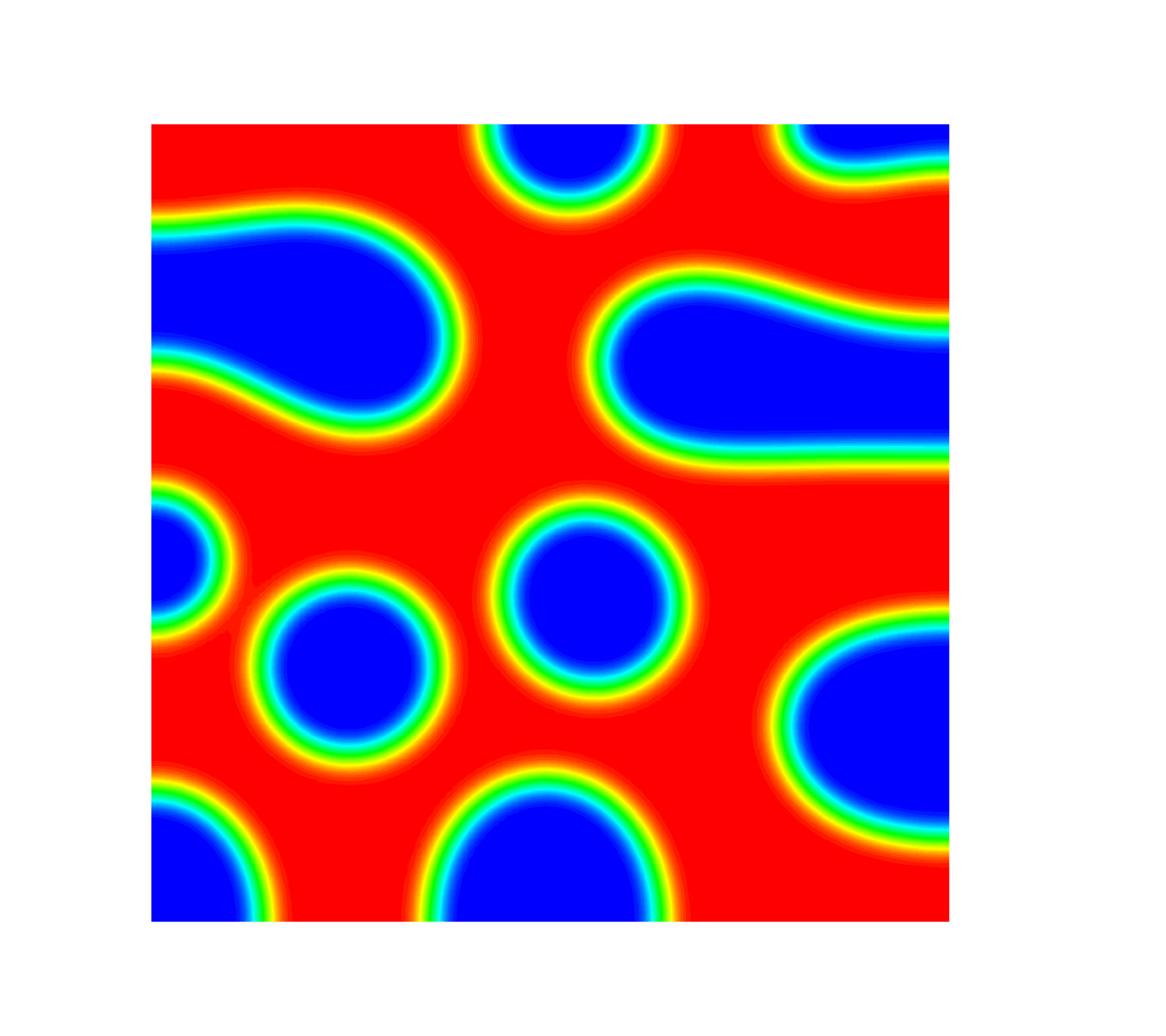}}
		\subfigure[$t=3.5$] {\includegraphics[width=1.6in,trim=70 20 60 20,clip]{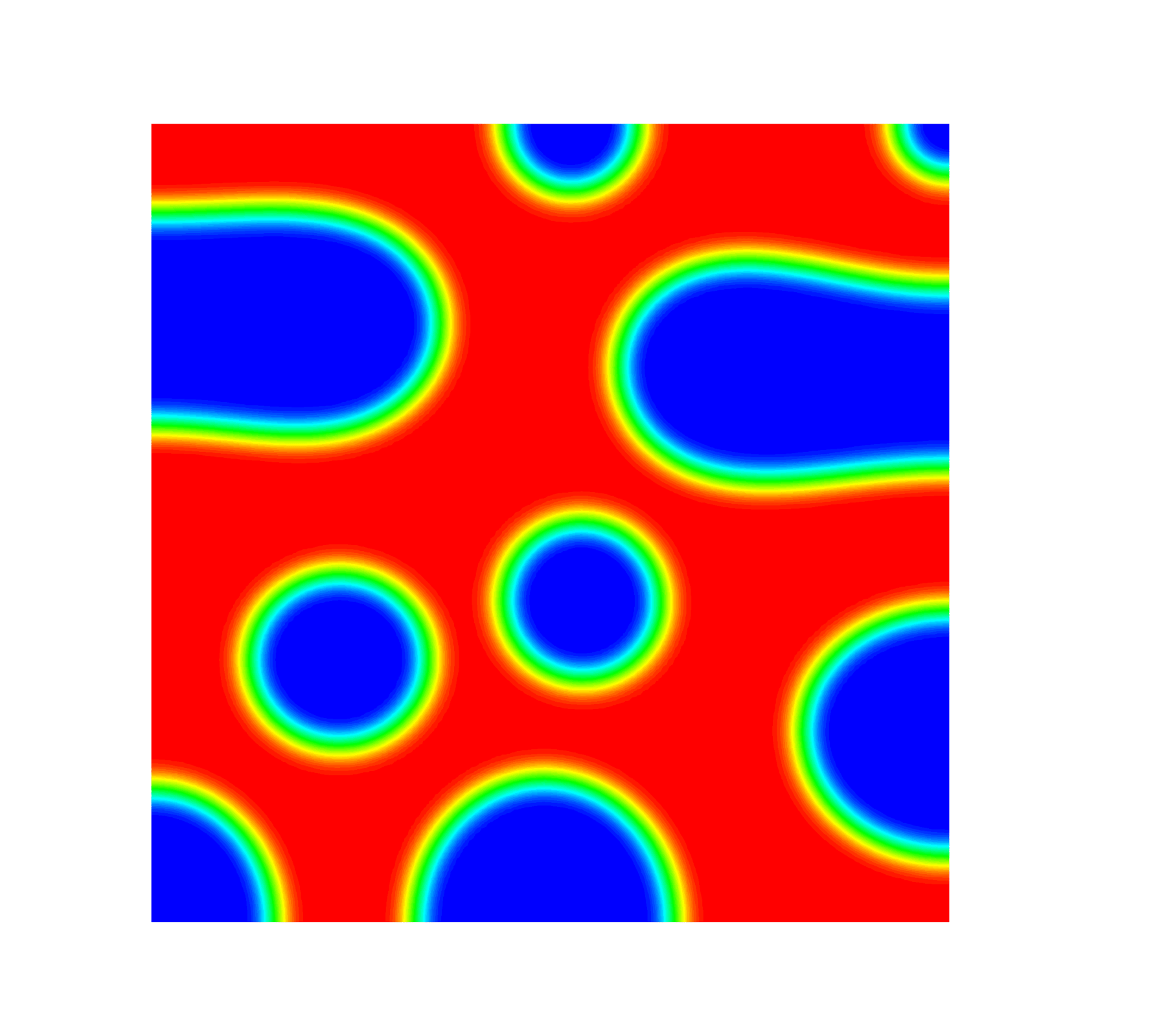}}
		\subfigure[$t=5$] {\includegraphics[width=1.6in,trim=70 20 60 20,clip]{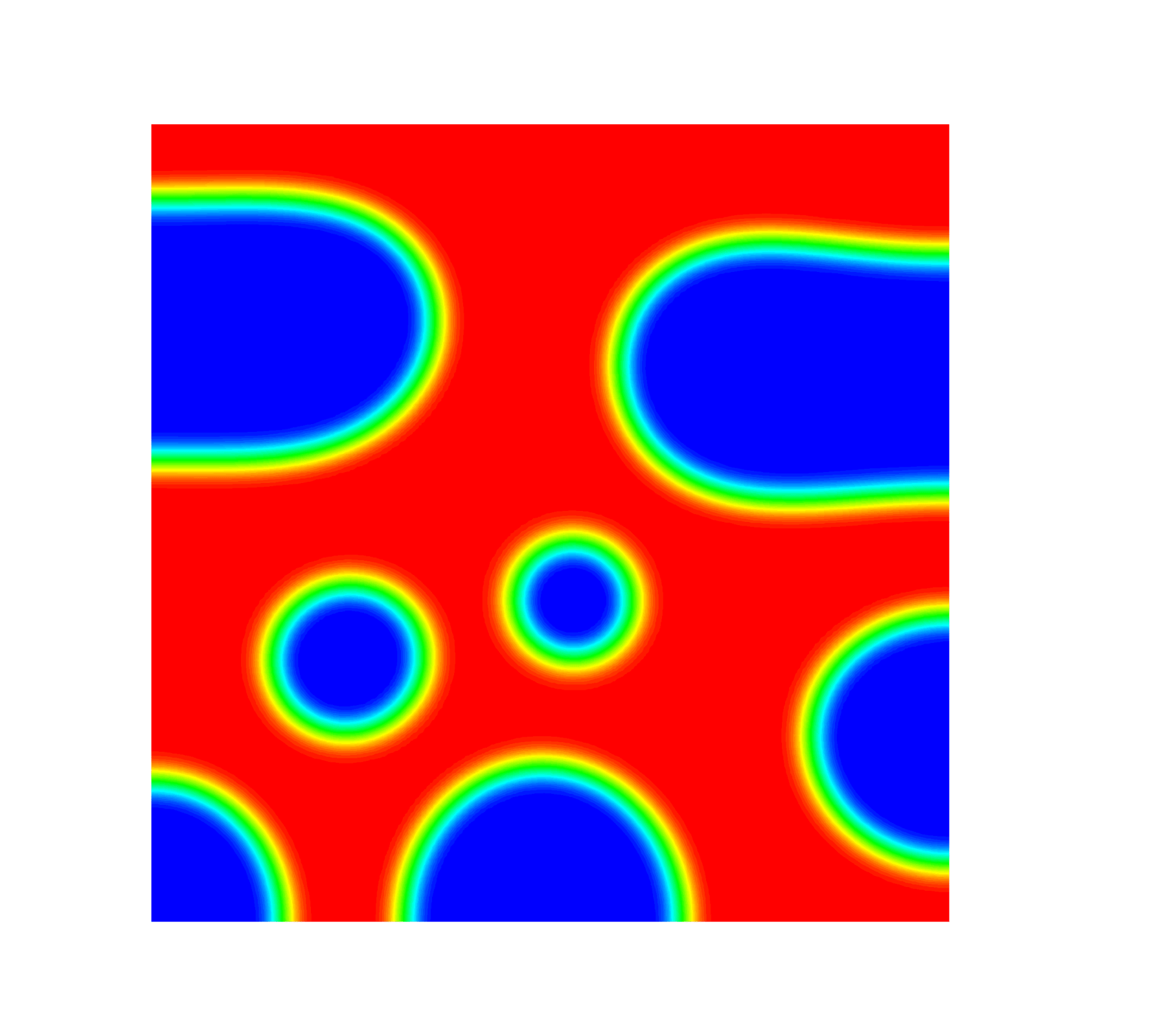}}
		\subfigure[$t=6$] 	{\includegraphics[width=1.6in,trim=70 20 60 20,clip]{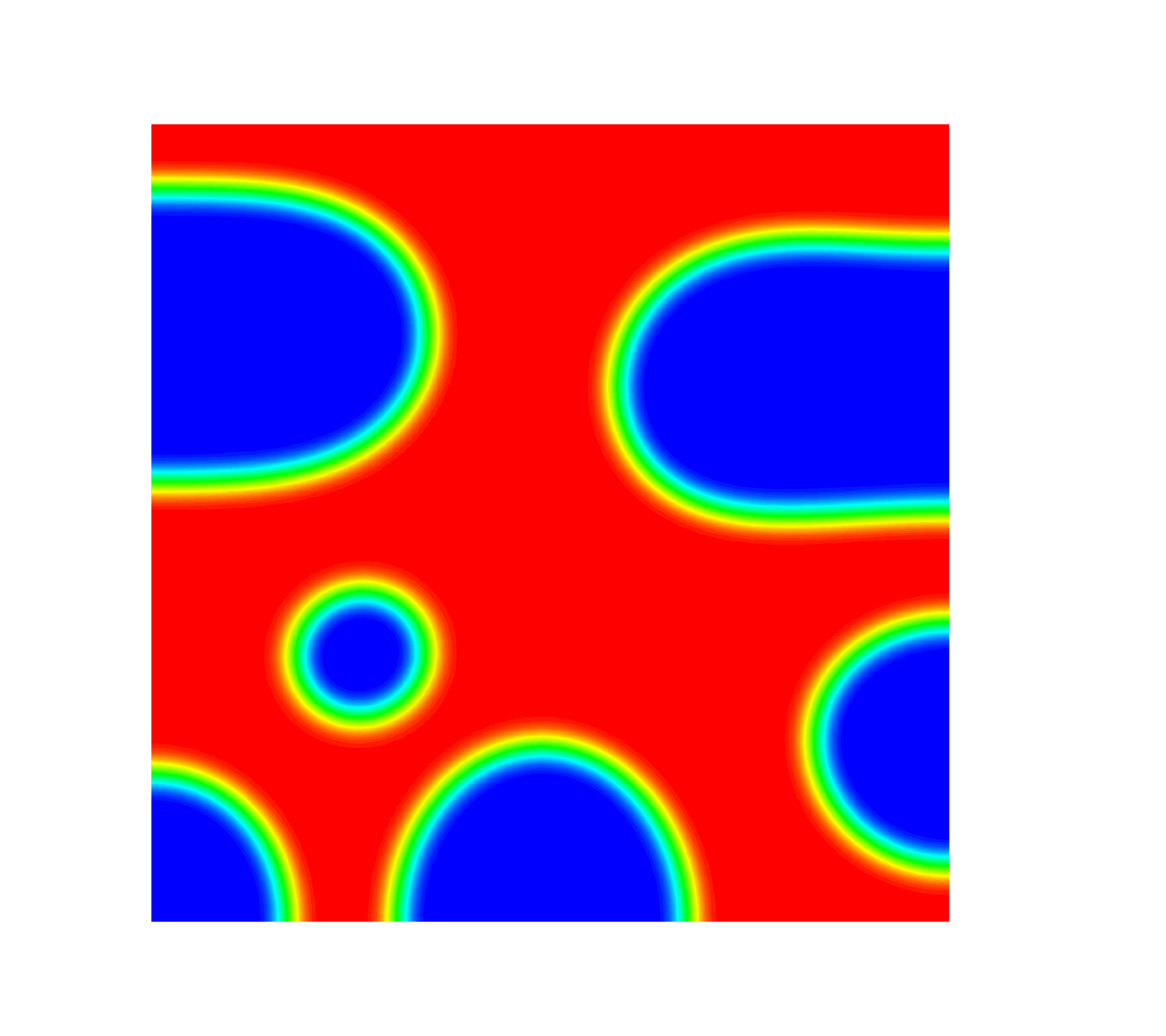}}
		\subfigure[$t=7.1$] {\includegraphics[width=1.6in,trim=70 20 60 20,clip]{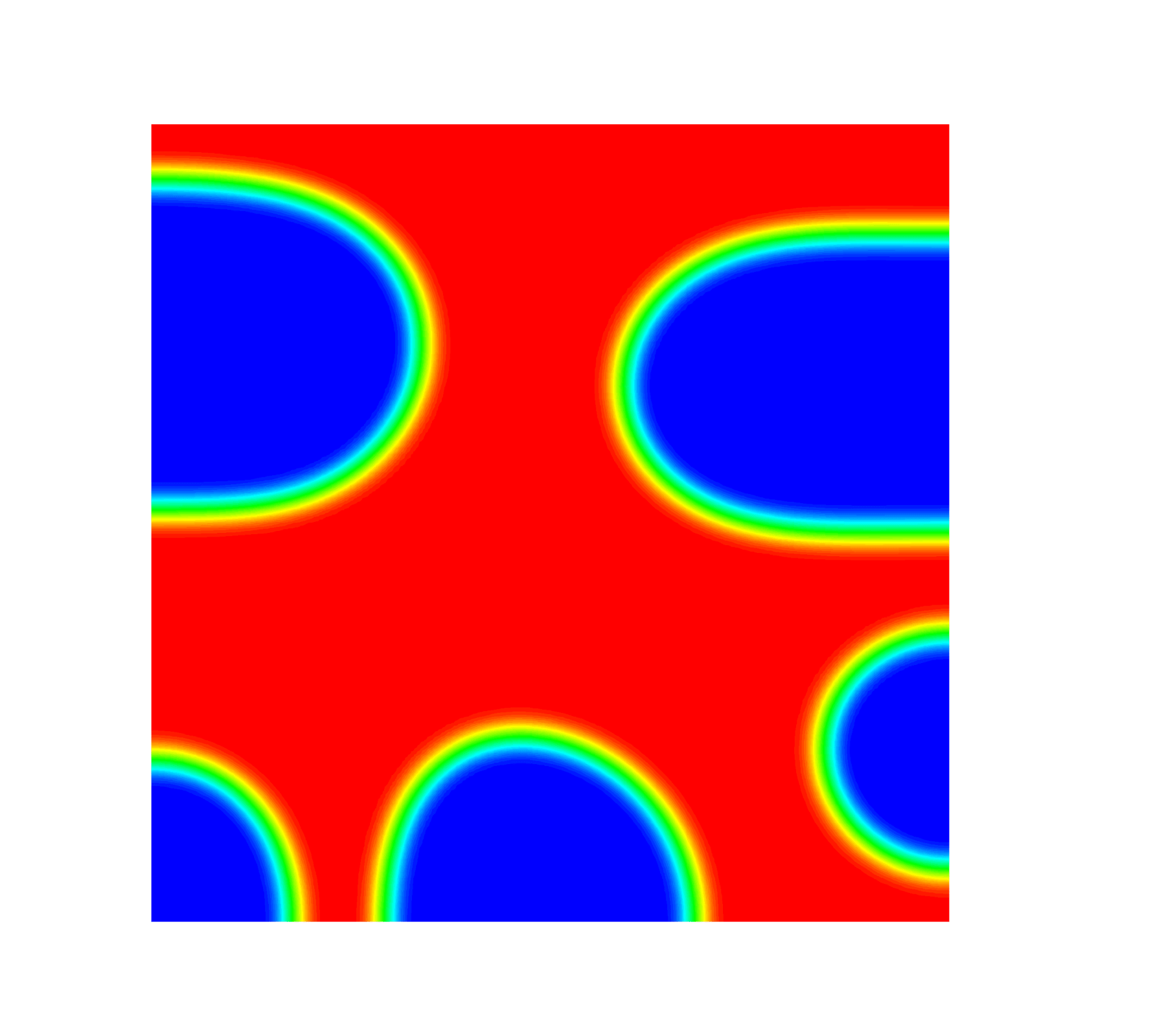}}
		\subfigure[$t=10$] 	{\includegraphics[width=1.6in,trim=70 20 60 20,clip]{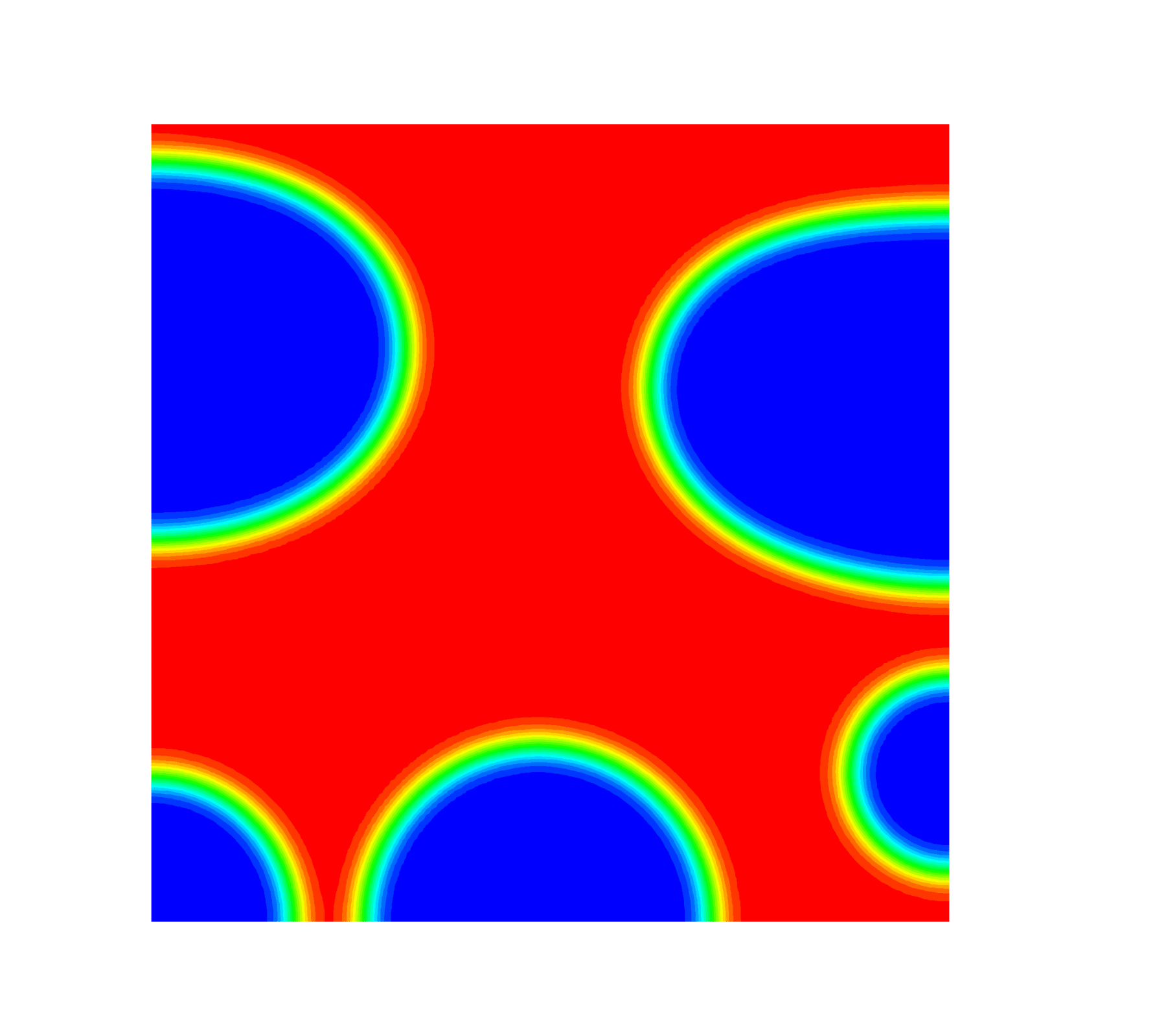}}
	\end{center}
	\caption{The evolution of phase variable for coarsening process.}
	\label{phase_spinodal}
\end{figure}

Figure \ref{phase_spinodal}  shows  the morphological patterns of phase function   during coarsening process.  We observe that the binary mixture undergoes phase transition behaviour, and  eventually forms binary component structure with specific interface, reflecting that the phase diagram gradually reaches a stable equilibrium state by evolving from a homogeneous state.

\begin{figure}[h!]
	\centering
	\begin{tabular}{cc}
		\subfigure[Energy]{\includegraphics[width=0.45\textwidth]{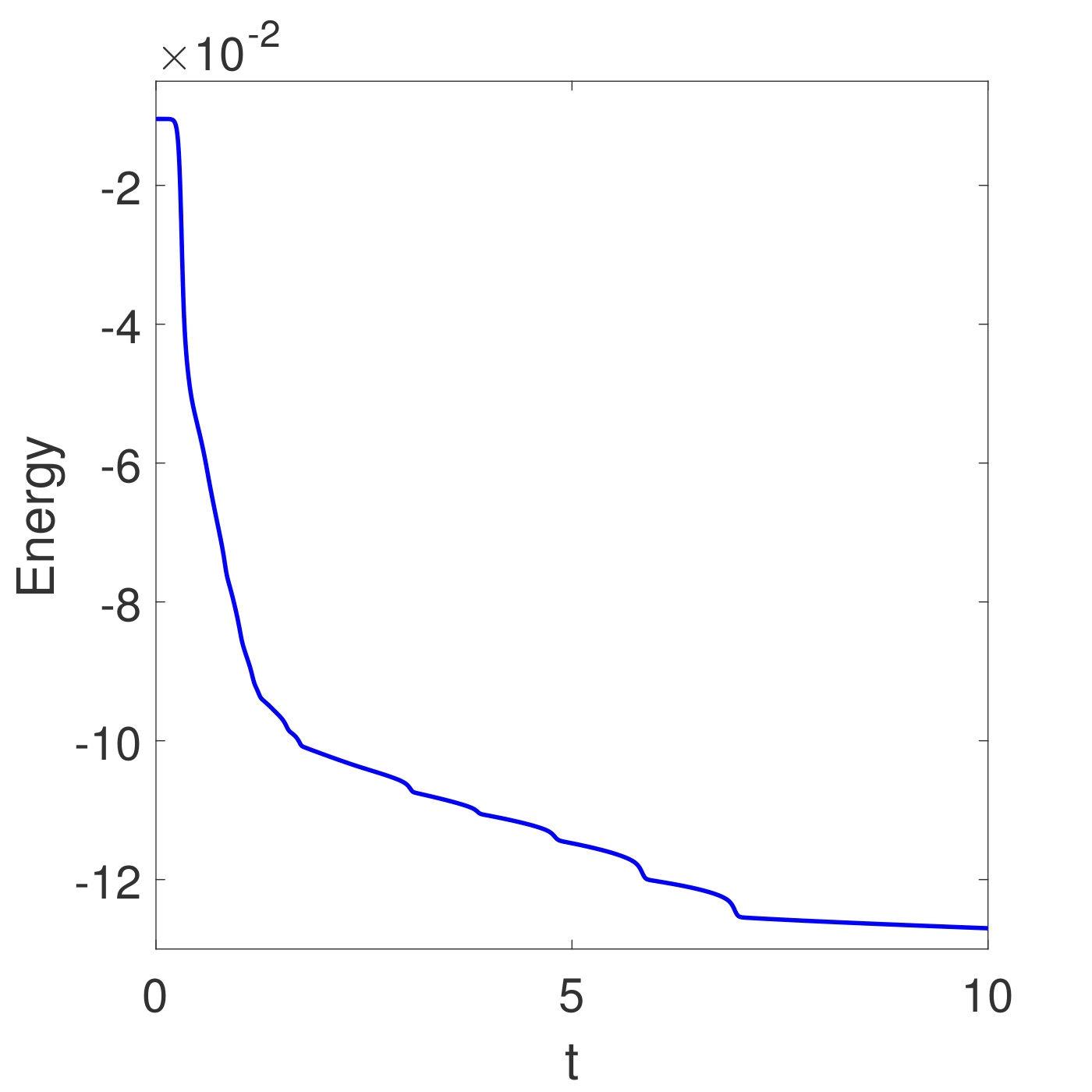}\label{Energy_spinodal}}& 	\subfigure[Mass]{\includegraphics[width=0.45\textwidth]{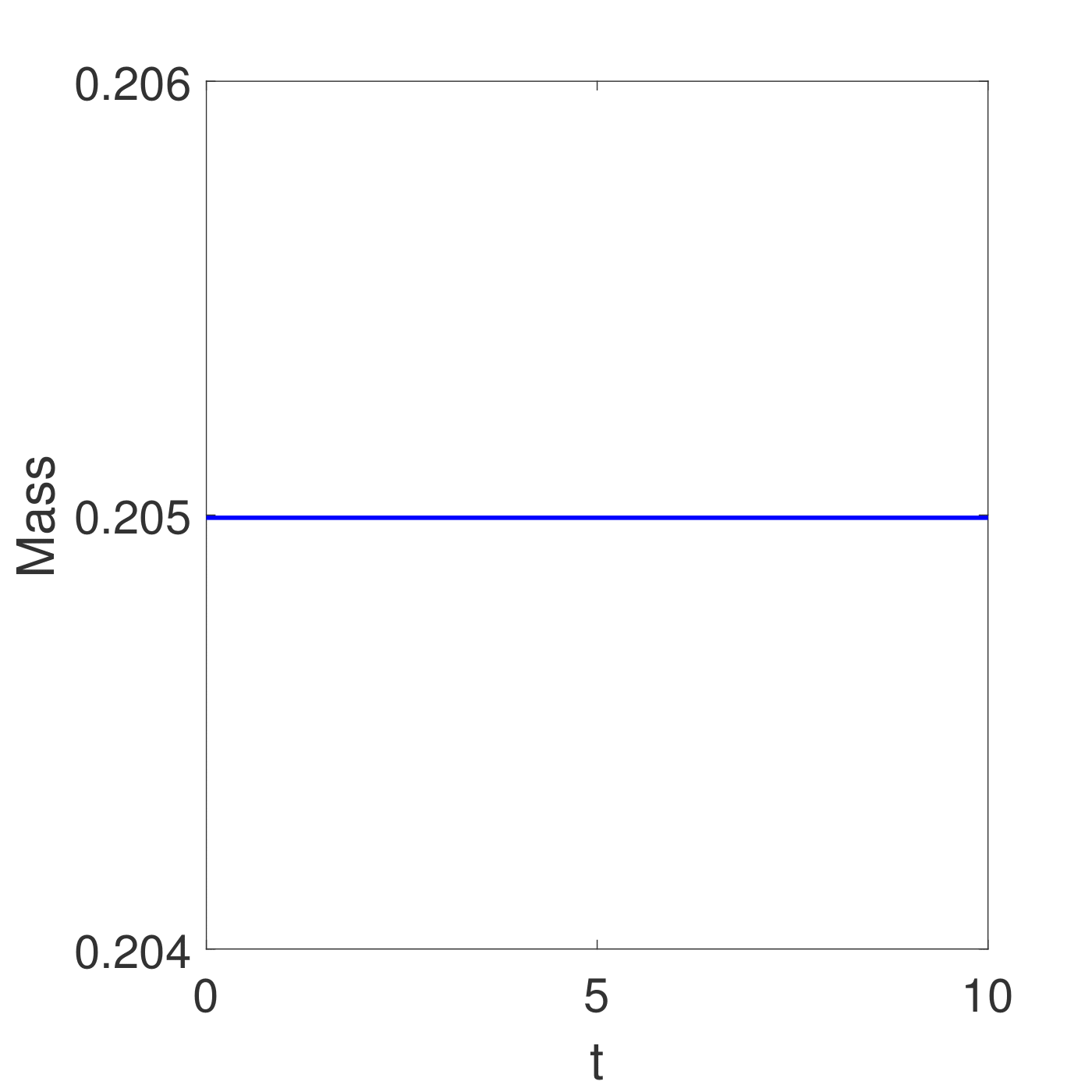}\label{Mass_spinodal}}
	\end{tabular}
	\caption{Evolution of the discrete energy and mass for coarsening process.}
	\label{EnergyMass_spinodal}
\end{figure}

\begin{figure}[h!]
	\centering
	\begin{tabular}{cc}
		\subfigure
		{\includegraphics[width=0.35\textwidth]{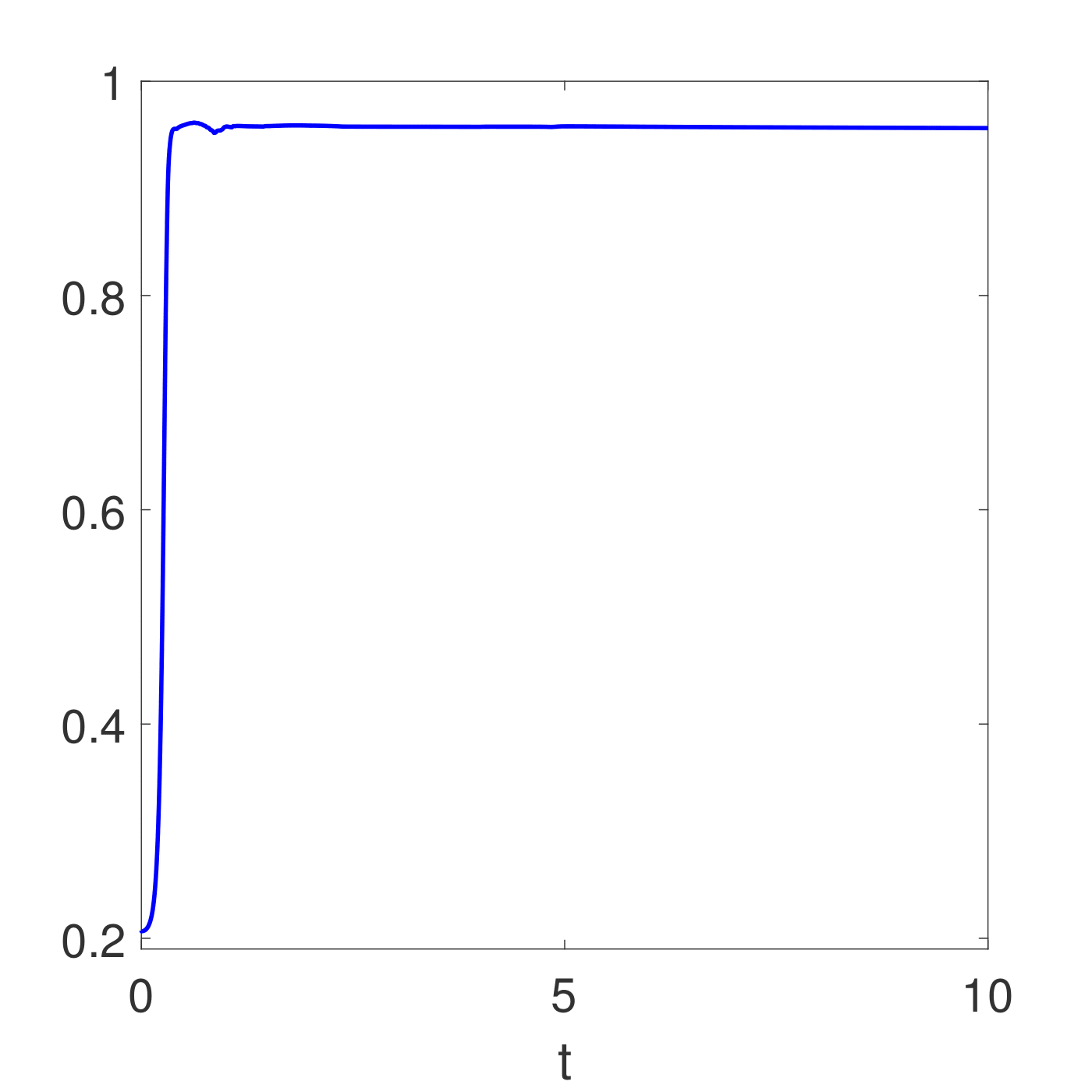}\label{phi_max}}& 	
		\subfigure
		{\includegraphics[width=0.35\textwidth]{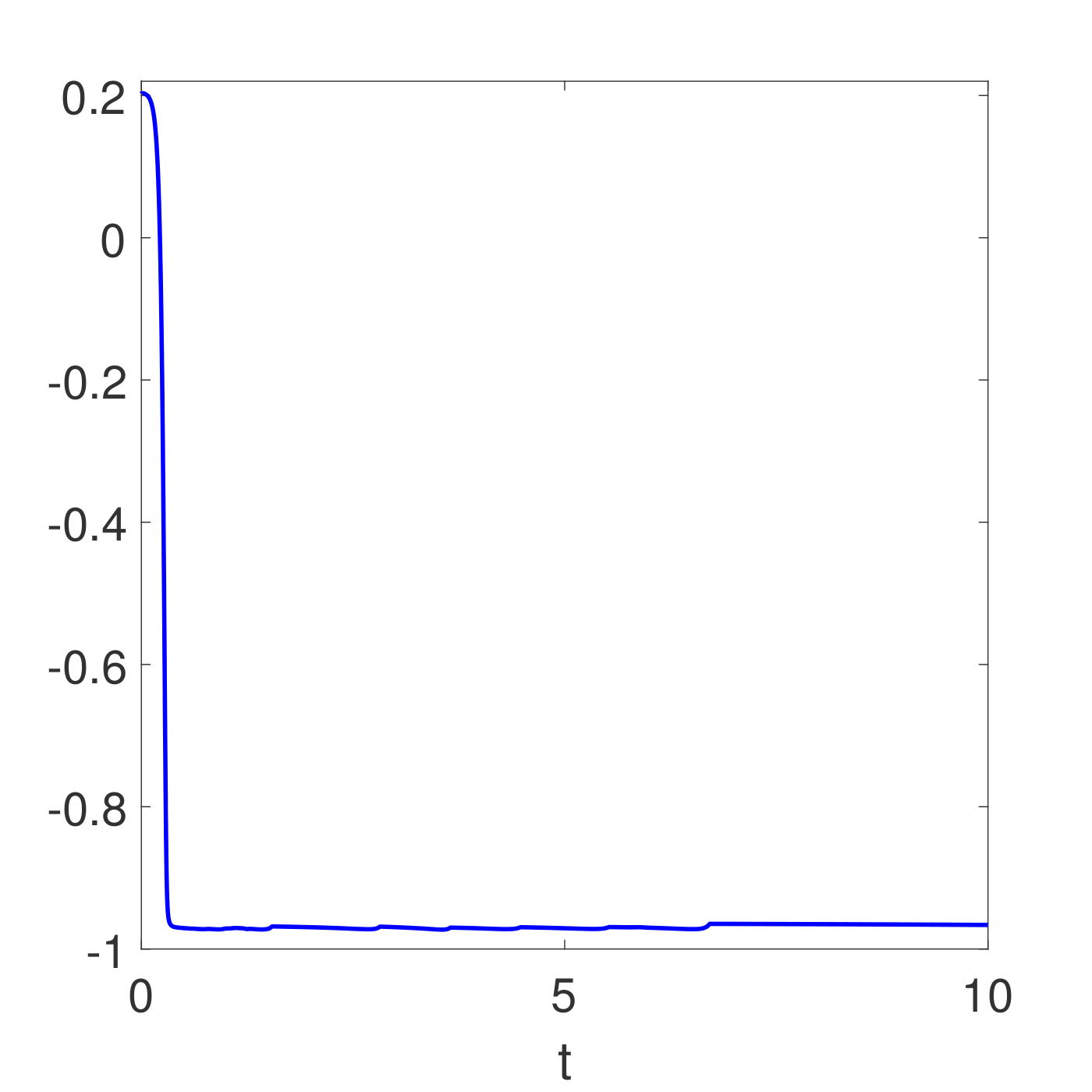}\label{phi_min}}
	\end{tabular}
	\caption{Evolution of the maximum value (left) and mininum value (right) of the phase variable.}
	\label{phi_mm}
\end{figure}

 Figure \ref{Energy_spinodal}  depicts the evolution curves of discrete energy. As expected, the discrete energy  is  indeed   non-increasing with respect to  time  that confirms that our numerical method is unconditionally stable. In order to further illustrate the efficiency of numerical method, we calculate the discrete mass by $\int_{\Omega}\phi_h^nd\bm{x}$ as plotted in Figure  \ref{Mass_spinodal}, in which the  desired mass conservation indicates the mass preserving of developed numerical method. 
 Figure \ref{phi_mm} presents the evolution of maximum and minimum of numerical solution for phase variable. It is observed
 that the maximum value is about $0.9613$, while the minimum value is $-0.9728$. The numerical value of $\phi$ lies
 completely within interval $(-1,1)$ guaranteeing the positivity of $\ln(1+\phi)$ and $\ln(1-\phi)$, which provides the
 evidence of the robustness of numerical method from another perspective, and is consistent well with the theoretical analysis in Theorem \ref{USS}.
 These dynamics of phase diagram show statistically similar features in the numerical solutions obtained in \cite{WChen_JJing_CWang_XWang_JSC_2022}. The reasonable dynamical
 behaviors of phase separation confirms the  suitability and effectiveness of numerical algorithm  presented.

\subsection{Rotational flow}
\begin{figure}[ht]
	\begin{center}
		\subfigure[$t=0$] {\includegraphics[width=1.6in,trim=70 20 60 20,clip]{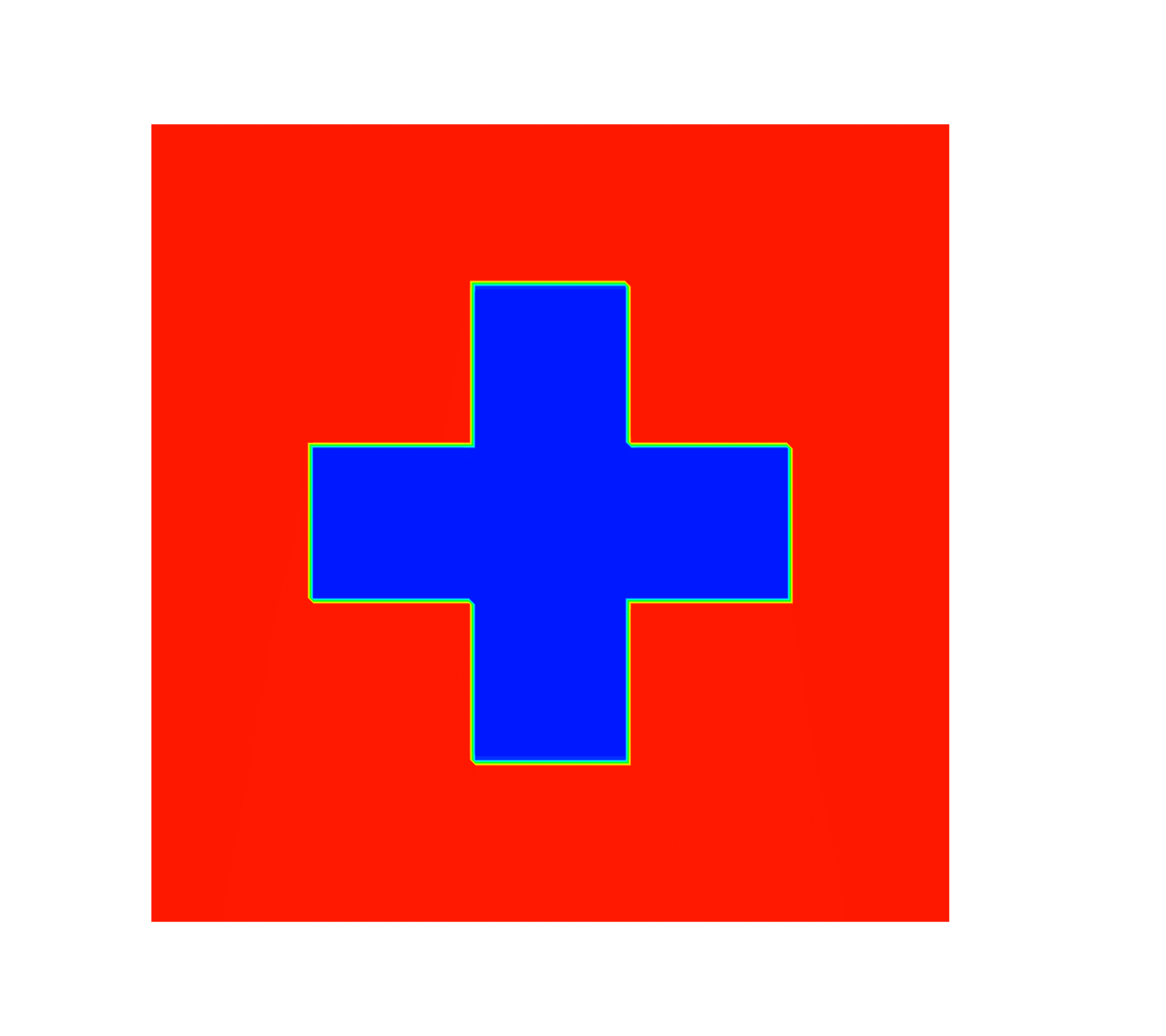}\label{initial_cross}}
		\subfigure[$t=0.2$] 	{\includegraphics[width=1.6in,trim=70 20 60 20,clip]{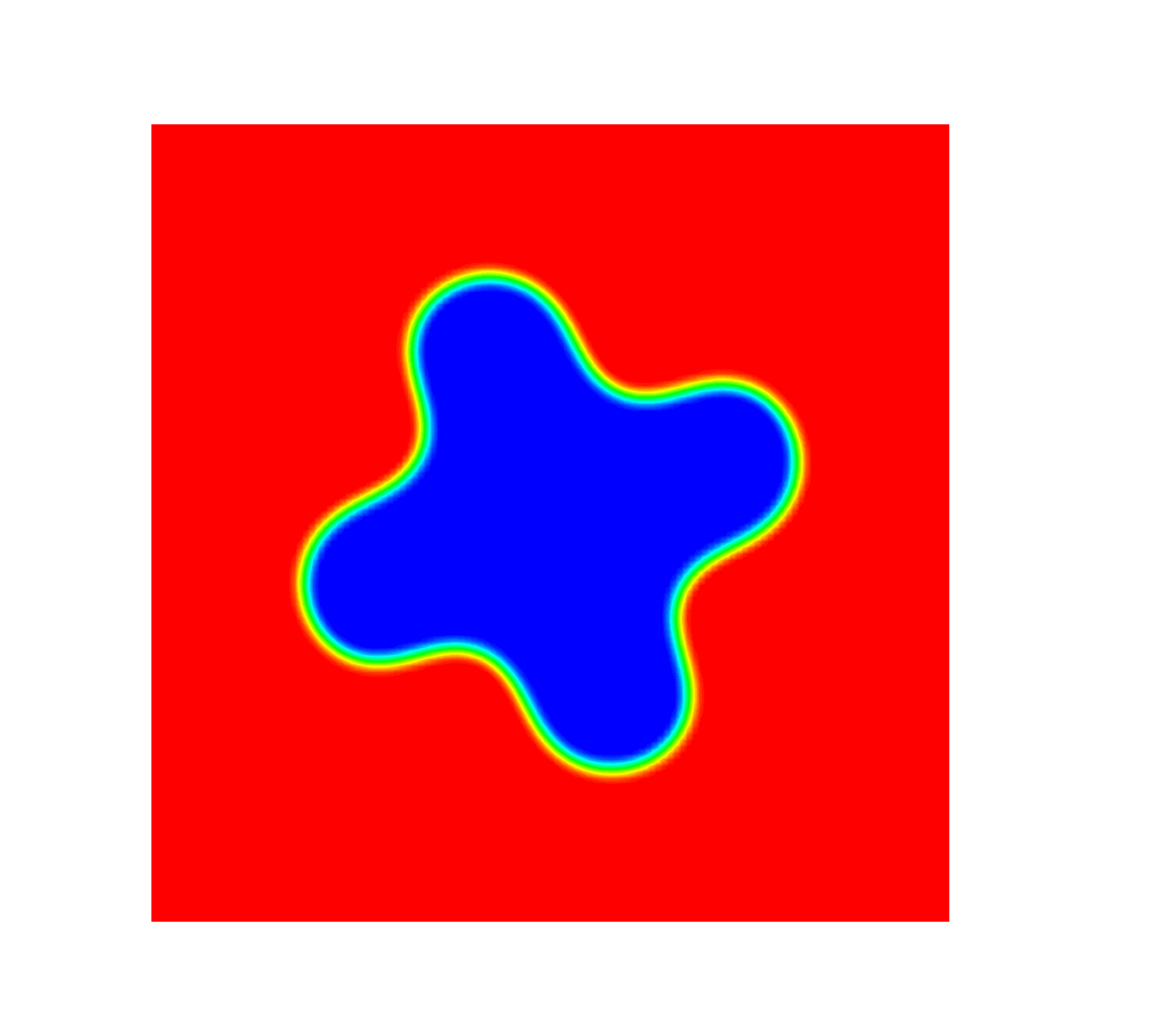}}
		\subfigure[$t=0.35$] 	{\includegraphics[width=1.6in,trim=70 20 60 20,clip]{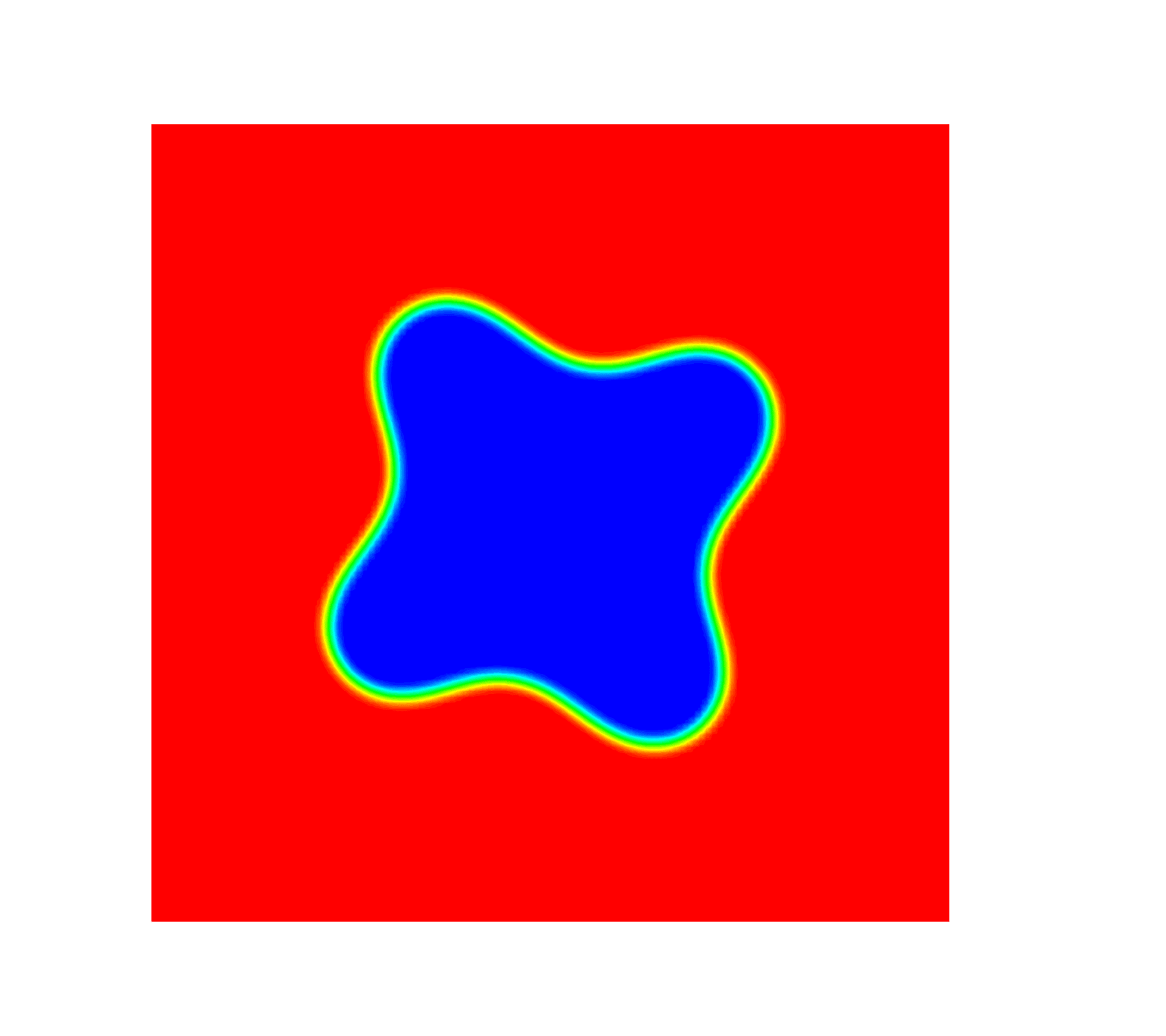}}
		\subfigure[$t=0.5$] 	{\includegraphics[width=1.6in,trim=70 20 60 20,clip]{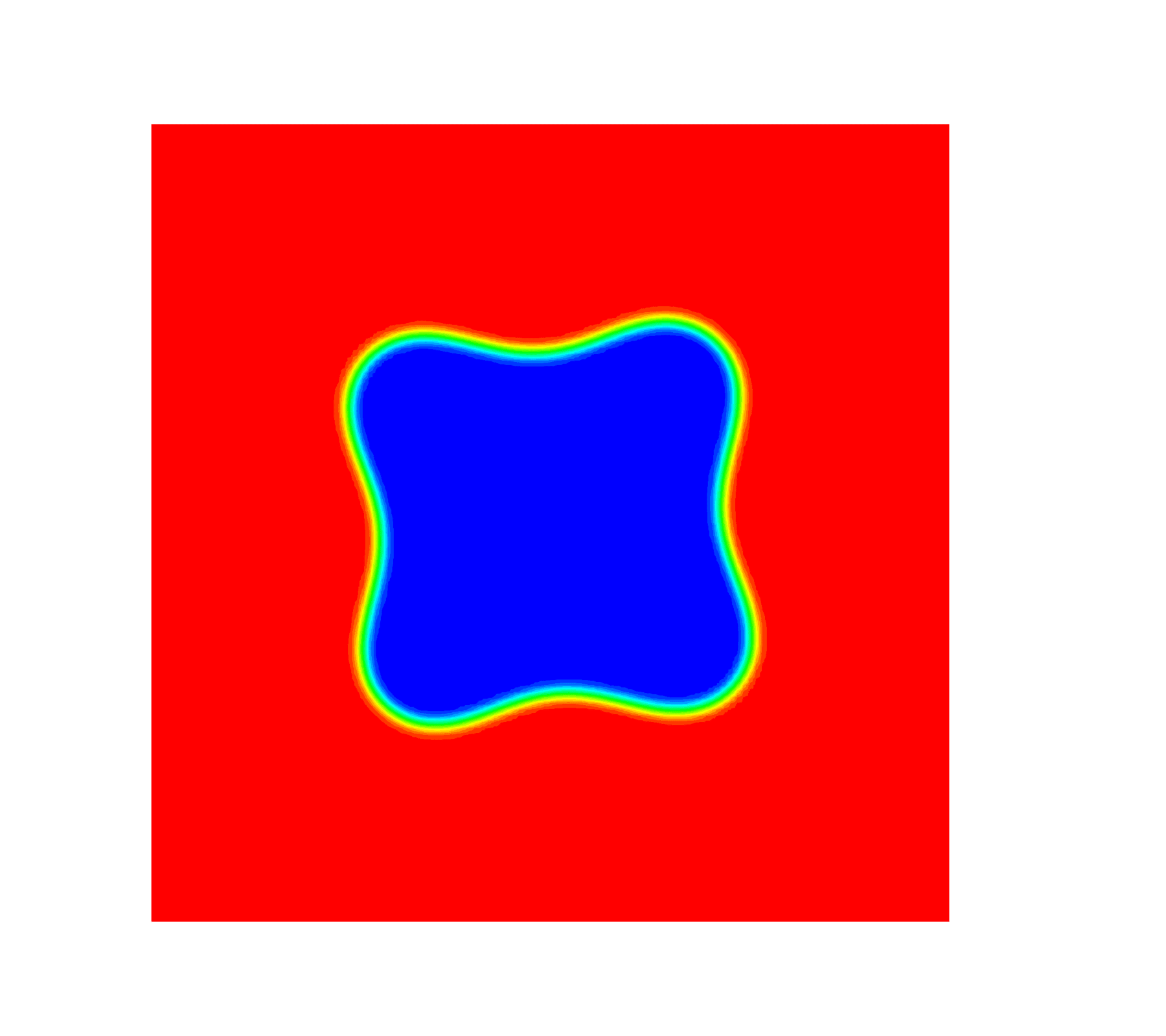}}
		\subfigure[$t=1.0$] 	{\includegraphics[width=1.6in,trim=70 20 60 20,clip]{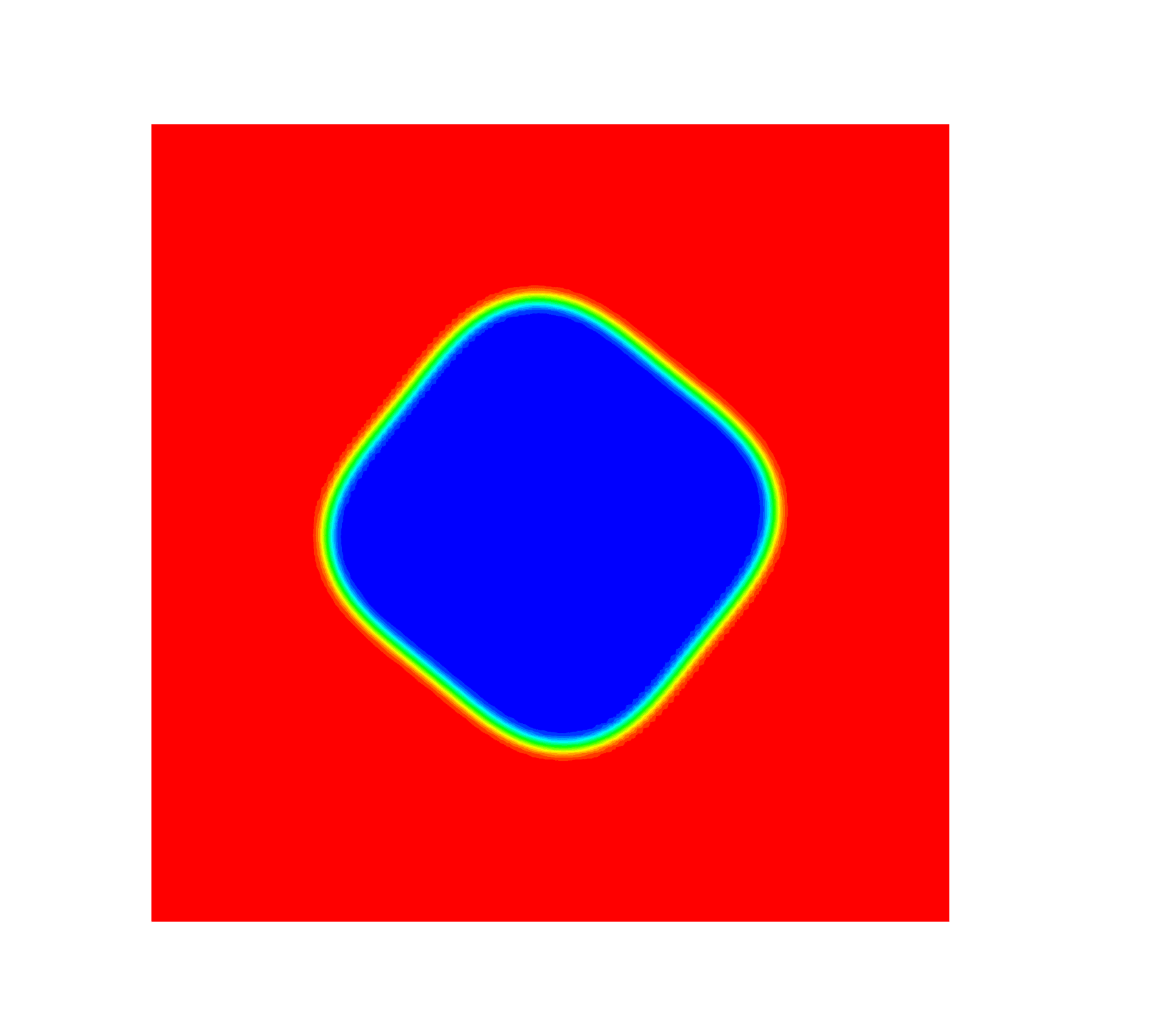}}
		\subfigure[$t=1.25$] 	{\includegraphics[width=1.6in,trim=70 20 60 20,clip]{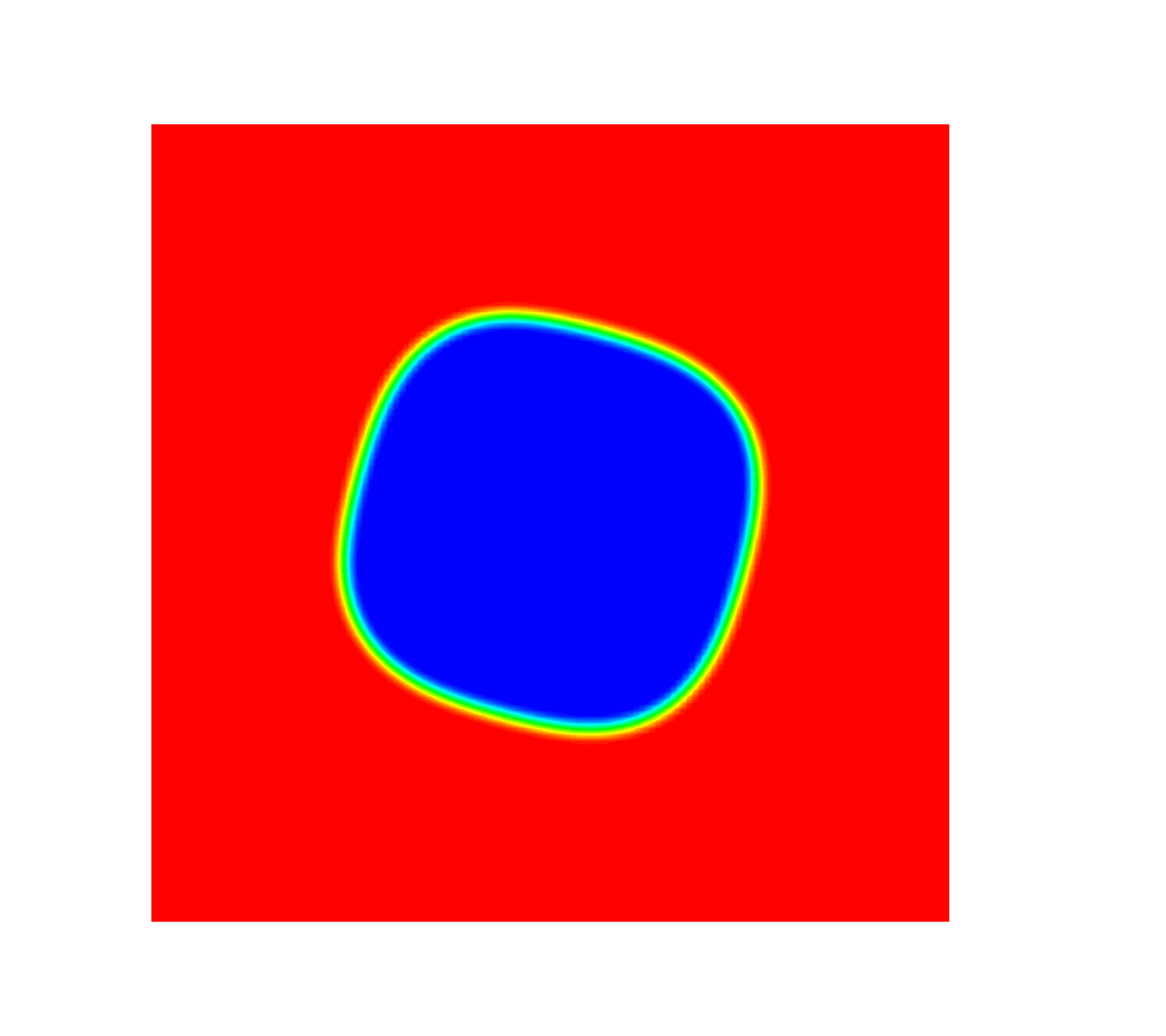}}
		\subfigure[$t=1.5$] 	{\includegraphics[width=1.6in,trim=70 20 60 20,clip]{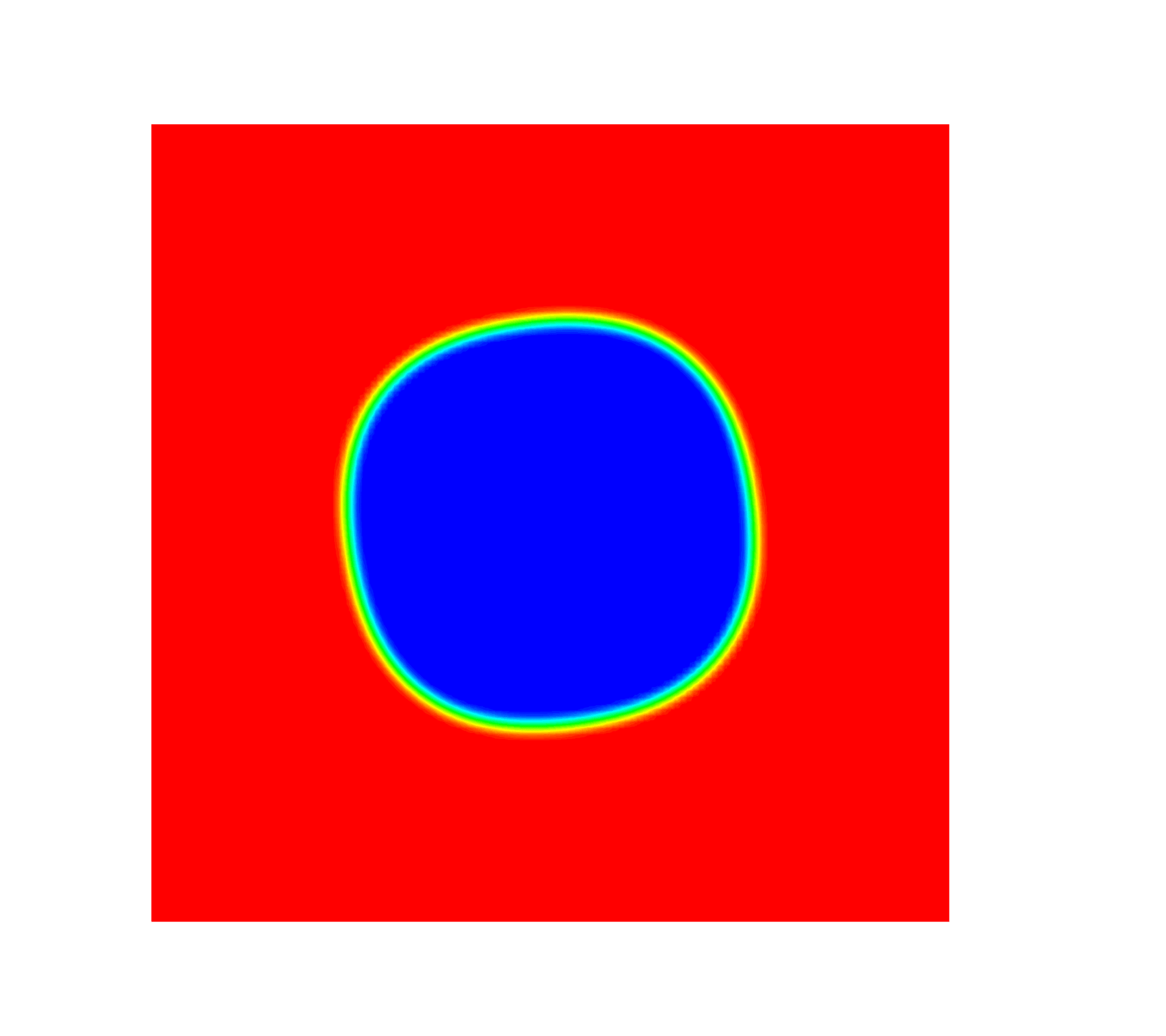}}
		\subfigure[$t=3.0$] 	{\includegraphics[width=1.6in,trim=70 20 60 20,clip]{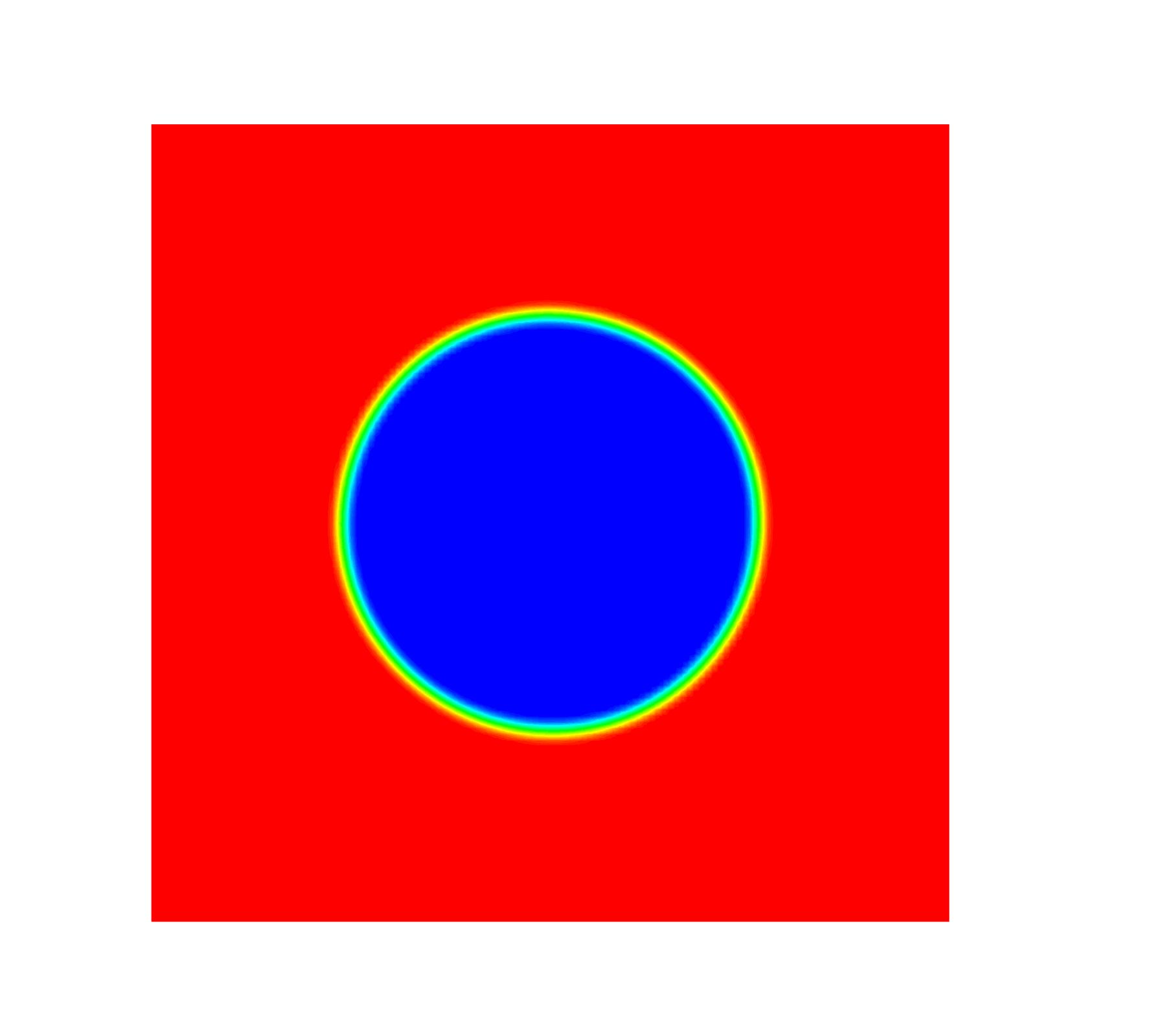}}
	\end{center}
	\caption{Snapshots of phase variable for cross-shaped relaxation.}
	\label{rotaional_cross}
\end{figure}

In this test, we simulate effect of rotational external forces on the shape relaxation for a cross-shaped droplet into a circle. The computational domain is $\Omega=[0,1]^2$ with the rotational boundary conditions $\Au=(1-2y,2x-1)$, and the parameters are  $\sw=1$, $Pe=200$, $Re=1$, $\eta=1$, $M=1$ and $\epsilon=0.01$. The initial phase variable is set to $\phi=0.9$ in the polygon area indicated by the blue color and $\phi=-0.9$ in the remainder area indicated by the red color, which is plotted in Figure \ref{initial_cross}.

\begin{figure}[ht]
	\begin{center}
		\subfigure {\includegraphics[width=1.7in,trim=70 20 30 20,clip]{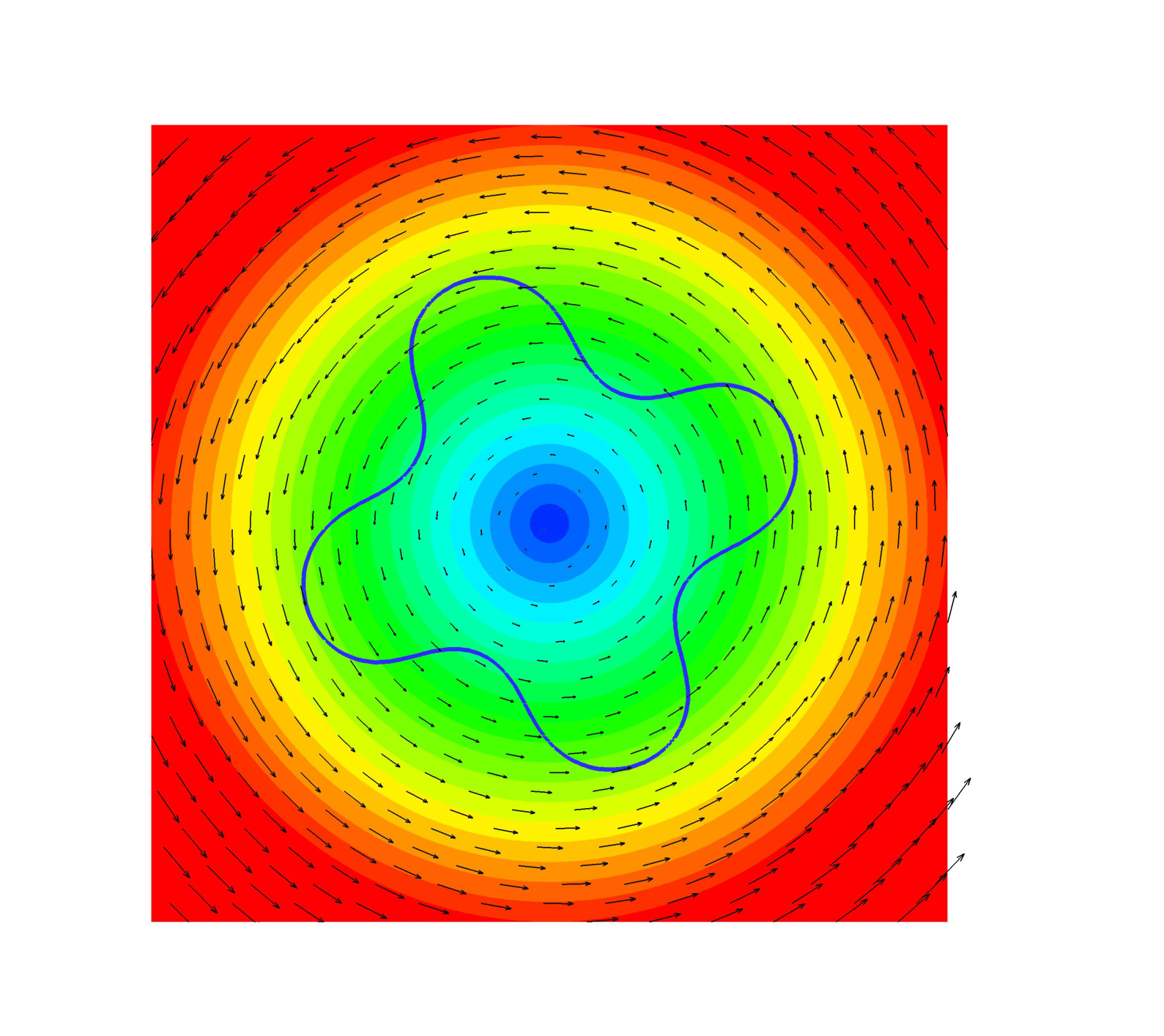}}\hskip  -3mm
		\subfigure{\includegraphics[width=1.7in,trim=70 20 30 20,clip]{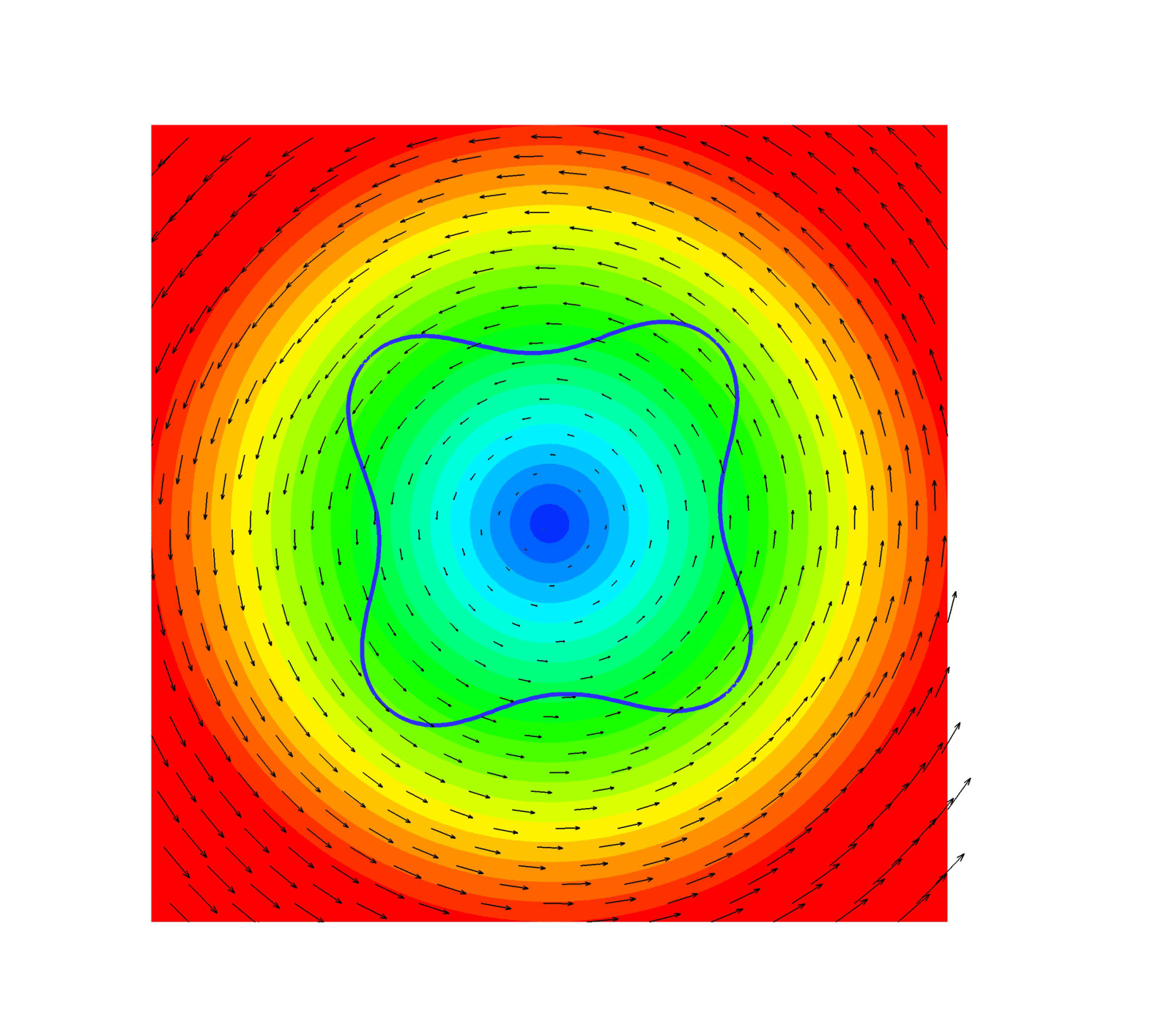}}\hskip  -3mm
		\subfigure{\includegraphics[width=1.7in,trim=70 20 30 20,clip]{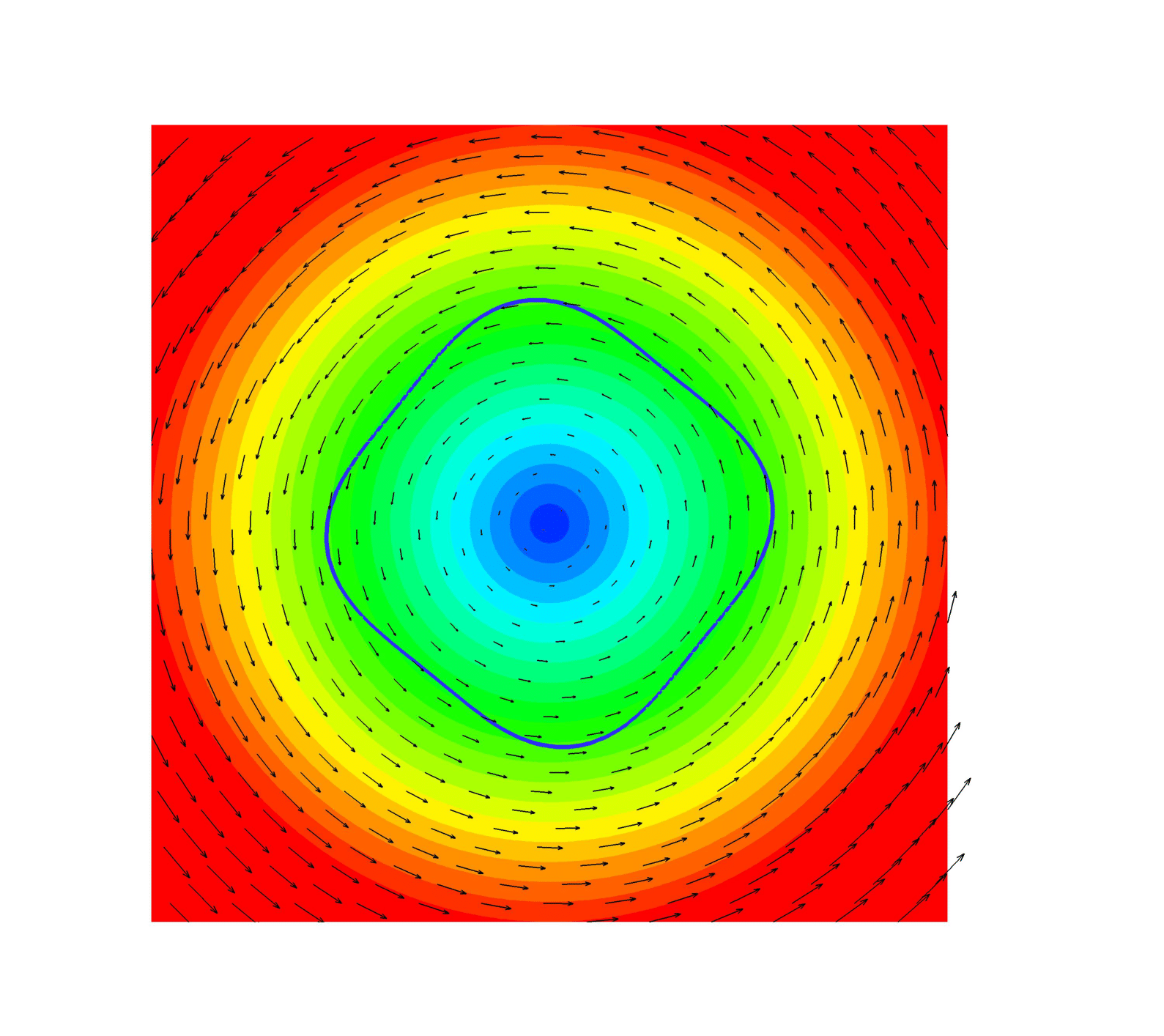}}\hskip  -3mm
		\subfigure{\includegraphics[width=1.7in,trim=70 20 30 20,clip]{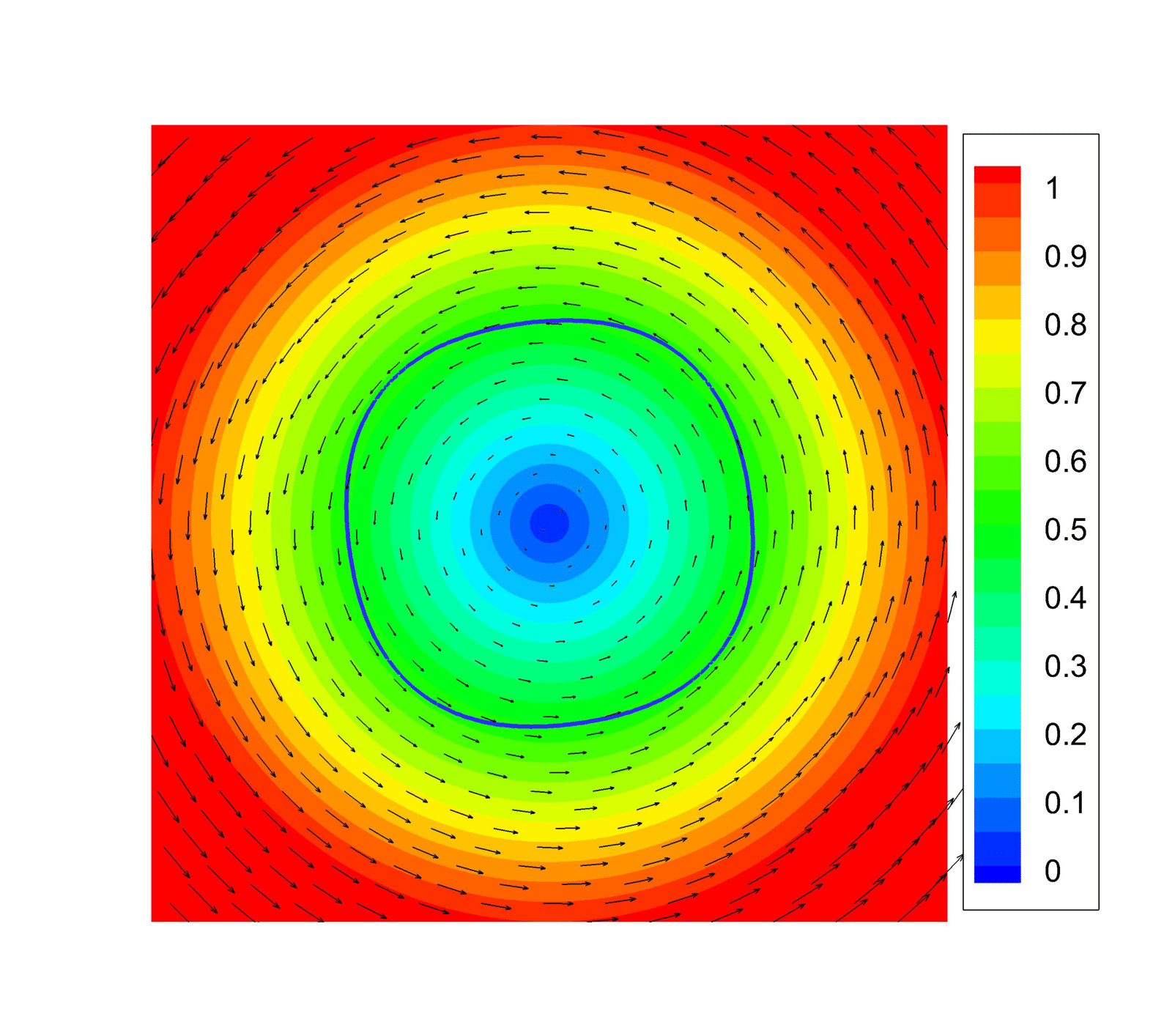}}
	\end{center}
	\caption{Velocity field  for cross-shaped relaxation at time $t=0.2, 0.5, 1.0, 1.5$.}
	\label{Velocity_cross}
\end{figure}

 The dynamical morphotype of phase variable  are depicted in  Figure \ref{rotaional_cross}.
 The characterized    velocity field are graphed in Figure \ref{Velocity_cross} at different time $t=0.2, 0.5, 1.0, 1.5$. 
Due to the effect of surface tension and imposed rotational forces on the boundary,  the isolated cross shape relaxes and gradually changes to a circular shape as  the flow moves counterclockwise. 

\newpage

\subsection{Lid-driven cavity flow}
In this test, we investigate the lid-driven cavity with an initial data for the half full concentration with each phase on unit domain  
The parameters are picked as  $Re=500$, $\eta=1$,  $\epsilon=0.01$, $\sw=5\times 10^5$, $Pe=100$, $M=0.5$. A uniformed grid $128 \times 128$ is used.  We pose the no-slip boundary    conditions for velocity expect for the upper boundary conditions of  $\Au=[16x^2(x-1)^2,0]^T$. 

\begin{figure}[ht]
	\begin{center}
		\subfigure[$t=0.0$] 	{\includegraphics[width=1.7in,trim=70 0 60 20,clip]{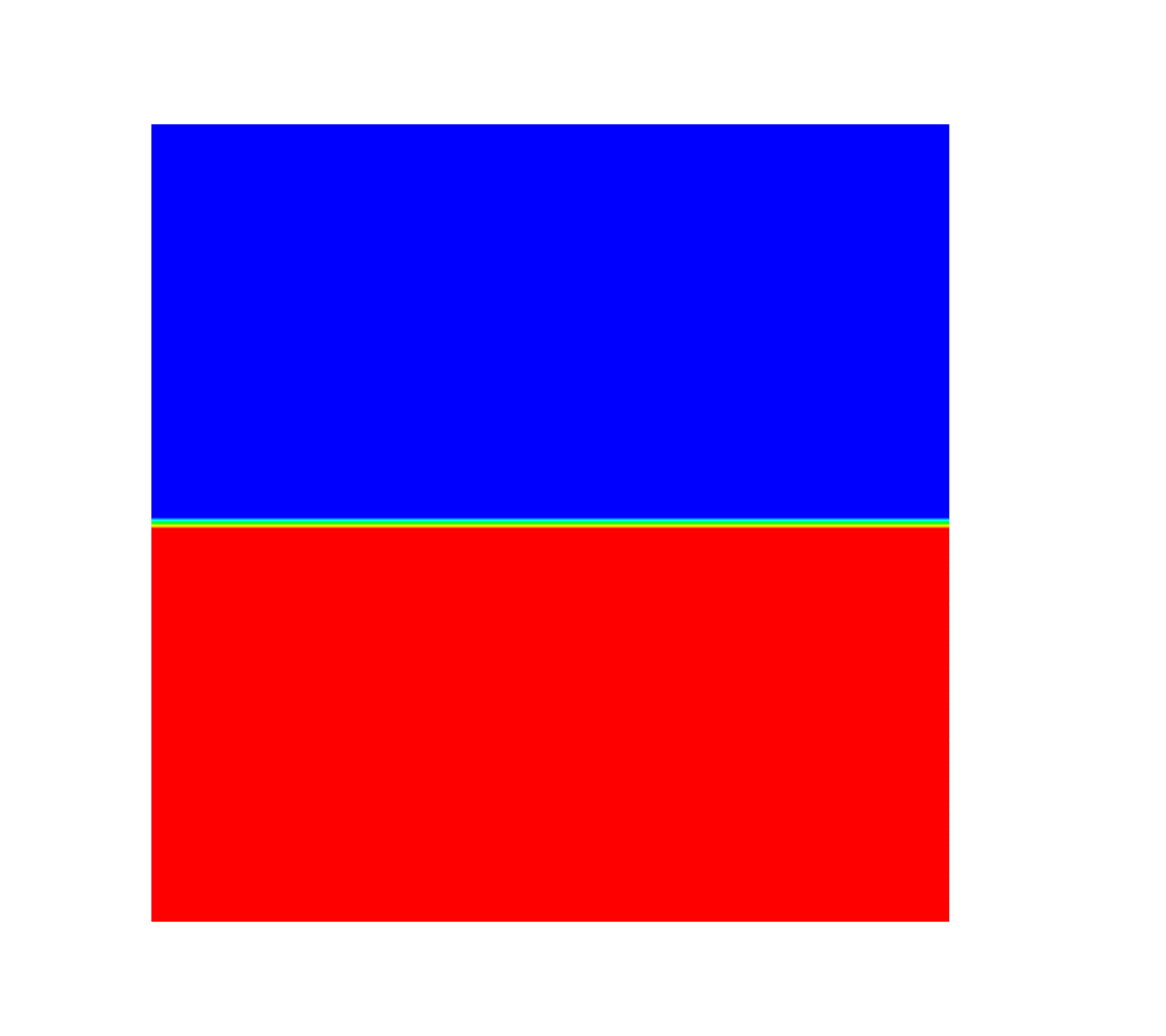}}\hskip -3mm 
		\subfigure[$t=5.0$] 	{\includegraphics[width=1.7in,trim=70 0 60 20,clip]{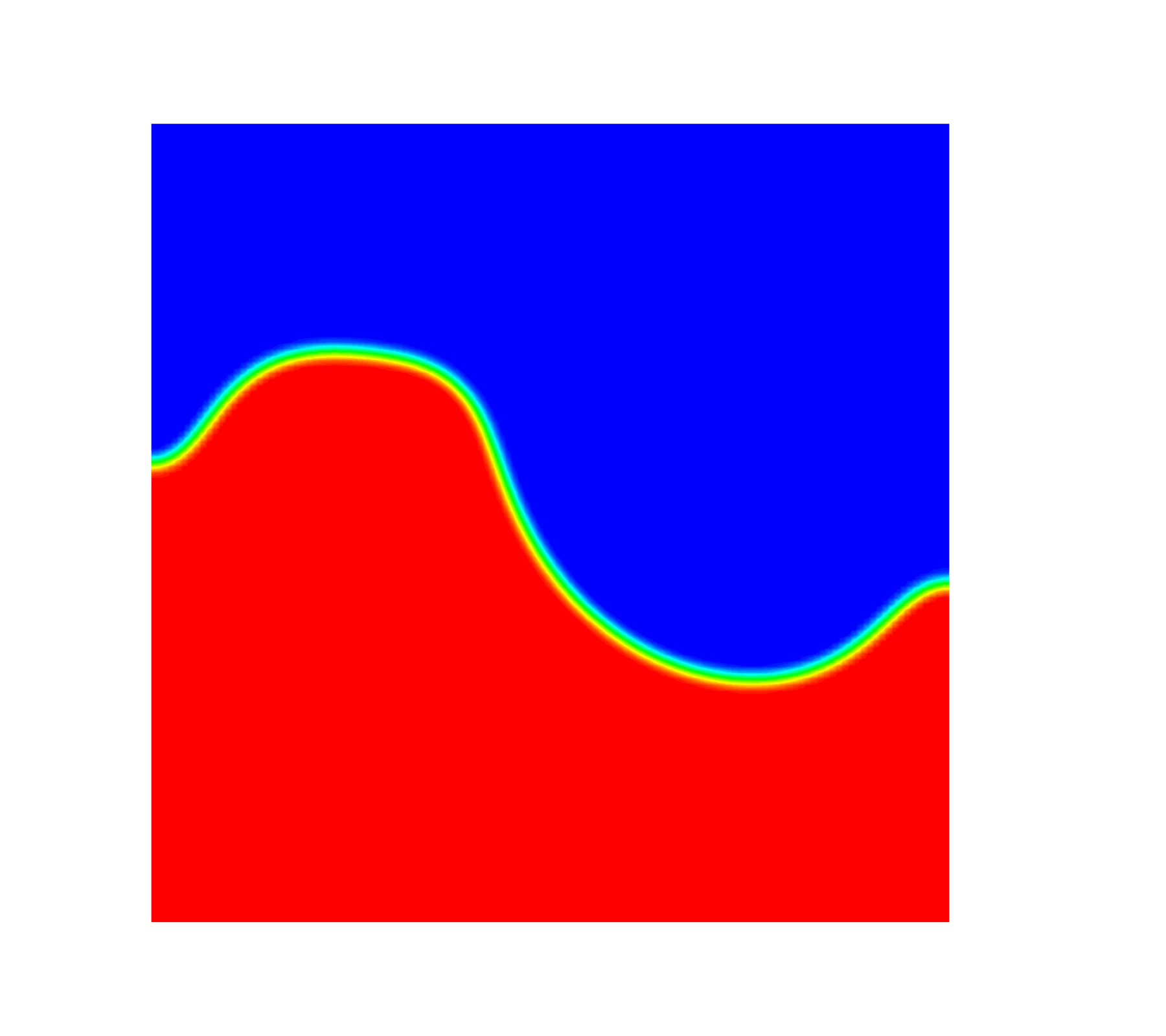}}\hskip -3mm 
		\subfigure[$t=7.0$] 	{\includegraphics[width=1.7in,trim=70 0 60 20,clip]{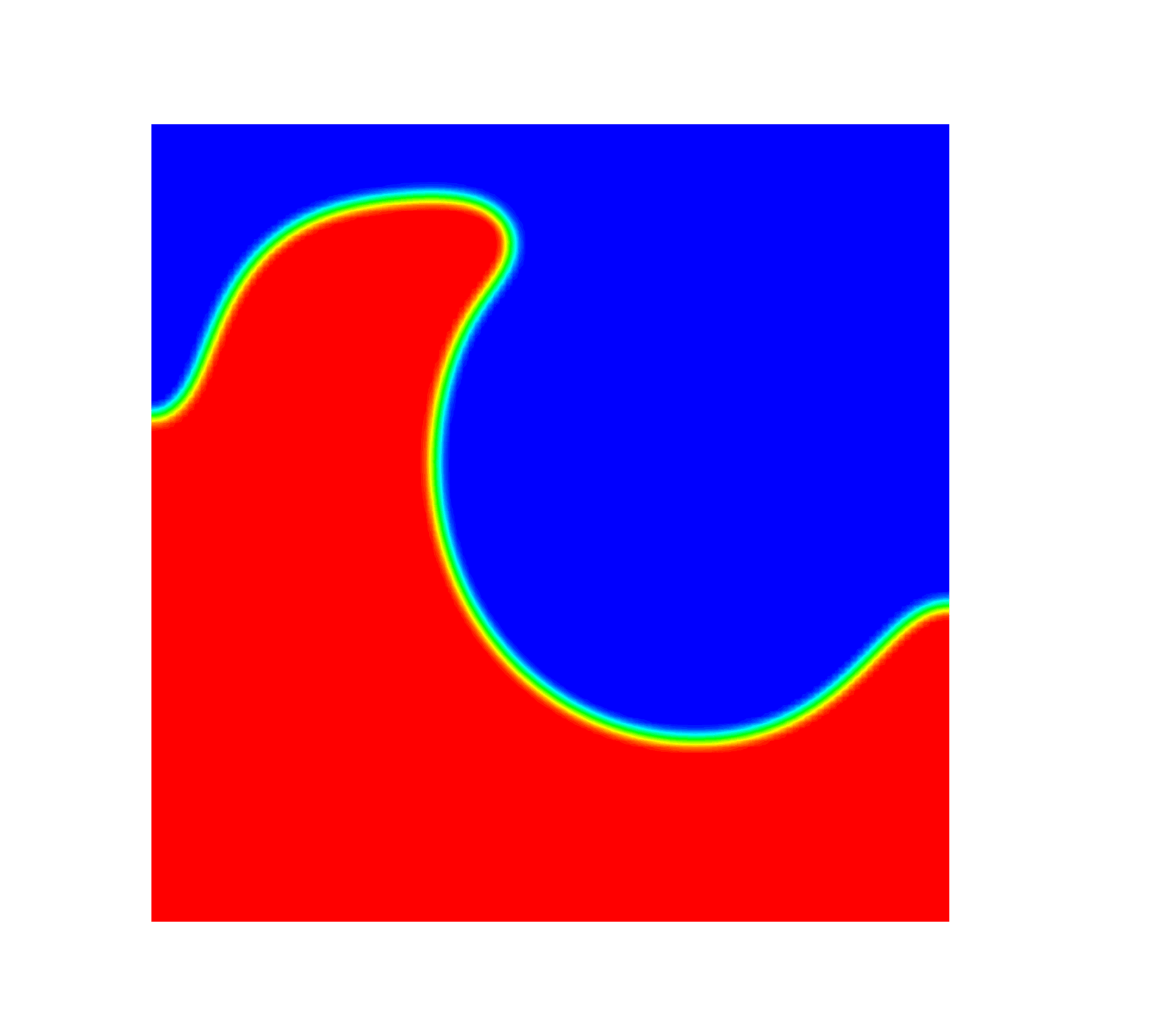}}\hskip -3mm 
		\subfigure[$t=8.0$] 	{\includegraphics[width=1.7in,trim=70 0 60 20,clip]{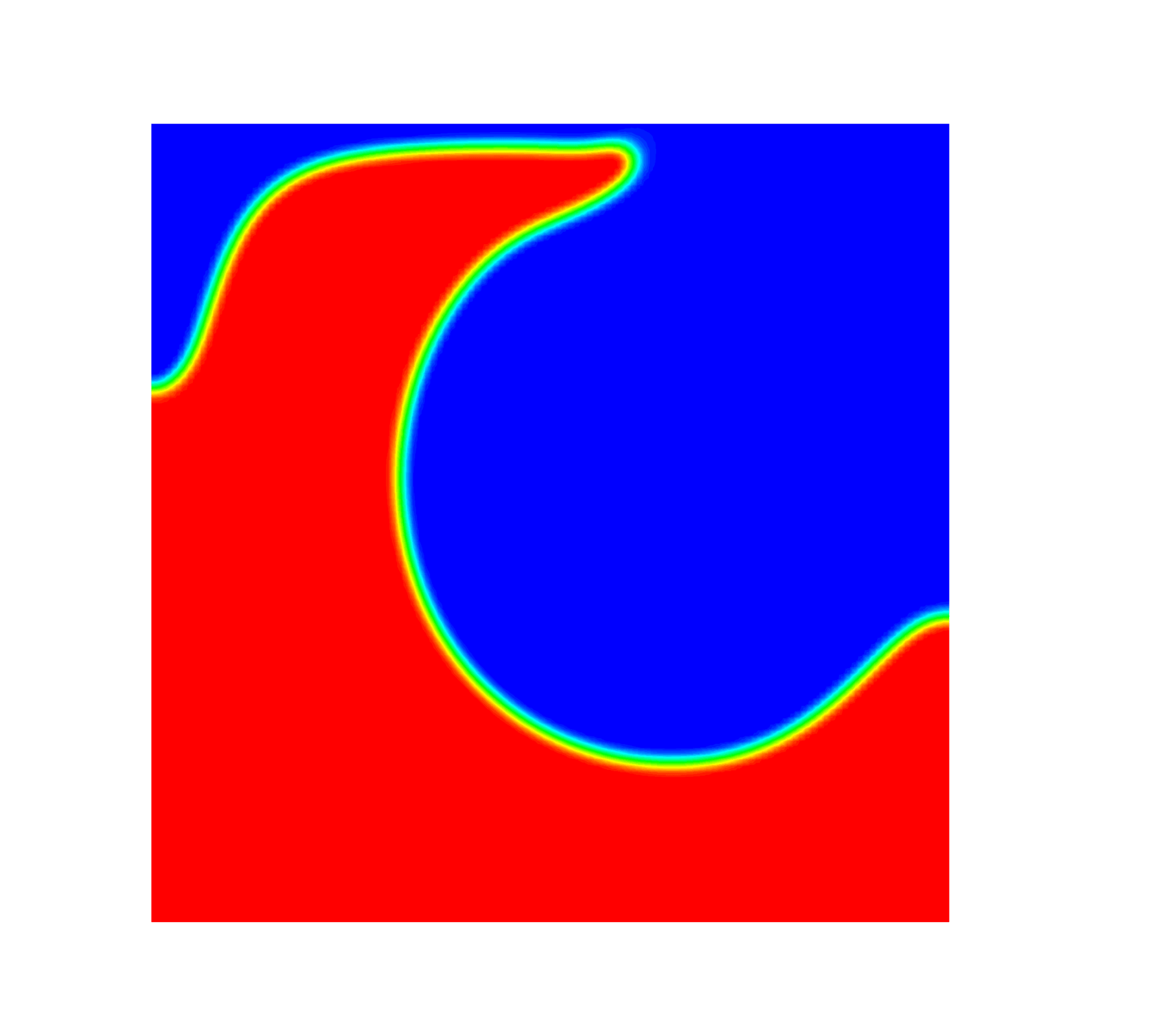}} 
		\subfigure[$t=9.0$] 	{\includegraphics[width=1.7in,trim=70 0 60 20,clip]{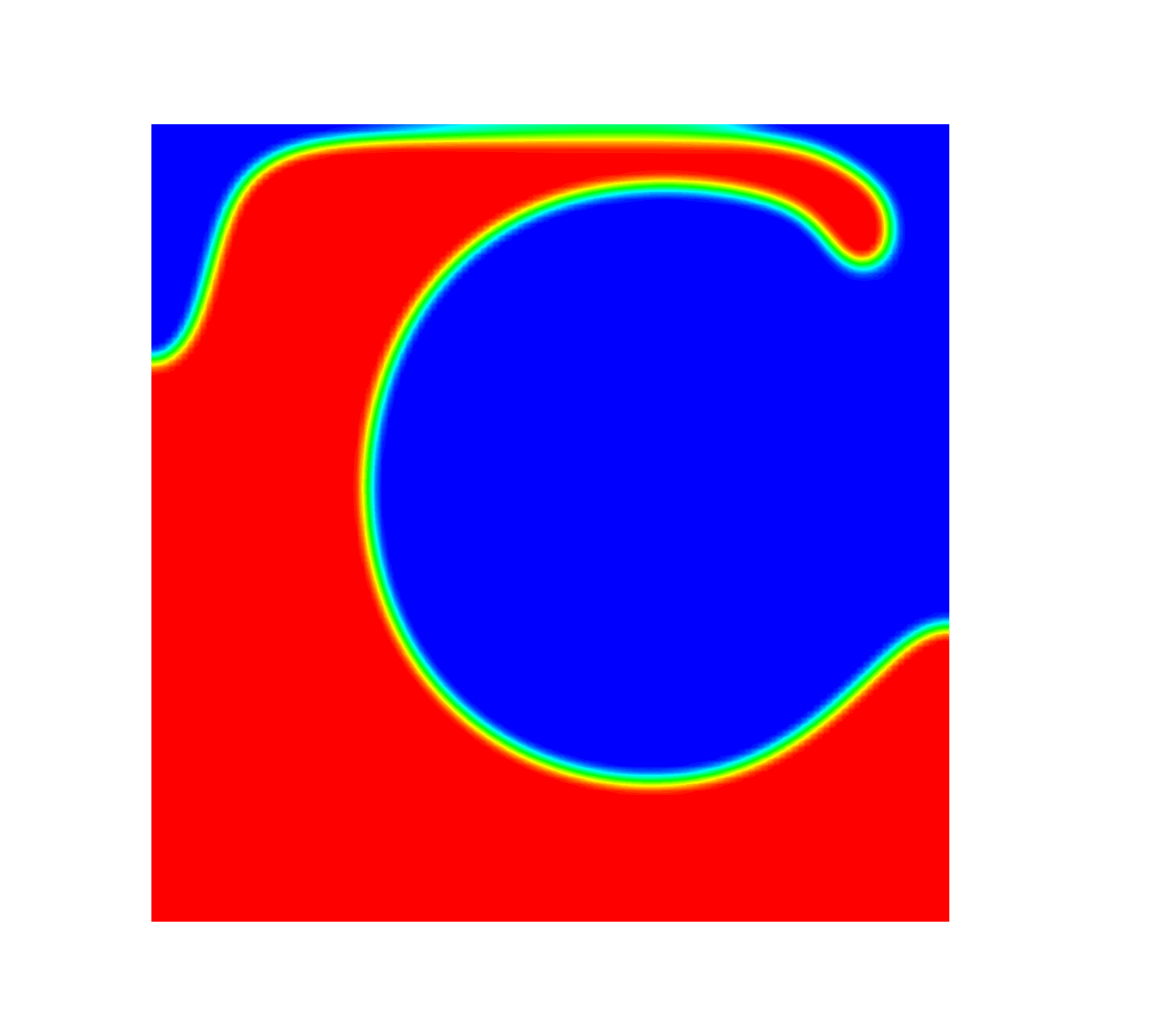}}\hskip -3mm 
		\subfigure[$t=10.0$] 	{\includegraphics[width=1.7in,trim=70 0 60 20,clip]{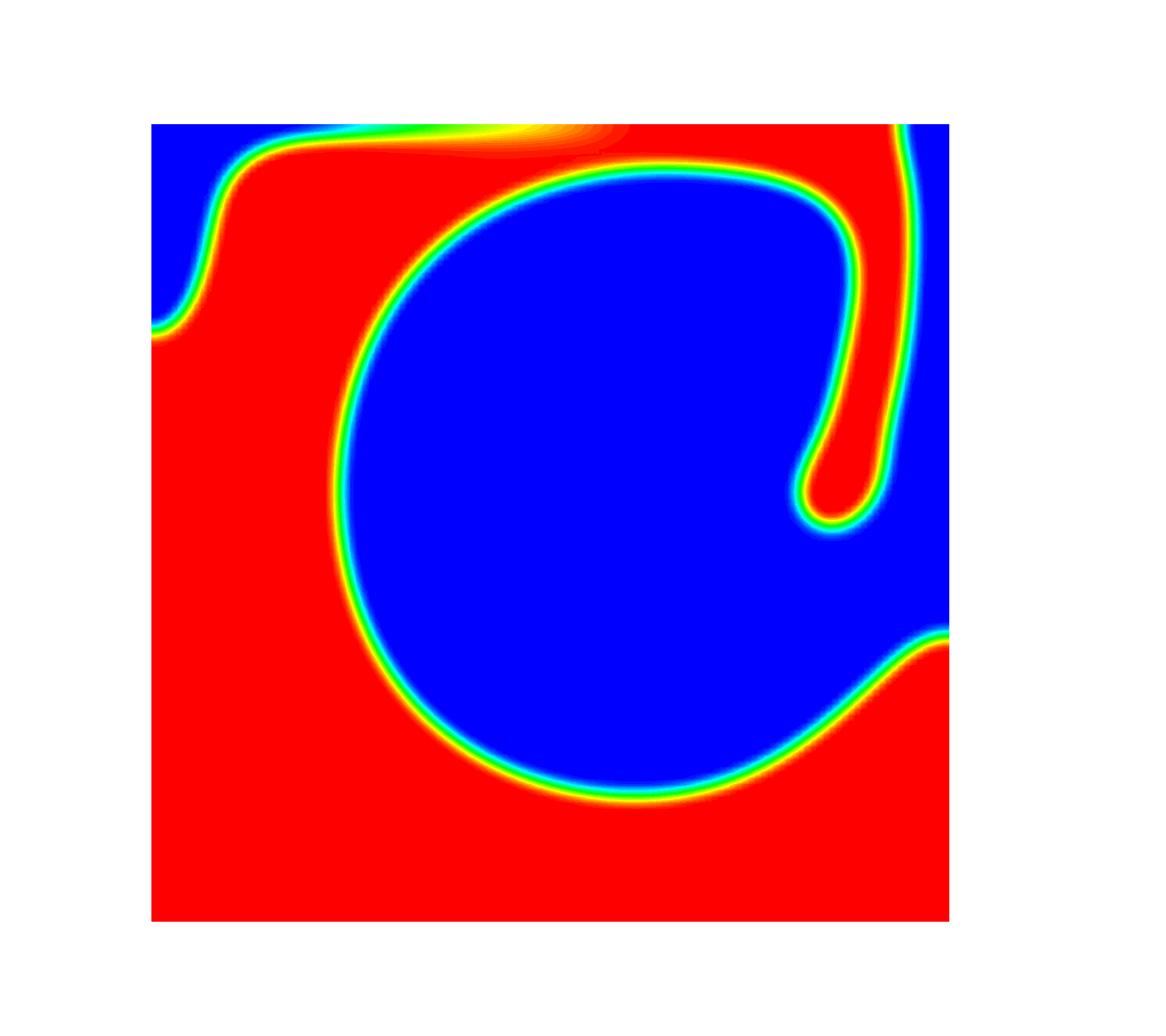}}\hskip -3mm 
		\subfigure[$t=12.0$] 	{\includegraphics[width=1.7in,trim=70 0 60 20,clip]{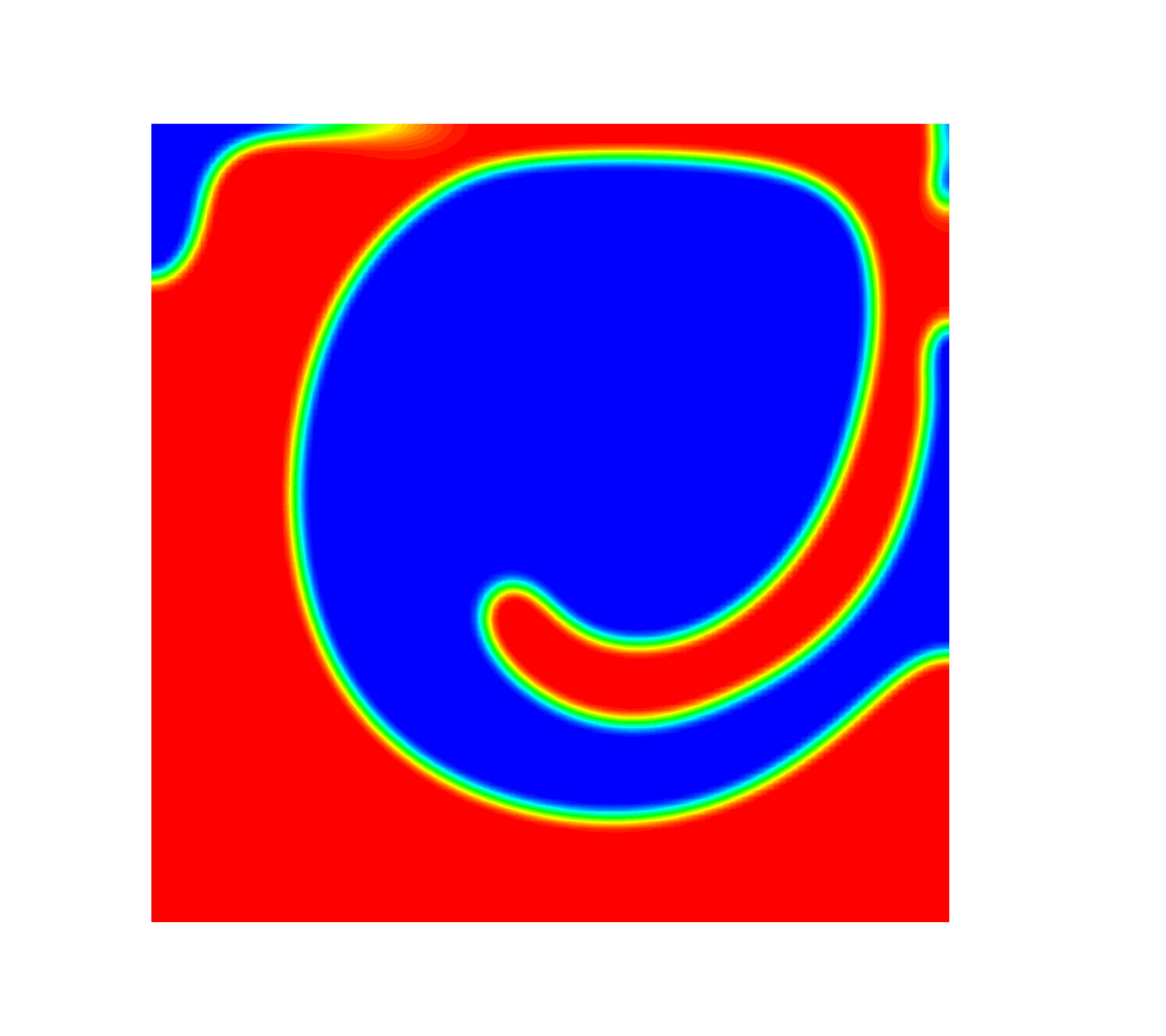}}\hskip -3mm 
		\subfigure[$t=15.0$] 	{\includegraphics[width=1.7in,trim=70 0 60 20,clip]{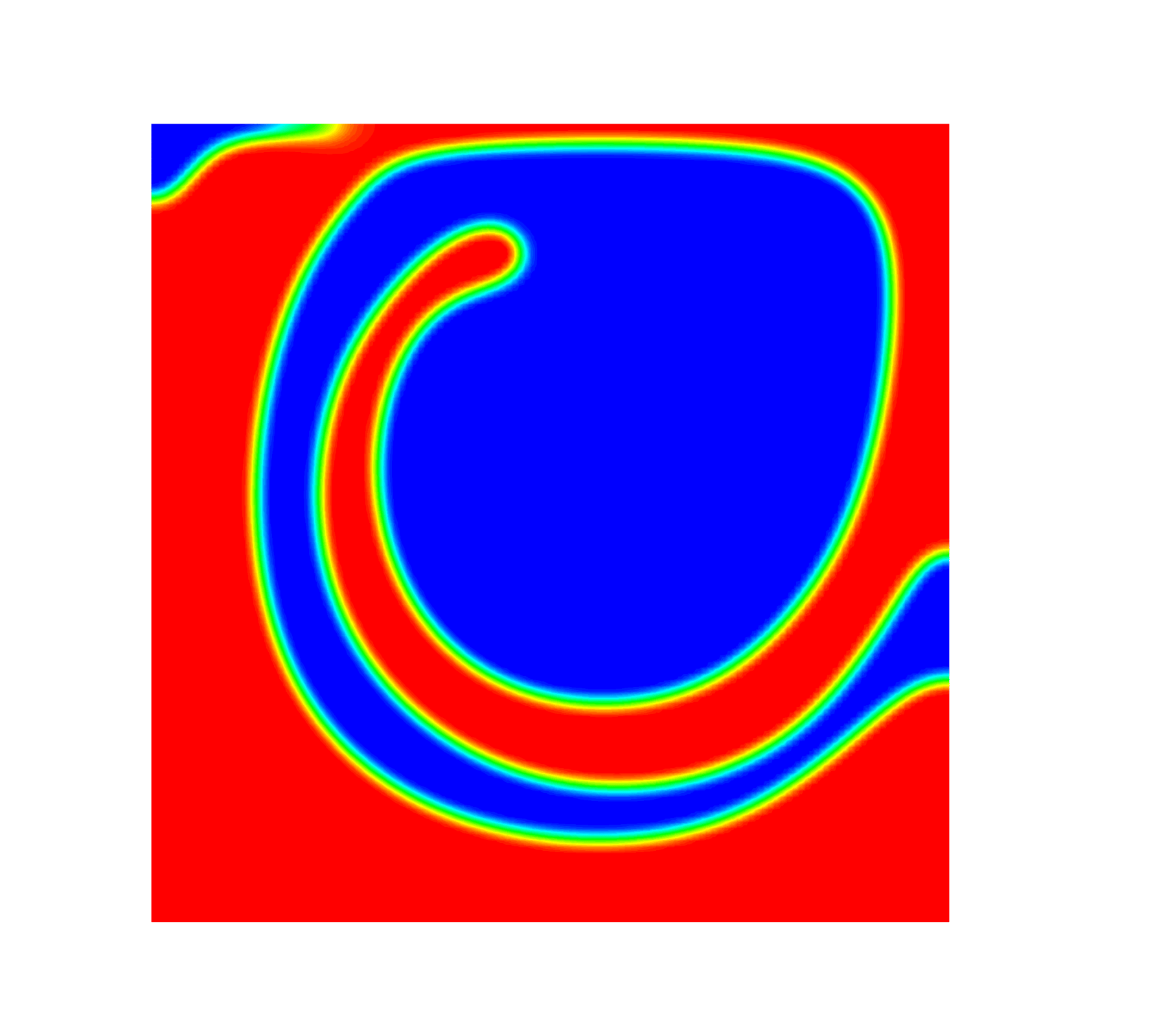}}
	\end{center}
	\caption{The evolution of phase variable for coarsening process.}
	\label{lid_driven} 
\end{figure}

Figure \ref{lid_driven} shows the evolution of the concentration $\phi$, in which red and blue indicate different   phases of  binary fluid. We can clearly observe that  a concave interface is formed as the fluids tend to mix, confirming the capability of our numerical method to cope with topological changes.  

\newpage
\subsection{Rayleigh-Taylor instability}
We simulate the Rayleigh-Taylor instability of binary fluid with different densities subject to buoyancy-driven flow. We consider two-phase flows with small density ratios by applying the Boussinesq approximation.

The momentum equation is rewritten as follows:
\begin{equation}\label{NS_Buo}
\rho_0\left(\frac{\partial \Au}{\partial t}-\frac{1}{Re} \nabla \cdot\big(\eta(\phi) \nabla \Au\big)\right)+\Au\cdot \nabla \Au+\nabla p+\frac{\epsilon^{-1}}{\sw} \phi \nabla \mu=-(1+\phi)g(\rho_1-\rho_0)-(1-\phi)g(\rho_2-\rho_0),
\end{equation} 
where $\rho_1$ and $\rho_2$ are the densities of the heavier fluid on top and lighter fluid, respectively, $\rho_0=\frac{\rho_1+\rho_2}{2}$ is considered to be the background density, and $\bm{g}=[0,g]^T$ is the gravitational acceleration. We simplify the right-hand side of \eqref{NS_Buo} to $-\bm{g}\phi(\rho_1-\rho_2)$, which is taken to be the external buoyancy term.

\begin{figure}[!ht]
    \centering
    \subfigure[$\eta=0.01$, $t=0.2,0.6,0.7,0.8,0.9,1.1$]{
        \includegraphics[height=2in,trim=1.0in 0.8in 4.2in 0.9in,clip]{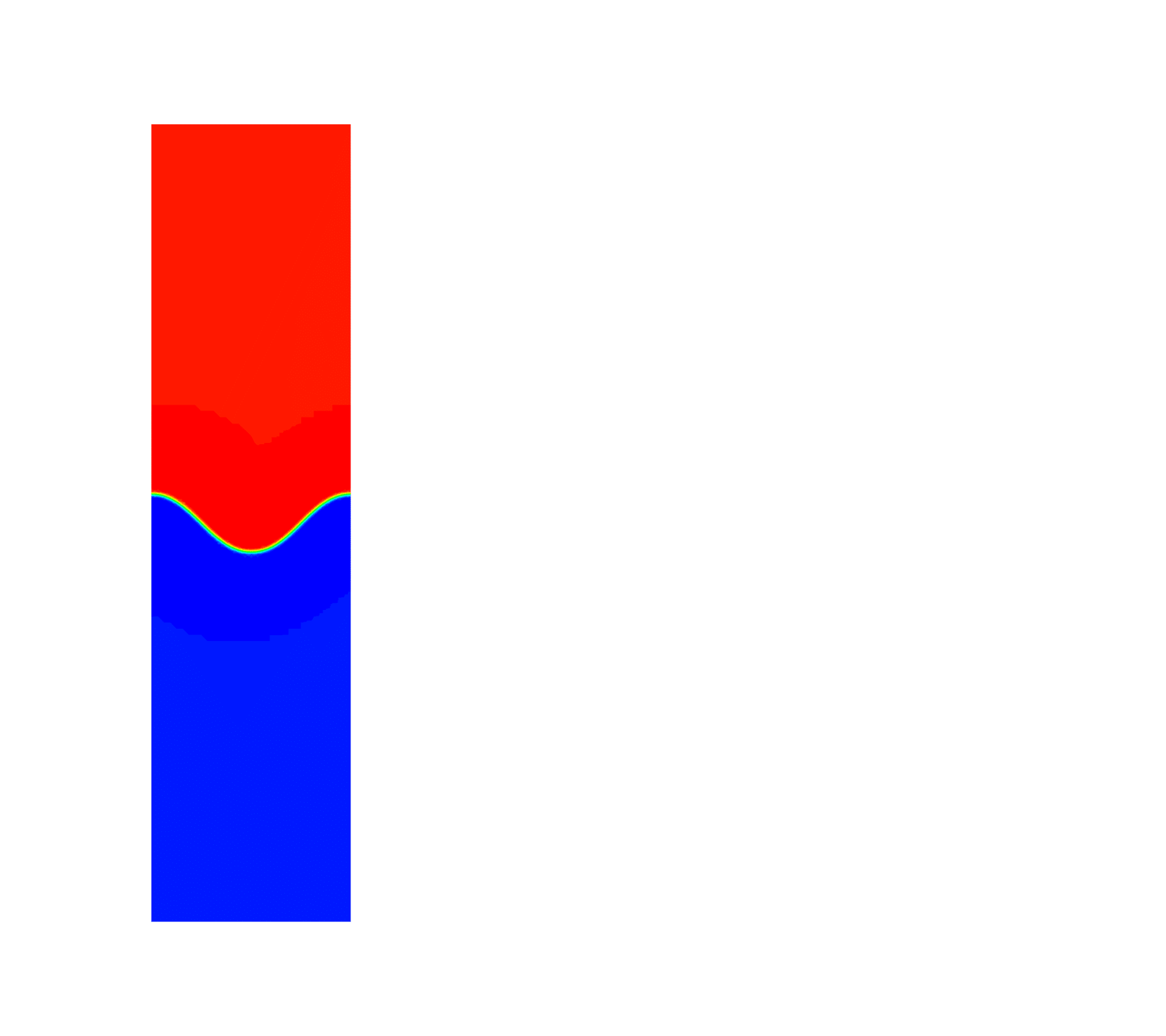}\hspace{-0.15in}
        \includegraphics[height=2in,trim=1.0in 0.8in 4.2in 0.9in,clip]{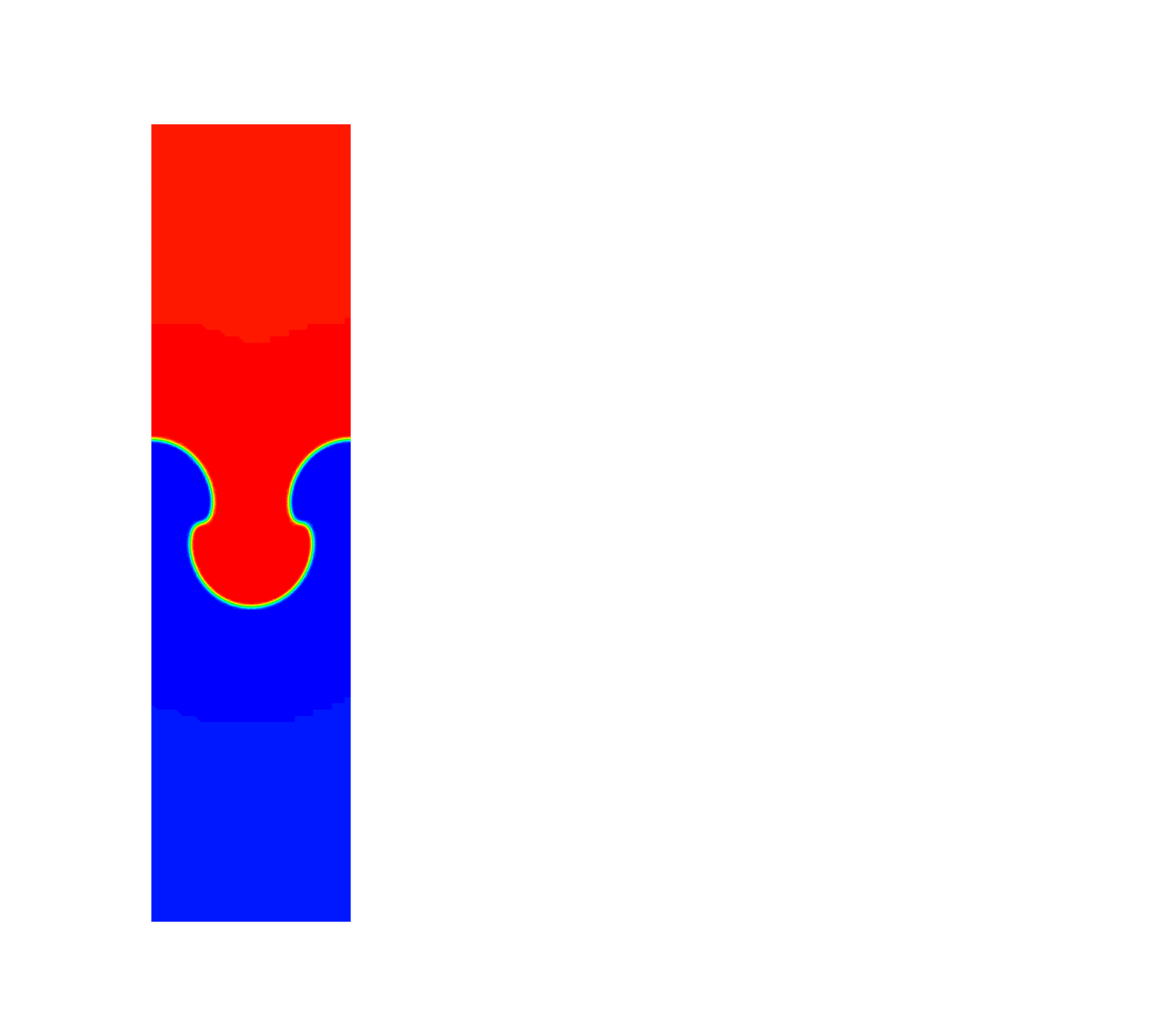}\hspace{-0.15in}
        \includegraphics[height=2in,trim=1.0in 0.8in 4.2in 0.9in,clip]{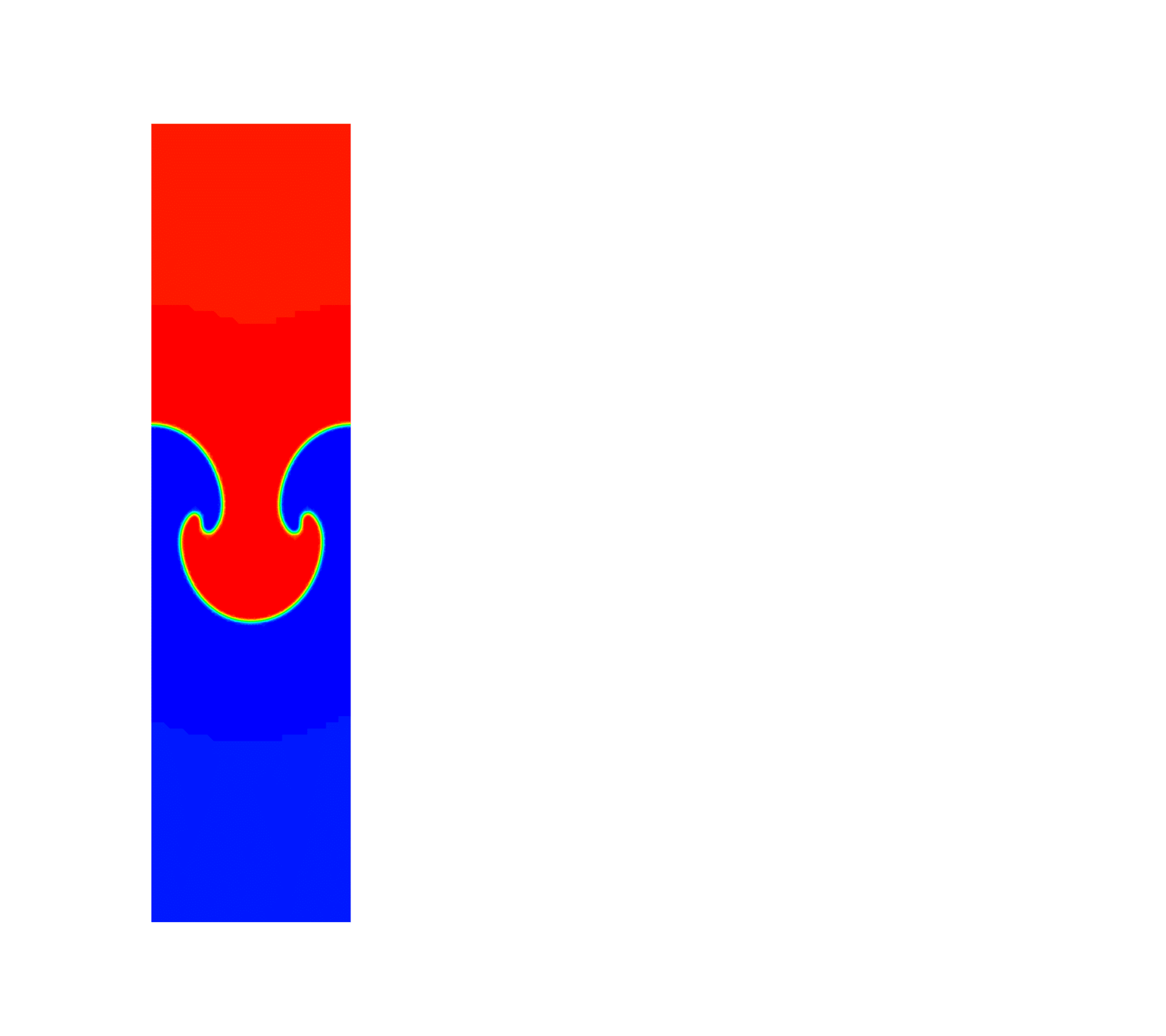}\hspace{-0.15in}
        \includegraphics[height=2in,trim=1.0in 0.8in 4.2in 0.9in,clip]{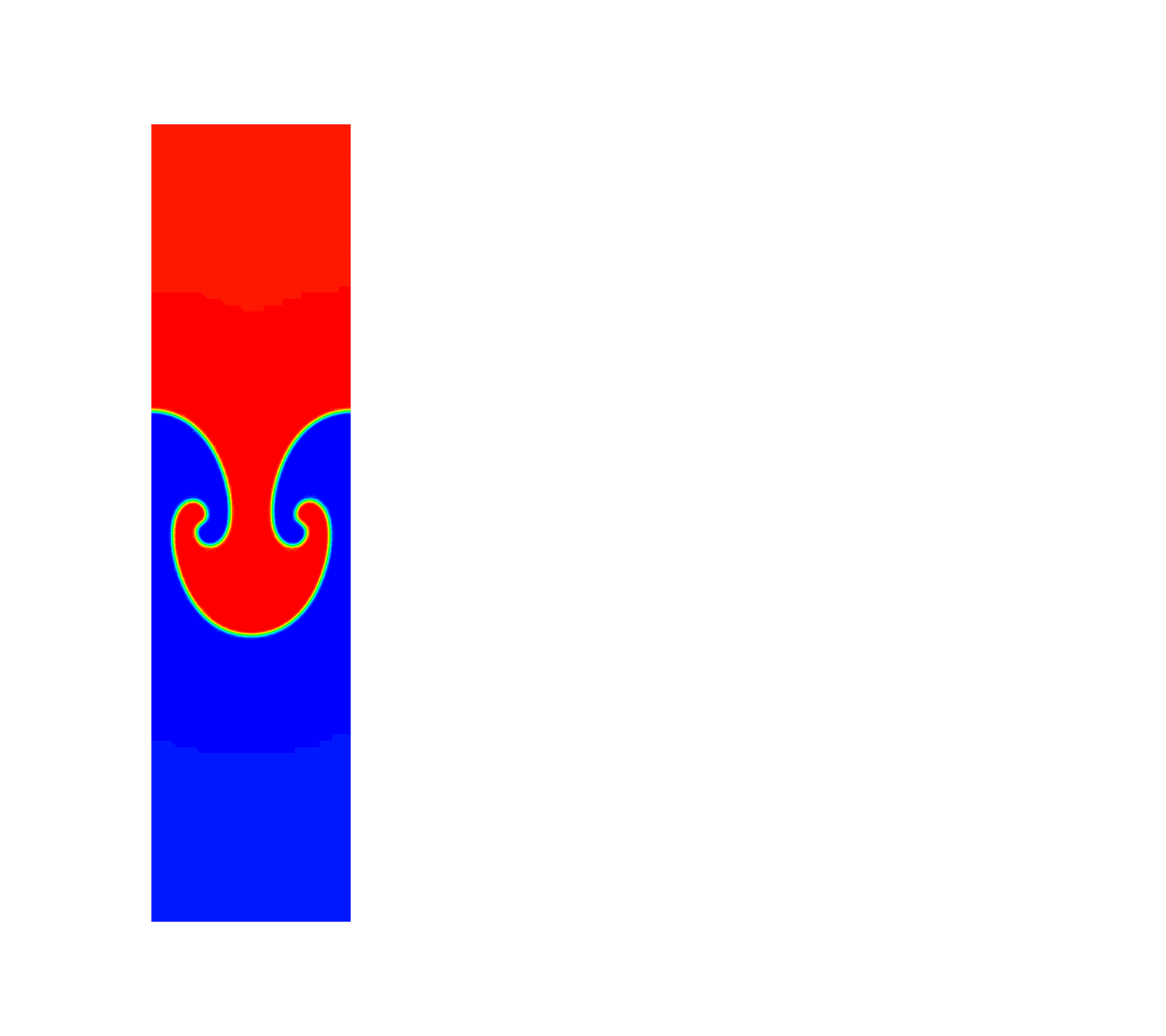}\hspace{-0.15in}
        \includegraphics[height=2in,trim=1.0in 0.8in 4.2in 0.9in,clip]{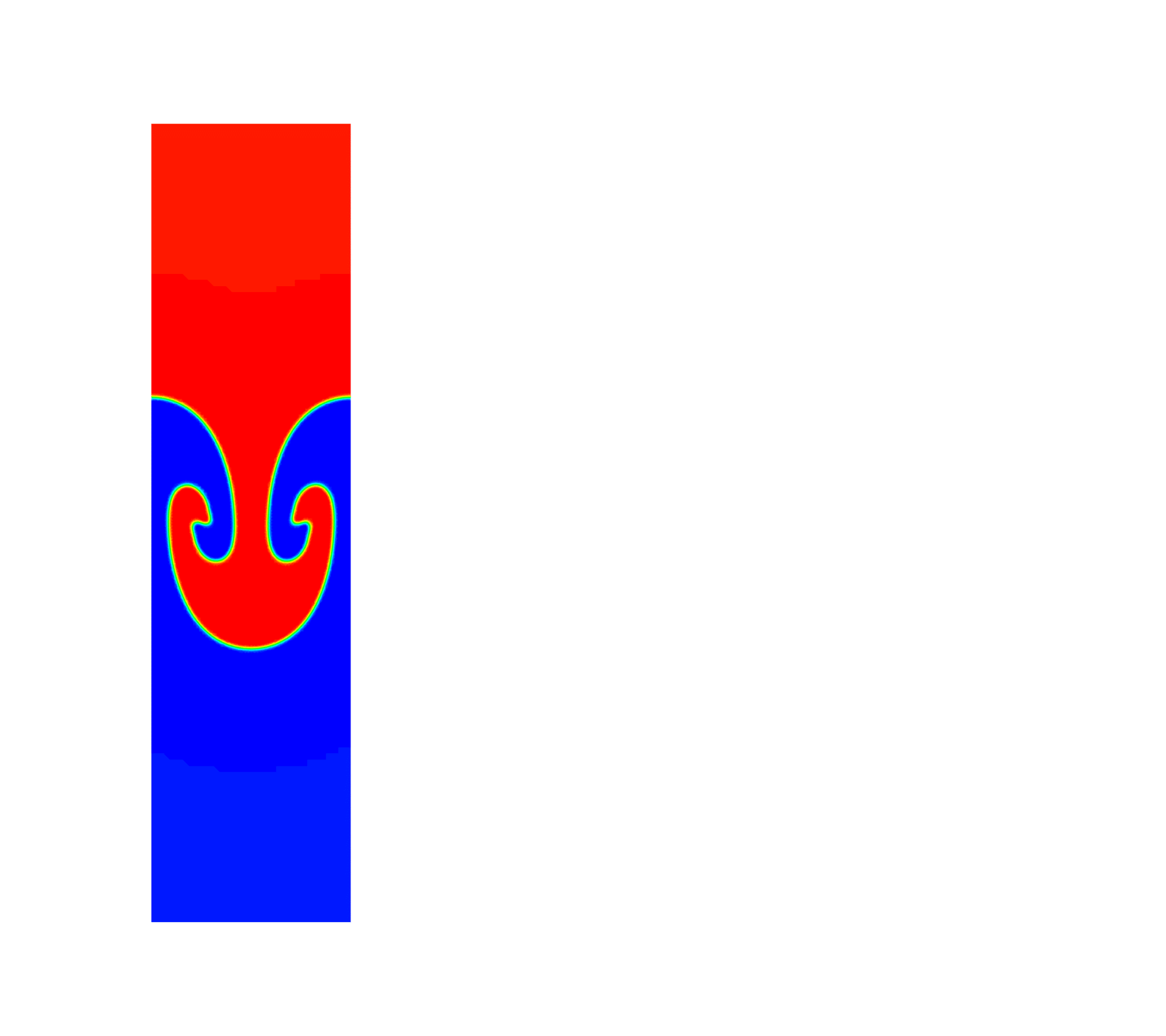}\hspace{-0.15in}
        \includegraphics[height=2in,trim=1.0in 0.8in 4.2in 0.9in,clip]{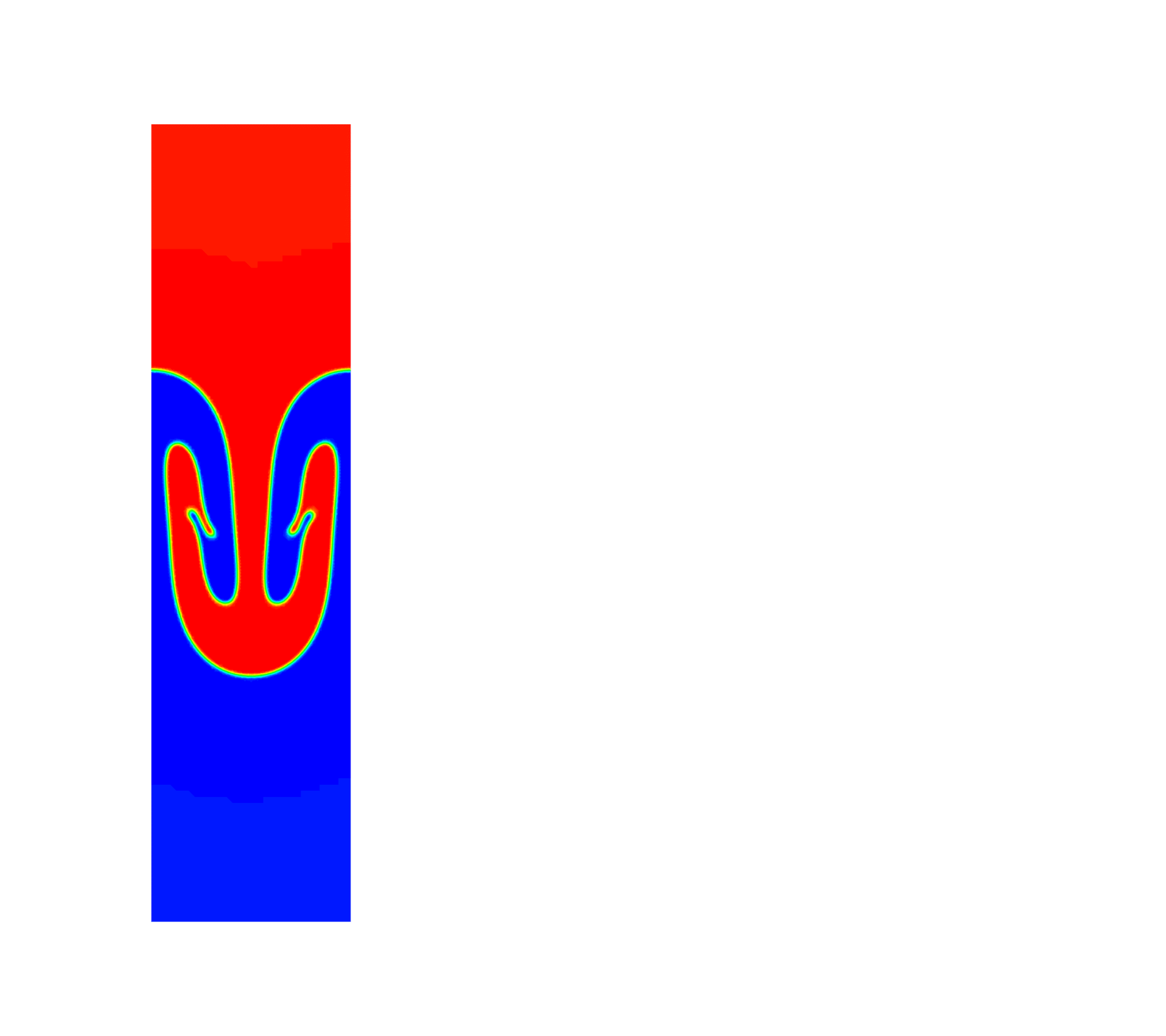}
        \label{nu10-2}}
    \subfigure[$\eta=0.001$, $t=0.3,0.6,0.7,0.8,0.9,1.2$]{
        \includegraphics[height=2in,trim=1.0in 0.8in 4.2in 0.9in,clip]{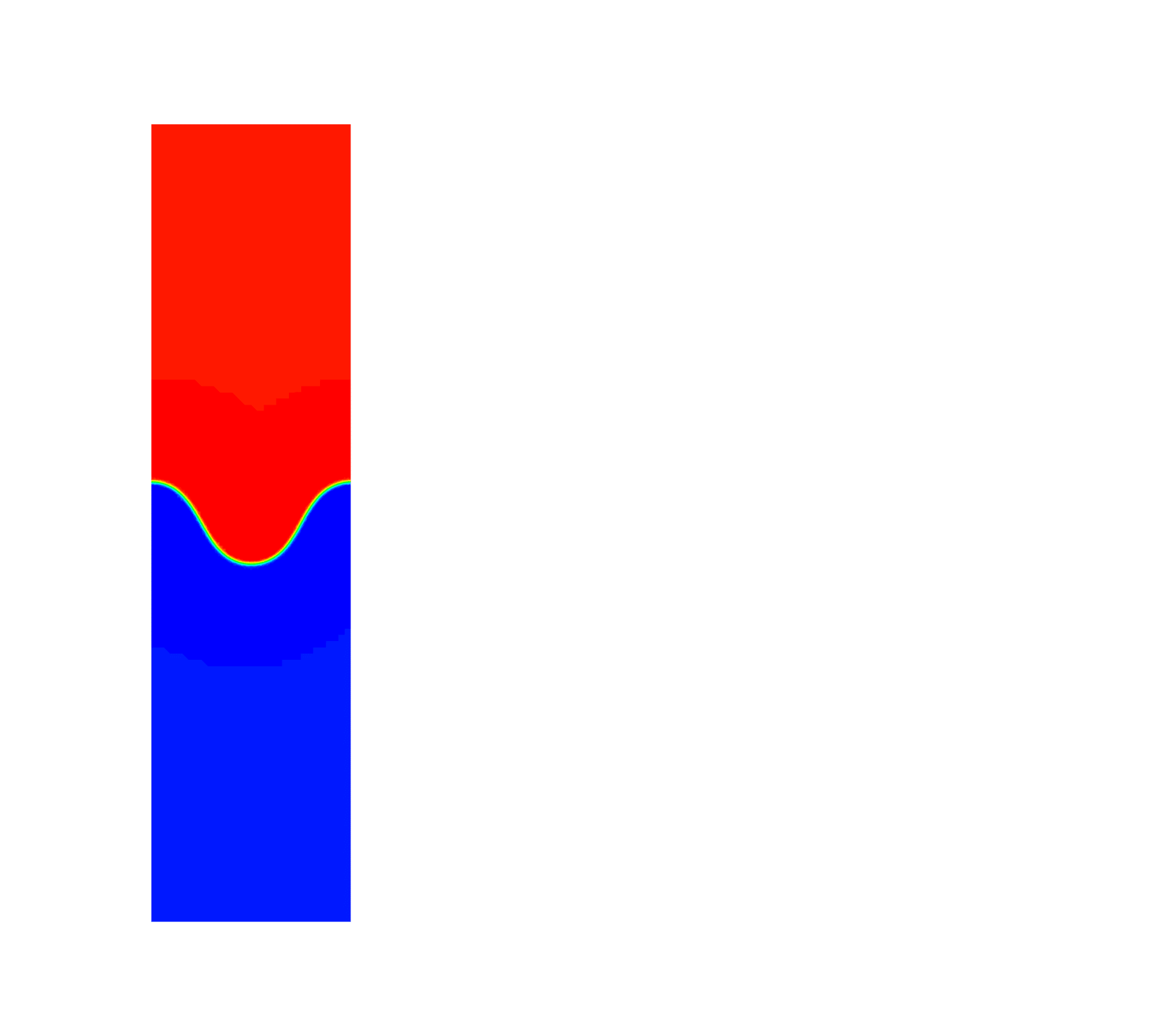}\hspace{-0.15in}
        \includegraphics[height=2in,trim=1.0in 0.8in 4.2in 0.9in,clip]{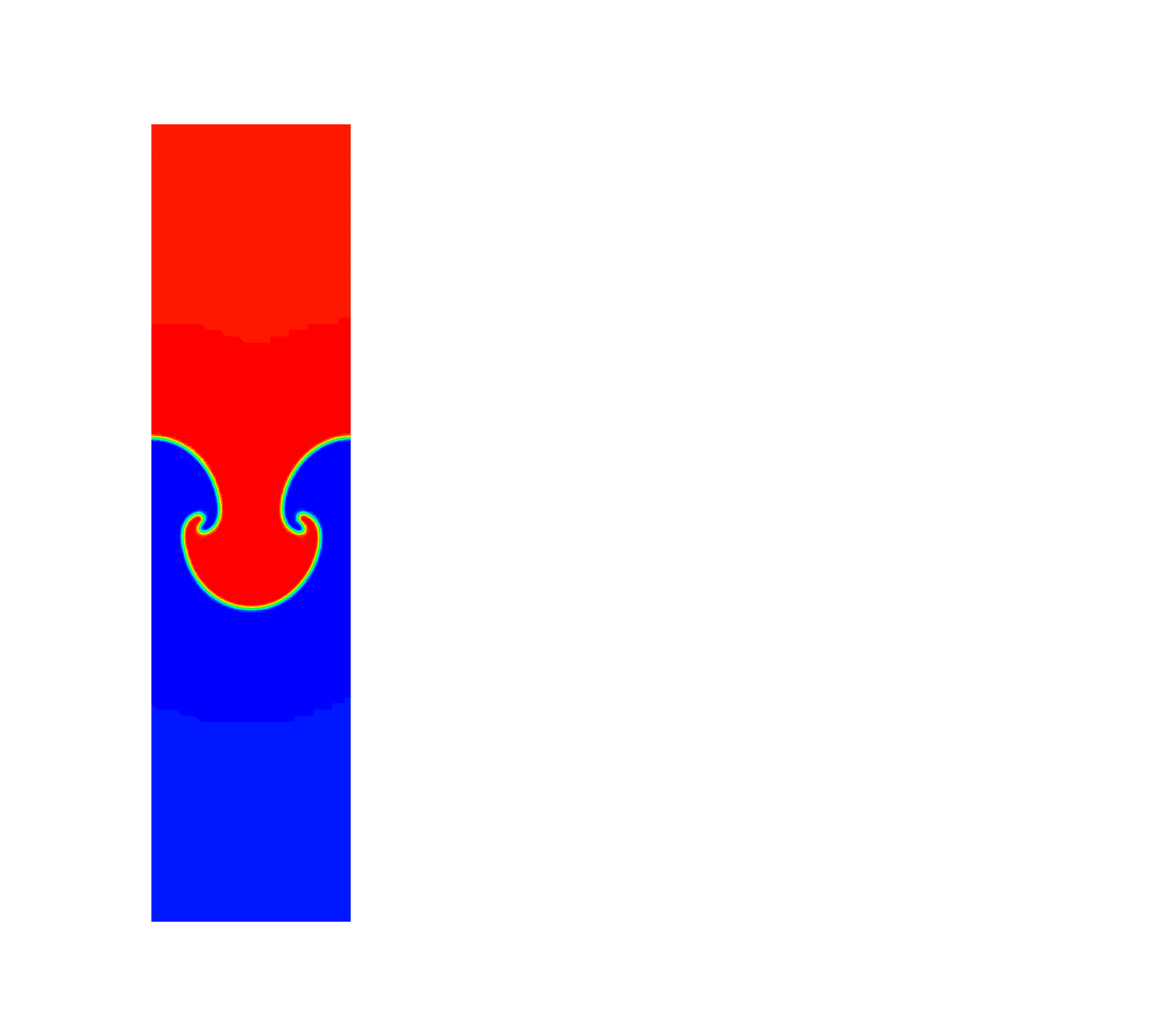}\hspace{-0.15in}
        \includegraphics[height=2in,trim=1.0in 0.8in 4.2in 0.9in,clip]{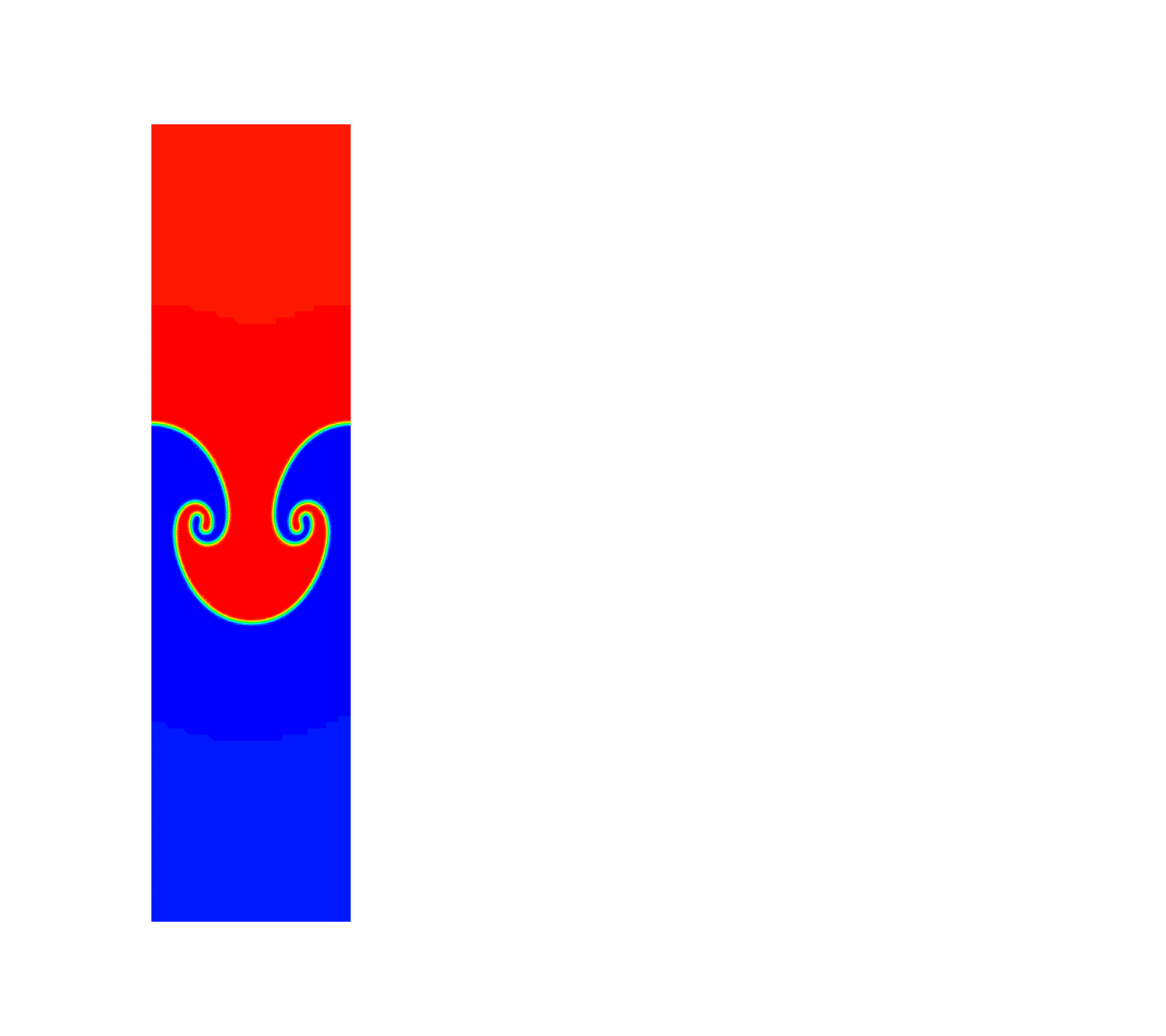}\hspace{-0.15in}
        \includegraphics[height=2in,trim=1.0in 0.8in 4.2in 0.9in,clip]{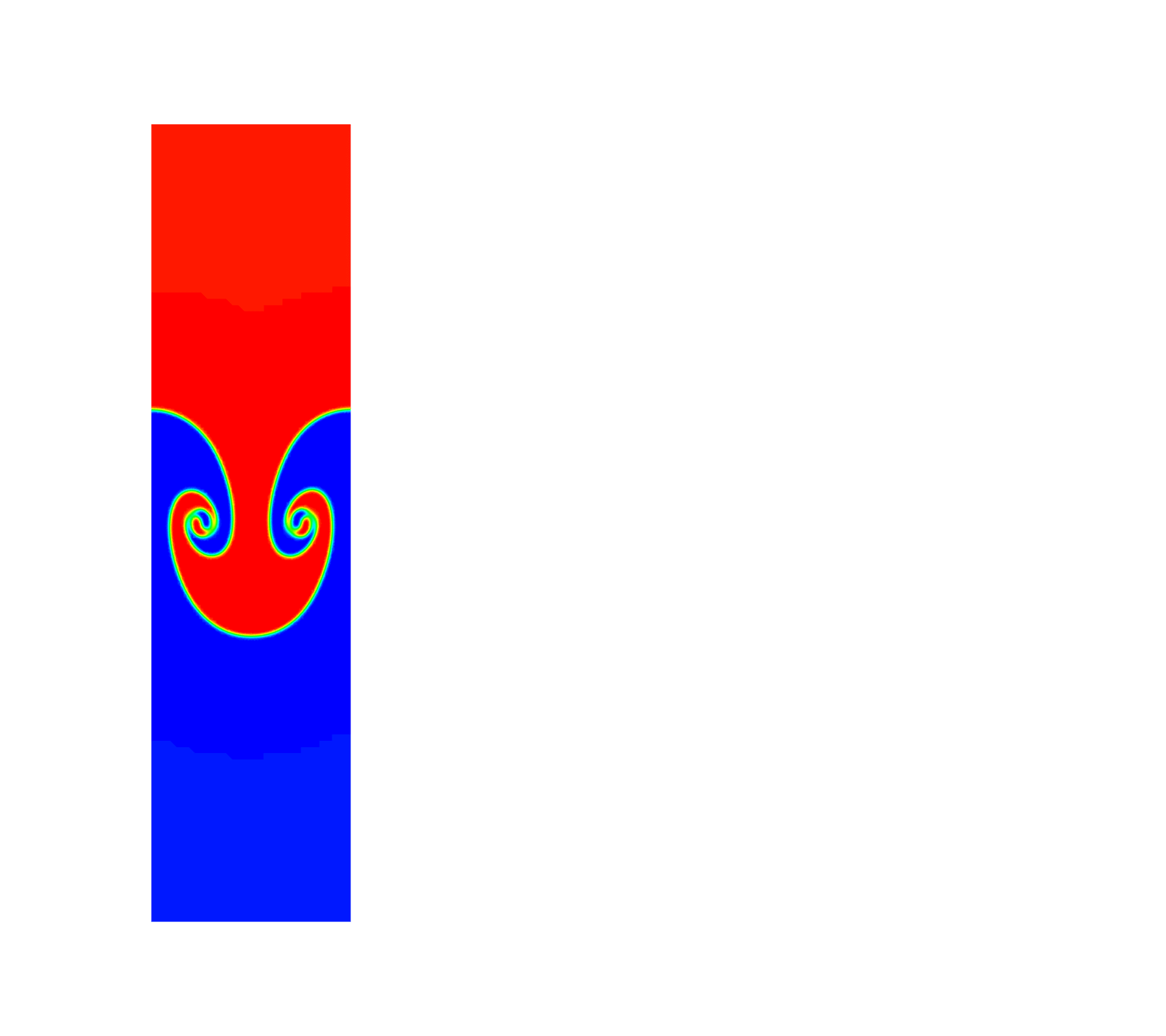}\hspace{-0.15in}
        \includegraphics[height=2in,trim=1.0in 0.8in 4.2in 0.9in,clip]{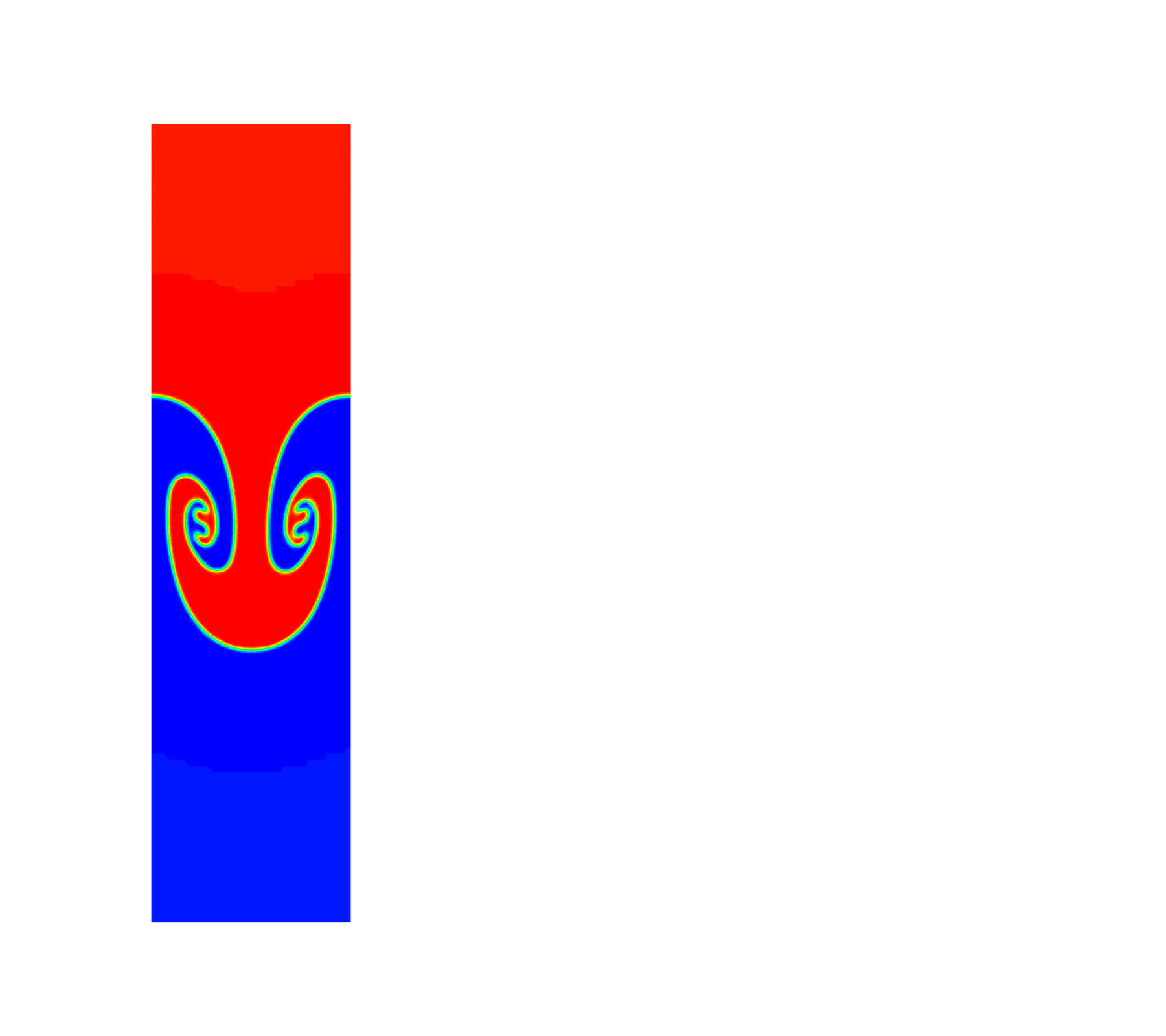}\hspace{-0.15in}
        \includegraphics[height=2in,trim=1.0in 0.8in 4.2in 0.9in,clip]{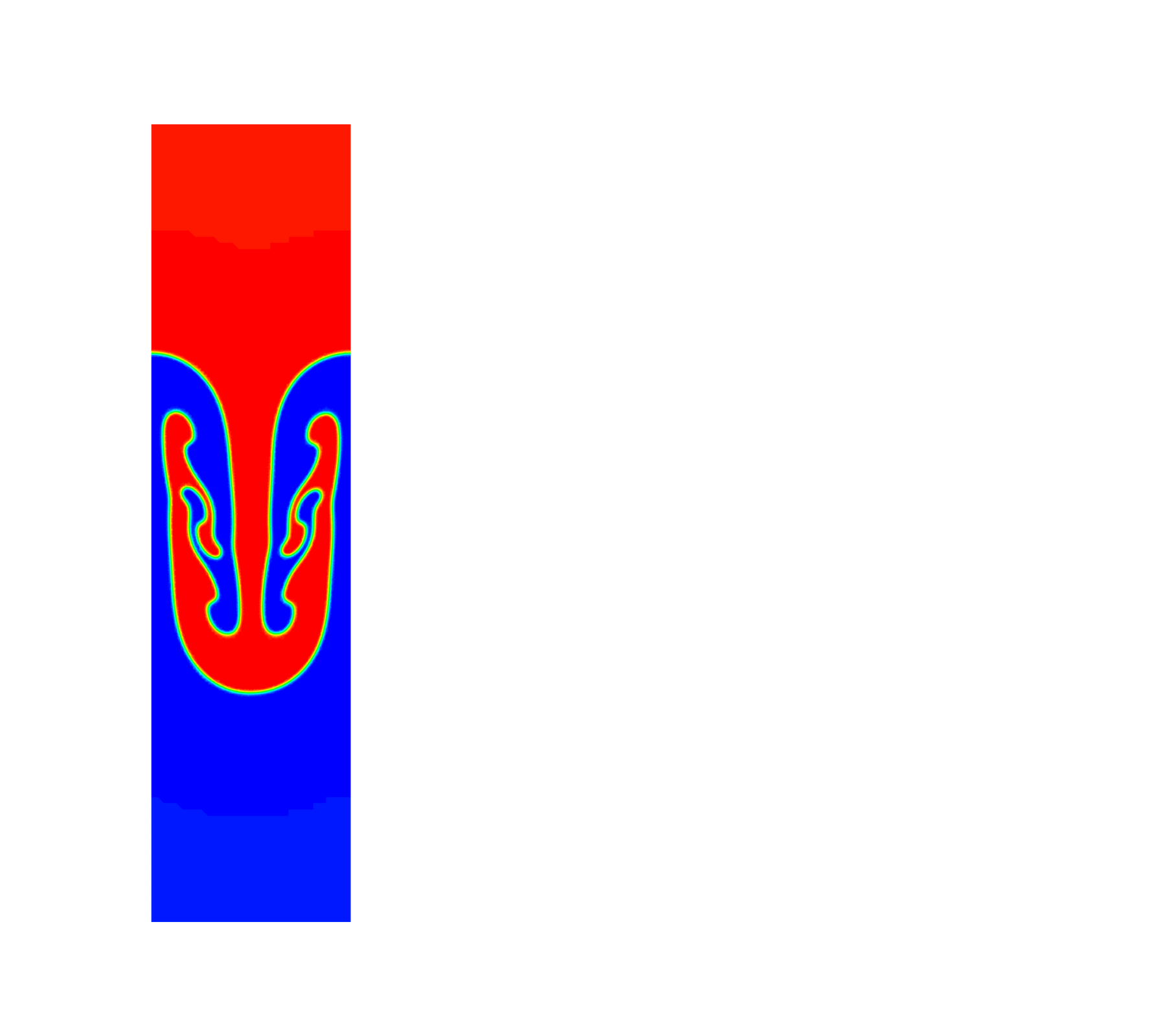}
        \label{nu10-3}
    }
    \caption{Snapshots of phase variable for Rayleigh-Taylor instability.}
    \label{Rayleigh_TaylorFigure}
\end{figure}

The initial position of the phase variable is set as $\phi_0=0.9\tanh\left(\dfrac{y-2-0.1\cos(2\pi x)}{\sqrt{2}\epsilon}\right)$ on domain $\Omega=[0,1]\times [0,4]$. We utilize the no-slip boundary conditions at the top and bottom boundaries and the free-slip boundary condition on the vertical walls for velocity, whereas the no-flux boundary conditions are enforced at all boundaries for the phase variable and chemical potential. Choose $\rho_1=3$, $\rho_2=1$, $g=10$, $Re=1$, $\epsilon=0.01$, $Pe=10$, $M=0.1$, and set the surface tension coefficient to zero, i.e., $\frac{\epsilon^{-1}}{\sw}=0$. The characterized snapshots of the phase variable are presented in Figure \ref{Rayleigh_TaylorFigure} corresponding to varied viscosities $\eta=0.01, 0.001$ at different times. 

Figure \ref{Rayleigh_TaylorFigure} clearly shows the anticipated prominent instability of the interface between immiscible fluids, specifically the heavier fluid on top penetrating the lighter fluid and eventually developing a spike falling downwards in both cases. Comparing the interface front in Figures \ref{nu10-2} and \ref{nu10-3}, the tail of the spike rolls up, forming a pair of counter-rotating vortices for smaller viscosity.

\newpage

\section{Conclusions}

In this paper, we propose a high order, bound-preserving and unconditionally energy stable finite element method for the Flory-Huggins Cahn-Hilliard-Navier-Stokes system. The $Q_k$ finite element with mass-lumping is employed for spatial discretization, while convex-concave splitting and pressure correction method are used for time-marching.   The unique solvability and unconditional stability are  rigorously established.  Moreover, by a key discrete $L^1$ estimate of the singular term we show that the solution is bounds-preserving.    The accuracy and robustness are verified by numerical simulations including  lid-driven cavity flow and Rayleigh-Taylor instability.  The method is also applicable to other phase filed fluid models with the Flory-Huggins potential.

\section*{Acknowledgments} 
	The work of  Y. Gao is partially supported by the NSFC, PR China under grants No. 12371406 and No. 11931013,  Guangdong Basic and Applied Basic Research Foundation under   grant No. 2023A1515010697.
D. Han and S. Sarkar are  supported by the National Science Foundation under  grant No. DMS-2310340. We thank Blanca Ayuso De Dios and Guosheng Fu for helpful discussions.

\def\cprime{$'$}

\end{document}